\documentclass[11pt]{amsart}
\usepackage{amssymb}
\usepackage{verbatim}
\usepackage[colorlinks]{hyperref}
\usepackage[colorinlistoftodos]{todonotes}
\usepackage[margin=1in]{geometry}

\usepackage{pgf}
\usepackage{tikz}
\usepackage{tikz-cd}
\usepackage{ytableau}

\newcommand{\AH}{\beta}
\newcommand{\am}{a_-}
\newcommand{\anti}{\ominus}
\newcommand{\ap}{a_+}
\newcommand{\Aut}{\mathrm{Aut}}
\newcommand{\bA}{\mathbb{A}}

\newcommand{\bc}{\overleftarrow{c}}
\newcommand{\bG}{{\overleftarrow{\G}}}
\newcommand{\bGring}{B}
\newcommand{\bR}{\overleftarrow{R}}
\newcommand{\bS}{\overleftarrow{S}}
\newcommand{\C}{\mathbb{C}}
\newcommand{\cO}{\mathcal{O}}

\newcommand{\diag}{\mathrm{diag}}

\newcommand{\Eb}{{\overline{E}}}

\newcommand{\Fl}{\mathrm{Fl}}
\newcommand{\Flb}{\overline{\Fl}}
\newcommand{\fmin}{f_{\mathrm{min}}}
\newcommand{\For}{\mathrm{For}}
\newcommand{\Frac}{\mathrm{Frac}}
\newcommand{\Fun}{\mathrm{Fun}}
\newcommand{\Gdiff}{\tilde \pi}
\newcommand{\G}{\mathfrak{G}}

\newcommand{\Ge}{\mathrm{\mathfrak{Ge}}}

\newcommand{\Gr}{\mathrm{Gr}}
\newcommand{\Grb}{\overline{\Gr}}
\newcommand{\Grd}{\Gr^\bullet}
\newcommand{\gtrlex}{>_{\mathrm{lex}}}

\newcommand{\hLa}{\hat{\La}}
\newcommand{\hs}{\hat{s}}
\newcommand{\Hom}{\mathrm{Hom}}
\newcommand{\id}{\mathrm{id}}
\newcommand{\K}{\mathbb{K}}
\newcommand{\Kt}{\tilde{\mathbb{K}}}
\newcommand{\la}{\lambda}
\newcommand{\La}{\Lambda}
\newcommand{\lad}{{\lambda^\bullet}}
\newcommand{\Lambdab}{{\overline{\Lambda}}}
\newcommand{\leleft}[1]{\,{}^{#1}\!\!\le}
\newcommand{\leright}[1]{\le^{#1}}
\newcommand{\loc}{\mathrm{loc}}
\newcommand{\move}{\mathrm{move}}

\newcommand{\omk}{\tilde \omega}
\renewcommand{\P}{{\mathbb Y}}
\newcommand{\pair}[2]{\langle #1\,,\,#2\rangle}
\newcommand{\parh}{\hat{\partial}}
\newcommand{\pih}{\hat{\pi}}
\newcommand{\pit}{\theta}
\newcommand{\Psit}{\tilde{\Psi}}
\newcommand{\psit}{\tilde{\psi}}

\newcommand{\pt}{\mathrm{pt}}
\newcommand{\Q}{\mathbb{Q}}
\newcommand{\res}{\mathrm{res}}

\newcommand{\RS}{\mathcal{RS}}
\newcommand{\Rt}{\tilde{R}}

\newcommand{\Spec}{{\rm Spec}}
\newcommand{\SSS}{{\mathbb{S}}}

\newcommand{\tLa}{\tilde{\Lambda}}
\newcommand{\tripair}[3]{\langle \,#1 | #2 | #3\,\rangle}
\newcommand{\ts}{\tilde{s}}
\newcommand{\vdim}{\mathrm{vdim}}

\newcommand{\wt}{\mathrm{wt}}
\newcommand{\xm}{x_-}
\newcommand{\xp}{x_+}
\newcommand{\Y}{\mathbb{Y}}
\newcommand{\Z}{\mathbb{Z}}

\newtheorem{prop}{Proposition}
\newtheorem{thm}[prop]{Theorem}
\newtheorem{lem}[prop]{Lemma}
\newtheorem{cor}[prop]{Corollary}
\newtheorem{conjecture}[prop]{Conjecture}
\theoremstyle{remark}
\newtheorem{ex}[prop]{Example}
\newtheorem{rem}[prop]{Remark}
\numberwithin{prop}{section}

\pgfmathsetmacro{\boxsize}{1}
\pgfmathsetmacro{\halfboxsize}{0.5}
\newcommand{\bbox}[2]{
\draw[thin] (#1,#2)-- (#1,#2+\boxsize)-- (#1+\boxsize,#2+\boxsize)-- (#1+\boxsize,#2)-- (#1,#2);
}

\newcommand{\rightelbow}[3]{
\draw[#3] (#1+\boxsize,#2) ++ (90:\halfboxsize) arc (90:180:\halfboxsize);
}
\newcommand{\leftelbow}[3]{
\draw[#3] (#1,#2+\boxsize) ++ (0:\halfboxsize) arc (0:-90:\halfboxsize);
}
\newcommand{\cross}[4]{
\draw[#3] (#1+\halfboxsize,#2) --  (#1+\halfboxsize,#2+\boxsize);
\draw[#4] (#1,#2+\halfboxsize) --  (#1+\boxsize,#2+\halfboxsize);
}
\newcommand{\horline}[3]{
\draw[#3] (#1,#2+\halfboxsize) --  (#1+\boxsize,#2+\halfboxsize);
}
\newcommand{\vertline}[3]{
\draw[#3] (#1+\halfboxsize,#2) --  (#1+\halfboxsize,#2+\boxsize);

}

\def\half{blue}

\title{Back stable $K$-theory Schubert Calculus}
\author{Thomas Lam}
\address{Department of Mathematics \\
University of Michigan \\
530 Church St. \\
Ann Arbor 48109 USA}
\email{tfylam@umich.edu}
\thanks{T.L. was supported by NSF DMS-1464693 and NSF DMS-1953852}
\author{Seung Jin Lee}
\address{Department of Mathematical Sciences \\ Research institute of Mathematics \\ Seoul National University \\ Gwanak-ro 1, Gwanak-gu \\ Seoul 151-747 Republic of Korea}
\email{lsjin@snu.ac.kr}
\author{Mark Shimozono}
\address{Department of Mathematics \\
460 McBryde Hall, Virginia Tech\\
 255 Stanger St. \\
Blacksburg, VA, 24601, USA }
\email{mshimo@math.vt.edu}
\begin{document}
\begin{abstract}
We study the back stable $K$-theory Schubert calculus of the infinite flag variety.  We define back stable (double) Grothendieck polynomials and double $K$-Stanley functions and establish coproduct expansion formulae.  Applying work of Weigandt, we extend our previous results on bumpless pipedreams from cohomology to $K$-theory.  We study finiteness and positivity properties of the ring of back stable Grothendieck polynomials, and divided difference operators in $K$-homology.
\end{abstract}
\maketitle
\setcounter{tocdepth}{1}
\numberwithin{equation}{section}
\tableofcontents

\section{Introduction}
In \cite{LLS:back stable} we initiated the study of back stable Schubert calculus.  In the present work, we generalize from cohomology to K-theory, and study back stable $K$-theory Schubert calculus.  We assume the reader has some familiarity with our earlier work, and in the introduction we will emphasize some of the differences.

\subsection{Back stable Grothendieck polynomials}

Whereas the back stable Schubert polynomials of \cite{LLS:back stable} are limits of Schubert polynomials, in this paper we will be concerned with the back stable Grothendieck polynomials
$$
\bG_w(x) := \lim_{p \to -\infty, q \to \infty} \G_w(x_p,x_{p+1},\ldots,x_q)
$$
which are limits of the Grothendieck polynomials $\G_w(x)$ of Lascoux and Sch\"utzenberger \cite{LaSc}.  Grothendieck polynomials represent classes of Schubert structure sheaves in the $K$-theory of the flag variety. 

The same definition of back stable Grothendieck polynomials has appeared in the work of Marberg and Pawlowski \cite{MP}, who studied the principal specializations ($x_i \mapsto q^{i-1})$ of these power series.  We remark that back stable Schubert polynomials have been generalized to the involution setting by Pawlowski \cite{Paw}.  

\subsection{Infinite flag varieties}
Infinite-dimensional flag varieties, such as flag varieties of Kac-Moody groups, come in a number of algebro-geometric variants.  Of interest to us is a thin flag variety $\Fl$ that is an ind-finite variety and a thick flag variety $\Flb$ that is an infinite-dimensional scheme.  There are other versions such as semi-infinite flag varieties and geometric models based on loop groups which will not feature in this work.  The $K$-groups in this work are certain Grothendieck groups of coherent sheaves, and the choice of scheme structure plays a more significant role than in our earlier work in cohomology.  For a discussion of $K$-groups on thin and thick flag varieties in the Kac-Moody setting, we refer the reader to Kumar \cite{Kum} and Baldwin and Kumar \cite{BK}.

Whereas in \cite{LLS:back stable} we only considered the thin infinite flag variety, in the present work we also consider the thick infinite flag variety.  Thick infinite flag varieties were studied by Kashiwara \cite{Ka} and we give a mostly elementary treatment in \S\ref{sec:geometry}.   A different and elegant approach to infinite flag varieties is also given by Anderson \cite{And:infinite}.  We show in Theorem~\ref{T:infinite basis} and Corollary~\ref{C:Schub basis} that back stable double Grothendieck polynomials represent classes of Schubert structure sheaves in the equivariant $K$-group $K_T(\Flb)$ of the thick infinite flag variety $\Flb$.  

\subsection{Coproduct formula}
Like their cohomological counterparts, back stable Grothendieck polynomials (and their double versions) satisfy a coproduct formula (Theorem~\ref{T:coproduct}), decomposing them into a symmetric part and a finite part:
\begin{equation}\label{E:coprodintro}
\bG_w = \sum_{\substack{u*v=w \\ v\in S_{\ne0}}} (-1)^{\ell(u)+\ell(v)-\ell(w)}
  G_u \G_v.
\end{equation}
Here, $u*v$ denotes the Demazure or Hecke product, $G_u$ denotes the $K$-Stanley symmetric function~\cite{Bu,FK}\footnote{The symmetric function $G_w$ is usually called a \emph{stable Grothendieck polynomial}, but to avoid conflicts in terminology we use a different name.}, and $\G_v$ is the Lascoux-Sch\"utzenberger Grothendieck polynomial.  We deduce~\eqref{E:coprodintro} from the coproduct formula in $K$-theory \emph{affine} Schubert calculus~\cite{LLS:coprod}.

\subsection{Double $K$-Stanley symmetric functions}
Back stable double Grothendieck polynomials $\bG_w(x;a)$ are defined as a limit of double Grothendieck polynomials $\G_w(x;a)$ in a similar manner.  We define the \emph{double $K$-Stanley symmetric functions} as the image $G_w(x||a):=\eta_a(\bG_w(x;a))$ under an algebra homomorphism $\eta_a: \bR(x;a)\to \hLa(x||a)$ (see \S\ref{sec:doubleK}) from back stable power series to symmetric power series.  The symmetric functions $G_w(x||a)$ appear to be novel.  (Setting $a = 0$, we recover the $K$-Stanley symmetric function $G_w$.)
We obtain as a subset the Grassmannian double $K$-Stanley symmetric functions $G_\la(x||a)$.  These functions form a basis of a ring $\Gamma(x||a)$ (Theorem~\ref{T:Gammaxa}), a double analogue of Buch's $\Gamma$ ring \cite{Bu}.  We explain determinantal formulae for $G_\la(x||a)$ in \S\ref{S:detformula}, connecting our construction to the literature \cite{And,HIMN}.

We conjecture (Conjecture~\ref{conj:kpos})
that the expansion coefficients of $G_w(x||a)$ in the basis $\{G_\la(x||a)\}$ have alternating signs.  The analogous expansion coefficients of double Stanley functions $F_w(x||a)$ into double Schur functions $\{s_\la(x||a)\}$ were shown to be positive in \cite[Theorem 4.22]{LLS:back stable} using the quantum equals affine phenomenon.

\subsection{$K$-bumpless pipedreams}
In \cite{LLS:back stable}, we defined \emph{bumpless pipedreams} to give explicit monomial expansions for back stable double Schubert polynomials.  Weigandt \cite{Wei} connected bumpless pipedreams to earlier \emph{alternating sign matrix} formulae of Lascoux \cite{Las}, and thereby obtained formulae for double Grothendieck polynomials in terms of $K$-bumpless pipedreams.  Weigandt's work immediately gives a formula for back stable double Grothendieck polynomials in terms of $K$-bumpless pipedreams, which we state in Theorem~\ref{thm:Kbumpless}. We also use $K$-bumpless pipedreams to give formulae for the double $K$-Stanley symmetric functions (Theorem~\ref{thm:Khalfplane} and Corollary~\ref{C:double G lambda by tableaux}), and expansion formulae (Theorem~\ref{thm:Krect} and Corollary~\ref{cor:Gw}) for $\bG_w(x;a)$ and $G_w(x||a)$ in terms of double $K$-Stanley symmetric functions.

Recently, bumpless pipedreams have found applications in the study of diagonal Gr\"obner degenerations of matrix Schubert varieties \cite{HPW,Kle,KW}, in the study of products of Schubert polynomials \cite{Hua1,Hua2}, and in other applications \cite{FGS,BS,Xio}.

\subsection{$K$-homology}
In \S\ref{sec:homology}, we study a basis dual to the Grassmannian double $K$-Stanley symmetric functions $G_\la(x||a)$.  We call these symmetric functions $g_\la(y||a)$ the \emph{$K$-Molev functions}. These symmetric functions are $K$-theory analogues of Molev's dual Schur functions $\hs_\la(y||a)$ \cite{Mol}.  At $a = 0$, the symmetric functions $g_\la(y||a)$ reduce to the dual stable Grothendieck polynomials $g_\la(y)$ studied by Lam and Pylyavskyy \cite{LP}; see also \cite{Len}.  Geometrically, the functions $g_\la(y||a)$ form a basis of the equivariant $K$-group $K^T(\Gr)$ of the thin infinite Grassmannian.  We show in Theorem~\ref{T:create Khom} that the $g_\la(y||a)$ can be obtained recursively by applying $K$-homology divided difference operators.  In Theorem~\ref{T:KL}, we sketch the relation between our $g_\la(y||a)$ and the deformation of symmetric functions studied by Knutson and Lederer \cite{KL}.

\subsection{The algebra of back stable Grothendieck polynomials}
We define the \emph{algebra of back stable Grothendieck polynomials} 
$$\bGring:= \bigoplus_w \Z \bG_w$$
to be the span of all back stable Grothendieck polynomials.  We show in Theorem~\ref{T:bG basis} that $\bGring$ is a ring, or equivalently, the product of back stable Grothendieck polynomials is finite.  This finiteness is quite nontrivial; for example, it fails to hold for the $K$-theory of the affine flag variety of ${\rm SL}_2$.  The ring $\bGring$ is an infinite flag variety version of Buch's $\Gamma$ ring spanned by $K$-Stanley functions \cite{Bu}, an analogous ring for the infinite Grassmannian.

We show in Theorem~\ref{T:Omega span} that after adjoining an element $\Omega = \Omega[x_-]$, the ring $\bGring$ can be decomposed as a tensor product $\Gamma[\Omega] \otimes R^+$, where $R^+$ is spanned by finite Grothendieck polynomials.

We conjecture that the similar finiteness also holds for back stable double Grothendieck polynomials.  Curiously, we show in that Proposition~\ref{P:bG double basis} that this would hold if the positivity conjecture (Conjecture~\ref{conj:kpos}) were true.

\subsection{Flagged Grothendieck polynomials}
In the course of studying the relationship between Grothendieck polynomials and their symmetrizations given by the $K$-Stanley functions, we found it natural to study a family of polynomials which interpolate between them. 
These are the flagged Grothendieck polynomials, which we introduce in \S \ref{SS:flagged Groth} 
for an arbitrary permutation using divided differences. We prove in Proposition~\ref{P:flagged G} a monomial expansion for flagged Grothendieck polynomials that generalizes the Fomin-Kirillov formula for Grothendieck polynomials.
For the special cases of vexillary and 321-avoiding permutations, the flagged Grothendieck polynomials were defined combinatorially in \cite{KMY:tableau complexes} \cite{KMY:Grobner} \cite{Ma}.

\subsection{K-classes of degeneracy loci}
Fulton \cite{Fu} realized the double Schubert polynomial as a universal formula for the cohomology class of a degeneracy locus $\Omega_w$ defined by rank conditions on composite maps in a sequence of maps of vector bundles living on a common base.   
Recently, Anderson and Fulton \cite{AF} studied what might be called the ``back stable limit" of this degeneracy locus construction, obtaining \emph{enriched Schubert polynomials}, which specialize (and are nearly equivalent) to the back stable Schubert polynomials of \cite{LLS:back stable}. 

Buch \cite{Bu:K quiver} observed that the $K$-class of the structure sheaf of $\Omega_w$ in a flag bundle, has a \emph{universal} formula given by the double Grothendieck polynomial $\G_w(x;a)$, based on the work of Fulton and Lascoux \cite{FL}, who showed that after a certain change of variable, the double Grothendieck polynomials map to classes of structure sheaves of opposite Schubert varieties in the flag variety.
Very recently, Buch \cite{Buch:email} gave us a detailed explanation of a limit of the degeneracy locus construction of double Grothendieck polynomials and suggested the result should be the back stable double Grothendieck polynomials. He also suggested to apply the degeneracy loci formulae in 
\cite{BKTY} and \cite{BKSTY}. The result coincides with one of our formulae for back stable Grothendieck polynomials, and is explained in \S \ref{S:degeneracy loci}.

\subsection{Further directions/relations} We were unable to pursue many obvious avenues of investigation, for example, the study of the ideal sheaf basis, and the relation to $K$-theory affine Schubert calculus \cite{LSS, Mor}.  We briefly discuss these ideas in \S\ref{S:further}.

\medskip
\noindent
{\bf Acknowledgements.}
We thank Anders Buch for a suggestion that led to Section~\ref{S:degeneracy loci}, and Anna Weigandt for explaining the relation between bumpless pipedreams and the work of Lascoux.

\section{Thick infinite flag scheme}\label{sec:geometry}

\subsection{Dynkin type $A_\Z$}\label{SS:Dynkin AZ}
The Dynkin diagram of type $A_\Z$ has Dynkin node set $I=\Z$ and simple bonds $(i,i+1)$ for all $i\in\Z$. It has weight lattice $X' = \bigoplus_{i\in I} \Z \La_i$ with basis of fundamental weights $\La_i$. For $i\in \Z$, the simple coroot $\alpha_i^\vee\in \Hom_\Z(X',\Z)$ is defined by $\pair{\alpha_i^\vee}{\La_j} = \delta_{ij}$ for $i,j\in \Z$.
For $i\in \Z$, let $\epsilon_i := \La_i-\La_{i-1}$. For $i\in \Z$, let $\alpha_i := \epsilon_i-\epsilon_{i+1}$ be the simple root. Let $Q := \bigoplus_{i\in I} \Z \alpha_i$ and $X := \bigoplus_{i\in \Z} \Z \epsilon_i$. We have 
\begin{align}\label{E:lattices}
	Q \subset X \subset X'.
\end{align}
Let the Weyl group $S_\Z$ be the subgroup of $\Aut(X')$ generated by the following reflections $s_i$:
\begin{align}\label{E:Weyl action}
	s_i\la = \la - \pair{\alpha_i^\vee}{\la} \alpha_i\qquad\text{for $i\in I$ and $\la \in X'$.}
\end{align}
We have
\begin{align*}
	s_i(\epsilon_j) = \begin{cases} \epsilon_{i+1} & \text{if $j=i$} \\
		\epsilon_i & \text{if $j=i+1$} \\
		\epsilon_j & \text{if $j\notin \{i,i+1\}$.}
		\end{cases}
\end{align*}
Thus the restriction of the action of $S_\Z$ to the basis $\{\epsilon_i\mid i\in I\}$ of the sublattice $X \subset X'$ is the permutation representation on the set $\Z$: $s_i$ exchanges $i$ and $i+1$ and fixes other integers.

The set $\Phi$ of roots are the elements of the form $w\alpha_i$ for $i\in I$ and $w\in S_\Z$; they have the form $\alpha_{ij} = \epsilon_i-\epsilon_j$ for $i,j\in \Z$ with $i\ne j$. Let $s_{ij}\in S_\Z$ be the associated reflection; acting on $\Z$ it exchanges $i$ and $j$ and fixes other integers. Let $\Phi^+$ be the set of positive roots, the $\alpha_{ij}$ with $i<j$. Let $\Phi^-=-\Phi^+$.

\subsection{Dynkin diagram automorphisms}
\label{SS:Dynkin autos}
Let $\Aut(A_\Z)$ denote the group of automorphisms of the diagram $A_\Z$, the permutations of the node set $I=\Z$ which preserve adjacency. $\Aut(A_\Z)$ is generated by the Dynkin shift $\gamma(i)=i+1$ for all $i\in I$ and the Dynkin reversal $\omega(i)=-i$ for all $i\in I$. This is an infinite dihedral group: $\omega^2=\id$ and 
\begin{align}\label{E:omega and gamma}
	\omega \gamma \omega = \gamma^{-1}.
\end{align}
The group $\Aut(A_\Z)$ acts by automorphisms on $X'$: we have $\zeta(\La_i)=\La_{\zeta(i)}$ for all $\zeta\in \Aut(A_\Z)$ and $i\in I$. In particular, we have 
\begin{align}
	\gamma(\epsilon_i) &= \epsilon_{i+1} \\
\label{E:omega epsilon}
\omega(\epsilon_i) &= - \epsilon_{1-i}\qquad \text{for  all $i\in I$.}
\end{align}
There is an induced action of $\Aut(A_\Z)$ by automorphisms on $S_\Z$: $\zeta(s_i) = s_{\zeta(i)}$ for all $\zeta\in \Aut(A_\Z)$ and $i\in I$. We have
\begin{align}
	\gamma(s_i) &= s_{i+1}\\
	\label{E:omega s}
	\omega(s_i) &= s_{-i} \qquad\text{for all $i\in I$.}
\end{align}

\subsection{Thin infinite flag variety}
Let $\C((t))$ and $\C((t^{-1}))$ denote the Laurent polynomial rings.  Denote
$$
E_a:= \left\{\sum_{i=a}^\infty c_i t^i\right\} \subset \C((t)) \qquad \text{and} \qquad F_a:=\left\{\sum_{i=a}^{-\infty} c_i t^i\right\} \subset \C((t^{-1})) 
$$
so that we have $\cdots \subset E_{1} \subset E_0 \subset E_{-1} \subset \cdots$ and $\cdots \subset F_{-1} \subset F_0 \subset F_1 \subset \cdots$.

A subspace $\Lambda \subset \C((t))$ (resp. $\Xi \subset \C((t^{-1}))$ is called {\it admissible} if we have $E_N \subset \Lambda \subset E_{-N}$ (resp. $F_{-N} \subset \Xi \subset F_{N}$) for some $N$.  
The {\it Sato Grassmannian} $\Grd$ (resp. $\Grd_-$) consists of the set of all admissible subspaces in $\C((t))$ (resp. $\C((t^{-1}))$) and has the structure of an ind-variety over $\C$.

The virtual dimensions $\vdim(\Lambda)$ and $\vdim(\Xi)$ are given by 
\begin{align*}
\vdim(\Lambda)&:=\dim(\Lambda/(\Lambda \cap E_0)) -\dim(E_0/(\Lambda \cap E_0)) \\
\vdim(\Xi)&:=\dim(\Xi/(\Xi \cap F_0))-\dim(F_0/(\Xi \cap F_0)). 
\end{align*}
The virtual dimension measures the size of the subspace $\Lambda$ (resp. $\Xi$)  relative to the standard subspace $E_0$ (resp. $F_0$).
Let $\Gr^{(d)}\subset \Grd$ (resp. $\Gr^{(d)}_- \subset \Grd_-$) consist of the admissible subspaces of virtual dimension $d$. Thus, $\Grd$ is the disconnected union of the $\Gr^{(d)}$ for $d\in\Z$.

An {\it admissible flag} in $\C((t))$ of virtual dimension 0 is a sequence
$$
\Lambda_\bullet = \{\cdots \subset \Lambda_1 \subset \Lambda_0 \subset \Lambda_{-1} \subset \cdots \}
$$
of admissible subspaces satisfying the conditions:  (1) $\vdim(\Lambda_i) = -i$, and (2) for some $N$, we have
$\Lambda_i = E_i$ for all $i$ with $|i| \geq N$.  Similarly, we define admissible flags in $\C((t^{-1}))$.  The {\it thin infinite flag variety} $\Fl$ (resp. $\Fl_-$) consists of the set of all admissible flags in $\C((t))$ (resp. $\C((t^{-1}))$) of virtual dimension 0, and has the structure of an ind-variety over $\C$.

\subsection{Thick infinite flag scheme}
We describe the thick infinite flag scheme in an elementary fashion; see \cite{Ka} for further details.  
We say that a subspace $V \subset \C((t^{-1}))$ is {\it opposed} to a subspace $\Xi\in \Grd_-$ if the composite map $V \subset \C((t^{-1}))\to \C((t^{-1}))/\Xi$ is an isomorphism, and that $V$ is \textit{thick} if it is opposed to some $\Xi\in\Grd_-$. Denote by $\Grb_\Xi$ the set of all $V  \subset \C((t^{-1}))$ opposed to $\Xi$. The thick infinite Grassmannian $\Grb$ is the space of all thick subspaces, that is,
$$
\Grb := \bigcup_{\Xi\in\Grd_-} \Grb_\Xi.
$$
Each $\Grb_\Xi$ is an affine space of countably infinite dimension. For example, for $\Xi=F_0$, the vectors
$t,t^2,\ldots,$ form a basis of $\C((t^{-1}))/F_0$. Thus any thick subspace $V$ opposed to $F_0$ 
has a basis of the form $t + v_1, t^2 + v_2, \ldots$ where $v_1,v_2,\ldots \in F_0$, and the $v_i$ are arbitrary.
Similarly, for any $\Xi\in\Grd_-$, $\Grb_\Xi$ can be identified with $\Xi^\infty \cong \prod_1^\infty \Xi$.  Identifying $\Xi^\infty$ with $\C^\infty \cong \prod_\Z \C \cong \Spec(k[x_i \mid i \in \Z])$, $\Grb_\Xi$ is endowed with the structure of an infinite-dimensional affine space.

As shown in \cite{Ka}, $\Grb$ has the structure of a separated scheme over $\C$, and $\Grb_\Xi \subset \Grb$ are affine open subschemes.  We also define the subscheme $\Grb^i \subset \Grb$ by 
$$
\Grb^i := \bigcup_{\vdim(\Xi) = i} \Grb_\Xi
$$
and we note that $\Grb^i \cap \Grb^j = \emptyset$ if $i \neq j$.

There is an injection $\Grd \hookrightarrow \Grb$, described as follows.  Let $\Lambda \in \Grd$.  Pick $N$ so that $E_{-N} \supset \Lambda \supset E_N$.  Let $\mathcal{B} \subset \Lambda \cap\bigoplus_{i=-N}^{N-1} \C t^i$ project to a basis of $\Lambda/E_N$.  Define 
$$
\Lambdab:= {\rm span}\,(\mathcal{B} \cup \{t^N,t^{N+1},\ldots\} )  \subset \C((t)).
$$
Then $\Lambdab$ is a thick subspace (opposed to $E_{-N-1} \oplus W$ for any complement $W$ of 
$\mathrm{span}({\mathcal{B}})$ in $\bigoplus_{i=-N}^{N-1} \C t^i$), 
and the map $\Lambda \mapsto \Lambdab$ induces an injection $\Grd \hookrightarrow \Grb$.

Let $\Xi_\bullet \in \Fl_-$ be an admissible flag with $\vdim(\Xi_i) = i$.  A sequence $V_\bullet$ of subspaces of $\C((t^{-1}))$
$$
\cdots \subset V_{1} \subset V_0 \subset V_{-1} \subset \cdots
$$
where $V_i$ is opposed to $\Xi_i$ is called a {\it thick flag opposed to $\Xi_\bullet$}.  Let $\Flb_{\Xi_\bullet}$ denote the set thick flags opposed to $\Xi_\bullet$.  Let $\ldots, w_{-1},w_0,w_1,\ldots \in \C((t^{-1}))$ be such that $\C((t^{-1}))/\Xi_i$ is spanned by $w_i,w_{i+1},\ldots$.  
 If $V_\bullet$ is opposed to $\Xi_\bullet$, then there exists vectors 
$$
v_k = w_k + \sum_{j=k-1}^{-\infty} a_{kj} w_j \in \C((t^{-1}))
$$
such that $V_i$ is the span of $v_i,v_{i+1},\ldots$.  The coefficients $a_{kj} \in \C$ are arbitrary and uniquely determine $V_\bullet$, endowing $\Flb_{\Xi_\bullet}$ with the structure of an infinite-dimensional affine space.
The {\it thick infinite flag scheme} 
\begin{equation}\label{eq:Flb}
\Flb = \bigcup_{\Xi_\bullet \in \Fl_-} \Flb_{\Xi_\bullet}
\end{equation}
is the space of all thick flags, and has the structure of a separated scheme.  The construction $\Lambda \mapsto \Lambdab$ induces an injection $\Fl \hookrightarrow \Flb$.  The thin standard flag $E_\bullet \in \Fl$ gives a thick standard flag $\Eb_\bullet \in \Flb$.

We use the flag $F_\bullet \in \Fl_-$ as the basepoint of $\Fl_-$.  For $w \in S_\Z$, we have a point $wF_\bullet \in\Fl_-$ given by
$$
(wF)_i = \prod_{j \in w((-\infty,i])} \C t^j .
$$
Similarly, we have $wE_\bullet \in \Fl$ and $w\Eb_\bullet \in \Flb$.

\begin{lem}
We have $\Flb = \bigcup_{w \in S_\Z} \Flb_{wF_\bullet}$.
\end{lem}
\begin{proof} Let $V_\bullet \in \Flb$.  Then there exists $N$ such that for $|i| \geq N$, we have that $V_i$ is opposed to $F_i$.  The statement then reduces to the corresponding statement for the flag variety of the finite-dimensional vector space $F_N/F_{-N}$.  Namely, any flag $G_\bullet$ in $\C^{2N}$ is opposed to (at least) one of the $(2N)!$ flags $wH_\bullet$, where $w \in S_{2N}$ is a permutation and $H_\bullet$ is some choice of basepoint flag.
\end{proof}

\begin{ex}
Define
$$
V_i = \begin{cases} {\rm span}(t^{i+1}, t^{i+2},\ldots) & \text{ if $i$ is even} \\
 {\rm span}(t^{i}, t^{i+2},\ldots) & \text{ if $i$ is odd}
\end{cases}
$$
so $V_0 =  {\rm span}(t,t^2,t^3,\ldots)$, $V_1 =   {\rm span}(t,t^3,\ldots)$, $V_2 =  {\rm span}(t^3,t^4,\ldots)$.  It is clear that $ \cdots \subset V_2 \subset V_1 \subset V_0 \subset \cdots$ and each $V_i$ is opposed to some admissible subspace.  So it belongs to the set that Kashiwara has in \cite[p. 190]{Ka}.   However, it does not belong to $\Flb$ from $\eqref{eq:Flb}$ since it is not true that $V_m$ is opposed to $F_m$ for large $m$, so it cannot be opposed to any admissible flag.
\end{ex}

We also define the thick partial flag schemes $\Flb_k$ for $k \in \Z$, considering of partial flags $V_\bullet$ in $\C((t^{-1}))$ indexed by $i \in \Z \setminus \{k\}$, with $V_i$ opposed to $\Xi_i$ for some $\Xi_\bullet \in \Fl_-$.  There is a natural morphism $\Flb \to \Flb_k$ defined by forgetting $V_k$.


\subsection{Infinite Schubert varieties}
\def\End{{\rm End}}
Let $B \subset \End(\C((t^{-1})))$ consist of the invertible linear transformations $\varphi: \C((t^{-1})) \to \C((t^{-1}))$ satisfying
$$
\varphi(t^i) \in c_i t^i + \bigoplus_{j>i} \C t^j  \qquad c_i \neq 0.
$$
Let $B_- \subset \End(\C((t^{-1})))$ consist of the invertible linear transformations $\varphi: \C((t^{-1})) \to \C((t^{-1}))$ satisfying
$$
\varphi(t^i) \in c_i t^i + \prod_{j<i} \C t^j  \qquad c_i \neq 0.
$$
Let $T = T_\Z = B \cap B_-$.  The groups $B, B_-$ act on $\Flb$.  The group $B$ acts on $\Fl$, but $B_-$ does not.  We have that $B\cdot E_\bullet = E_\bullet$ and $B \cdot \Eb_\bullet = \Eb_\bullet$.  The {\it (thin) Schubert cells} and {\it (thin) Schubert varieties} are defined to be 
$$
\Omega_w := B \cdot wE_\bullet \subset \Fl \qquad \text{and} \qquad X_w := \overline{B \cdot wE_\bullet} \subset \Fl.
$$
We have the decomposition $\Fl = \bigsqcup_{w \in S_\Z} \Omega_w$.  Furthermore, $\Omega_w \cong \C^{\ell(w)}$ and $X_w$ is an irreducible projective variety of dimension $\ell(w)$. 

An order ideal in $S_\Z$ is a subset $S\subset S_\Z$ such that if $x<y\in S_\Z$ with $y\in S$ then $x\in S$.
Let $K_T(\Fl)$ denote the $K$-group of $T$-equivariant coherent sheaves on $\Fl$.  Since $\Fl$ is an ind-scheme, the $K$-group $K_T(\Fl)$ is the inductive limit 
$$
\varinjlim_{S} K_T(\bigcup_{w \in S} X_w) = \bigoplus_{w \in S_\Z} K_T(\pt) [\cO_{X_w}].
$$

The {\it (thick) opposite Schubert cells} and {\it (thick) opposite Schubert varieties} are defined to be
$$
\Omega^w := B_- \cdot w\Eb_\bullet \subset \Flb \qquad \text{and} \qquad X^w := \overline{B_- \cdot w\Eb_\bullet} \subset \Flb.
$$
We have the decomposition $\Flb = \bigsqcup_{w \in S_\Z} \Omega^w$.  Furthermore, $\Omega^w$ is an affine space of countably infinite dimension and of codimension $\ell(w)$ in $\Flb$.  

For any finite order ideal $S\subset S_\Z$, let $\Omega^S = \bigcup_{w\in S} w \Omega^e = \bigsqcup_{w\in S} \Omega^{w}$. Let $K_T(\Omega^S)$ be the Grothendieck group of coherent $T$-equivariant $\cO_{\Omega^S}$-modules. One may show that $K_T(\Omega^S) \cong \bigoplus_{w\in S} K_T(\pt) [\cO_{X^w}]$ \cite[Lemma 2.3]{KS}.
Define $$K_T(\Flb) = \varprojlim_S K_T(\Omega^S) = \prod_{w\in S_\Z} K_T(\pt) [\cO_{X^w}].$$


\subsection{NilHecke ring, localization and the GKM ring}
We recall some results from \cite{KS}.  These results are stated in the Kac-Moody setting, but the proofs are valid for our thick infinite flag scheme $\Flb$.

We have $K_T^*(\pt)\cong R(T) \cong \Z[X]$ (see \S \ref{SS:Dynkin AZ} for the definition of $X$) with the image of $\la\in X$ written $e^\la$.  
The infinite symmetric group $S_\Z$ acts on $T$ and therefore on $R(T)=\Z[X]$. 
Let $R(T)_\loc = R(T)[(1-e^\alpha)^{-1}\mid \alpha\in \Phi]$ and let $\K_\loc = R(T)_\loc[S_\Z]$ be the twisted group ring with multiplication $(f v)(gw) = f v(g) vw$ for $f,g\in R(T)_\loc$ and $v,w\in S_\Z$. Then $\K_\loc$ acts on $R(T)_\loc$.
Let $D_i \in \K_\loc$ be the element $(1-e^{-\alpha_i})^{-1} (1 - e^{-\alpha_i} s_i)$ for $i\in \Z$.
Let $\K$ be the $K$-theoretic nilHecke algebra, the subring  of $\K_\loc$ generated by $R(T)$ and the $D_i$.  The algebra $\K$ acts naturally on $R(T)$. In this context $D_i$ is known as the \emph{Demazure operator} \cite{Dem}.

Let $\Fun(S_\Z,R(T))$ be the $R(T)$-algebra of functions $f:S_\Z \to R(T)$ under pointwise product, and similarly define $\Fun(S_\Z,R(T)_\loc)$.
The algebra $\K$ acts on $\Fun(S_\Z,R(T)_\loc)$ by
\begin{align}
	D_i(f)(w) &= \dfrac{1}{1-e^{-w(\alpha_i)}} f(w) - \dfrac{e^{-w(\alpha_i)}}{1-e^{-w(\alpha_i)}} f(ws_i).
\end{align}
The $T$-fixed points on $\Flb$ are the flags $w \Eb_\bullet$, and we identify $\Flb^T$ with $S_\Z$.  For each $w\in S_\Z$, localization at $w$ gives a $R(T)$-module homomorphism
$i_w^*: K_T(\Flb) \to R(T)$ \cite[\S 2]{KS}. It satisfies
\begin{align}\label{E:loc support}
	i_w^*([\cO_{X^v}]) = \begin{cases}
		\displaystyle{\prod_{\alpha\in \Phi^+ \cap w(\Phi^-)}} (1-e^{\alpha}) & \text{if $v=w$} \\
		0 & \text{unless $v\le w$.}
	\end{cases}
\end{align}
The map $S_\Z  \cong \Flb^T\hookrightarrow \Flb$ induces the $R(T)$-algebra homomorphism
$\res: K_T(\Flb) \to \Fun(S_\Z,R(T))$ given by $\res(f)(w) = i_w^*(f)$. Define
\begin{align}
	\psi^v = \res([\cO_{X^v}]) \qquad\text{for $v\in S_\Z$.}
\end{align}
Recall that for $i\in \Z$ we have the projection $p_i: \Flb \to \Flb_i$ to the $i$-th minimal thick partial flag scheme.
We have \cite[Cor. 3.3]{KS}
\begin{align}\label{E:pushpull}
	p_i^*p_{i*}(\cO_{X^v}) = \begin{cases}
		[\cO_{X^{vs_i}}] & \text{if $vs_i < v$} \\
		[\cO_{X^v}] & \text{if $vs_i>v$.}
	\end{cases}
\end{align}

The following result is due to Kashiwara; see \cite[Prop. 3.3]{LSS}.
\begin{prop} \label{P:pushpull D} The following diagram commutes:
	\begin{equation*}
		\begin{tikzcd}
		K_T(\Flb) \arrow[r,"\res"] \arrow[d,swap,"p_i^*p_{i*}"]& \Fun(S_\Z,R(T)) \arrow[r,hook] &\Fun(S_\Z,R(T)_\loc) \arrow[d,"D_i"] \\
		K_T(\Flb) \arrow[r,"\res"]& \Fun(S_\Z,R(T)) \arrow[r,hook] &  \Fun(S_\Z,R(T)_\loc)
		\end{tikzcd}
	\end{equation*}
\end{prop}

\begin{prop} \label{P:Schubert loc}
There are functions $\{\psi^v\mid v\in S_\Z\}\subset \Fun(S_\Z,R(T))$ which are uniquely determined by the following conditions:
	\begin{enumerate}
		\item $\psi^v(\id)=\delta_{\id,v}$ for $v\in S_\Z$.
		\item If $ws_i<w$ then
		\begin{align}
			\psi^v(w) &= \begin{cases} 
				\psi^v(ws_i) & \text{if $vs_i>v$} \\
				e^{-w(\alpha_i)} \psi^v(ws_i) + (1- e^{-w(\alpha_i)}) \psi^{vs_i}(w) & \text{if $vs_i<v$.}
			\end{cases}
		\end{align}
		\item We have
		\begin{align}\label{E:root values}
			\psi^v(w)\in \Z[Q]\qquad\text{for all $v,w\in S_\Z$.}
		\end{align}
	\end{enumerate}
\end{prop}
\begin{proof} Part (1) follows from the support condition \eqref{E:loc support}. Part (2) follows from 
	\begin{align}
		D_i(\psi^v) = \begin{cases}
			\psi^{vs_i} & \text{if $vs_i<v$} \\
			\psi^v & \text{otherwise}
		\end{cases}
	\end{align}
	which holds by equation \eqref{E:pushpull}, Proposition \ref{P:pushpull D}, and the definition of $\psi^v$. Equation \eqref{E:root values} holds by induction.
\end{proof}

Let $\Psi$ be the set of $\psi\in \Fun(S_\Z,R(T))$ such that
\begin{align}\label{E:GKM condition}
	1 - e^\alpha \mid \psi(s_\alpha w) - \psi(w)\qquad\text{for all $\alpha\in \Phi$ and $w\in S_\Z$.}
\end{align}
\begin{prop} \label{P:GKM basis} The space $\Psi$ is an $R(T)$-subalgebra of $\Fun(S_\Z,R(T))$. Moreover 
	\begin{align}\label{E:GKM basis}
		\Psi = \prod_{v\in S_\Z} R(T) \psi^v.
	\end{align}
\end{prop}
The ring $\Psi$ is called the GKM ring \cite{GKM}.  In the case of Kac-Moody flag varieties the analogous condition is due to Kostant and Kumar \cite{KK}. The action of $\K$ on $\Fun(S_\Z,R(T))$ preserves $\Psi$.

\begin{prop} The map $\res$ induces an isomorphism $K_T(\Flb)\cong \Psi$ of $R(T)$-algebras and $\K$-algebras where $D_i$ acts by $p_i^*p_{i*}$.
\end{prop}

If there is an action of a group $G$ on an algebra $R$ and an $R$-module $M$, then we say that $M$ is a \textit{$G$-equivariant} $R$-module if $g(a)\cdot g(m) = g(a\cdot m)$ for all $g\in G$, $r\in R$, and $m\in M$.

\begin{prop}  \label{P:Aut GKM} \
\begin{enumerate}
\item[(a)]
	The group $\Aut(A_\Z)$ acts on $\K$ by conjugation and acts on $R(T)$, making $R(T)$ into an $\Aut(A_\Z)$-equivariant $\K$-module.
	\item[(b)] $\Aut(A_\Z)$ acts on $\Fun(S_\Z,R(T))$ by conjugation, stabilizing the subring $\Psi$, making $\Fun(S_\Z,R(T))$ and $\Psi$ into $\Aut(A_\Z)$-equivariant $\K$-modules.
\item[(c)]
	For every $\zeta\in\Aut(S_\Z)$ and $v\in S_\Z$, we have $\zeta(\psi^v)=\psi^{\zeta(v)}$.
\end{enumerate}
\end{prop}

\subsection{Multiplicative formal group law}
Let $A$ be a ring. Define the binary operation $\oplus$ on $A$ by
\begin{align}
	\label{E:oplus} a \oplus b &:= a + b - ab
\end{align}
For $b\in A$ such that $1-b$ is a unit in $A$, define 
\begin{align}
	\label{E:oneg}
	\ominus b &:= \dfrac{-b}{1-b}\qquad\text{for $b\ne1$.} \\
	\label{E:ominus}
	a\ominus b &:= a \oplus(\ominus b) = \dfrac{a-b}{1-b}.
\end{align}
The operation $\oplus$ is commutative and associative with neutral element $0$ and $b\ominus b=0$ for $b\ne 1$.

\subsection{From exponentials to polynomials} 
\label{SS:change of variable}
Define the Laurent polynomial ring $R$ and its subrings $R^+$ and $R^-$, as follows:
\begin{align} \label{E:K ring}
	R &:= \Z[(1-x_i)^{\pm1}\mid i\in \Z] \\
	\label{E:R+}
	R^+ &:= \Z[(1-x_i)\mid i\in \Z_{>0}][(1-x_i)^{-1}\mid i \in \Z_{\le0}] \\
	\label{E:R-}
	R^- &:= \Z[(1-x_i)\mid i\in \Z_{\le 0}][(1-x_i)^{-1}\mid i \in \Z_{>0}].
\end{align}
Note that $\Z[\ominus(x_i)]=\Z[(1-x_i)^{-1}]$. 

Let $\Theta: R(T) \to R$ be the ring isomorphism given by $\Theta(e^{\epsilon_i}) = (1-x_i)^{-1}$ for $i\in \Z$. This is merely a renaming of the generators of a Laurent polynomial ring, but is convenient for combinatorial applications. We have
\begin{align}\label{E:Theta root}
\Theta(1 - e^{\alpha_{ij}}) = x_j \ominus x_i\qquad\text{for all $i,j\in\Z$ with $i\ne j$.}
\end{align}
 The isomorphism $\Theta$ is $S_\Z$-equivariant.


Let $\Kt$ be the algebra generated by $R$ and elements $\Gdiff_i$ for $i \in I$, equipped with an isomorphism $\Theta: \K \to \Kt$ which restricts to the isomorphism $\Theta: R(T) \to R$ and satisfies $\Theta(D_i) = \Gdiff_i$.
We have
\begin{align*}
	\Gdiff_i &=\Theta(D_i) \\
	&= \Theta \left((1-e^{-\alpha_i})^{-1}(1-e^{-\alpha_i} s_i)\right)\\
 &=(x_i\ominus x_{i+1})^{-1} ( 1 - (1-x_{i+1})^{-1}(1-x_i) s_i) \\
	&=(x_i-x_{i+1})^{-1} (1-x_{i+1}-(1-x_i)s_i) \\
	&= (x_i-x_{i+1})^{-1}(1-s_i)(1-x_{i+1}) \\
	&= A_i \circ (1-x_{i+1}),  
\end{align*}
where $A_i = (x_i-x_{i+1})^{-1}(1-s_i)$ is the usual divided difference operator\cite{LLS:back stable}.
Consider now the $\Aut(A_\Z)$-equivariant $\K$-module structure on $R(T)$.
By \eqref{E:omega epsilon}, we have
\begin{align}
	\label{E:omega on R(T)}
	\omega(e^{\epsilon_i}) &= e^{-\epsilon_{1-i}}\qquad\text{for $i\in I$,}
 \\
	\label{E:gamma on R(T)}
	\gamma(e^{\epsilon_i}) &= e^{\epsilon_{i+1}}\qquad\text{for $i\in I$.}
\end{align}
Via the isomorphism $\Theta$, we obtain an $\Aut(A_\Z)$-equivariant $\Kt$-module structure on $R$. The elements $\omega$ and $\gamma$ of $\Aut(S_\Z)$ induce automorphisms of $R$ and $\Kt$ that we denote by $\omk$ and $\gamma$ respectively.
  By \eqref{E:omega on R(T)}, we have $\omk((1-x_i)^{-1}) = 1 - x_{1-i}$. Thus for all $i\in \Z$, we have
\begin{align}
	\label{E:omk x}
	\omk(x_i) &= \ominus(x_{1-i})\\
	\omk(\Gdiff_i) &= \omk(\Gdiff_{-i})\\
	\gamma(x_i) &= x_{i+1} \\
	\gamma(\Gdiff_i) &= \Gdiff_{i+1}.
\end{align}
Furthermore, we have
\begin{align}\label{E:omk and gamma}
	\omk \gamma \omk = \gamma^{-1}
\end{align}
on $\Kt$ and on $R$, by \eqref{E:omega and gamma}.

Let $R(a)$ be the ring isomorphic to $R$, using the variables $a_i$ instead of $x_i$.
We write $\Theta_a:R(T)\to R(a)$ for the isomorphism involving the $a_i$ variables.


Let $(\Theta_a)_*:\Fun(S_\Z,R(T)) \cong \Fun(S_\Z,R(a))$ be defined by $(\Theta_a)_*(f)= \Theta_a \circ f$. 
Let $\Psit\subset\Fun(S_\Z,R(a))$ be the image of $\Psi$ under $(\Theta_a)_*$. Applying $\Theta_a$ to \eqref{E:GKM condition}, using \eqref{E:Theta root}, and observing that $a_j\ominus a_i$ and $a_j-a_i$ generate the same ideal of $R(a)$,
we see that $\Psit$ is the subring of functions $f: S_\Z\to R(a)$ such that
\begin{align}\label{E:transformed GKM}
	a_j-a_i \mid f(w) - f(s_{ij}w)\qquad\text{for all $w\in S_\Z$ and $i,j\in\Z$ with $i\ne j$.} 
\end{align}
This induces an isomorphism $(\Theta_a)_*:\Psi\to\Psit$.  Finally, for $v\in S_\Z$, let $\psit_v\in \Psit$ be defined by $\psit^v = (\Theta_a)_*(\psi^v)=\Theta_a\circ \psi^v$.

\section{Grothendieck polynomials}

\subsection{Demazure operators}

For $i\in \Z$ and $w\in S_\Z$, define 
\begin{align*}
	s_i * w := \begin{cases}
		s_i w & \text{if $s_i w > w$} \\
		w & \text{otherwise.}
		\end{cases}
\end{align*}
This defines a monoid $(S_\Z,*)$ called the $0$-Hecke monoid. The operation $*$ is called the \emph{Hecke product} or \emph{Demazure product}.  Recall the operators $\Gdiff_i\in\Kt$ defined in \S \ref{SS:change of variable}.

\begin{lem} \label{L:Gdiff}\
\begin{enumerate}
\item The $\Gdiff_i$ satisfy the braid relations $\Gdiff_i \Gdiff_{i+1} \Gdiff_i = \Gdiff_{i+1} \Gdiff_{i} \Gdiff_{i+1}$ , and $\Gdiff_i^2=\Gdiff_i$.
\item The $\Gdiff_i$ generate a monoid isomorphic to the $0$-Hecke monoid.
\end{enumerate}
\end{lem}

\subsection{Subgroups of permutations} Our notation for symmetric groups and partitions follows \cite{LLS:back stable}.  Define the following subgroups of $S_\Z$:
$$S_+ :=  \langle s_i\mid i\in \Z_{>0}\rangle, \qquad 
S_- :=  \langle s_i\mid i\in \Z_{<0}\rangle, \qquad
S_{\ne0} := S_- \times S_+ =  \langle s_i\mid i\in \Z\setminus\{0\}\rangle.$$

\subsection{Grothendieck polynomials}

For $w\in S_n$, the \emph{Grothendieck polynomial} $\G_w\in \Z[x_1,\dotsc,x_n]$ is defined by
\cite{LaSc} \cite[\S 4]{FK}
\begin{align}
\label{E:fG long}
   \G_{w_0} &:= x_1^{n-1} x_2^{n-2}\dotsm x_{n-1} \\
 \label{E:fG ddiff}
   \G_w &:= \Gdiff_i \G_{ws_i}\qquad\text{for $ws_i>w$.}
\end{align}

\begin{ex} The Grothendieck polynomials $\G_w$ for $w\in S_3$ are given below with $\G_{s_1s_2s_1}$ at the top.
	\begin{equation*}
		\begin{tikzcd}
			& x_1^2x_2 \arrow[dl,swap,"\Gdiff_1"] \arrow[dr,"\Gdiff_2"] \\
			x_1x_2 \arrow[d,swap,"\Gdiff_2"]&& x_1^2 \arrow[d,"\Gdiff_1"] \\[4pt]
			x_1\arrow[dr,swap,"\Gdiff_1"] && x_1\oplus x_2\arrow[dl,"\Gdiff_2"]  \\ 
			& 1
		\end{tikzcd}
	\end{equation*}
\end{ex}

\begin{rem} By Lemma \ref{L:Gdiff}, the polynomials $\G_w \in \Z[x_1,\dotsc,x_n]$ for $w \in S_n$ are well-defined. Moreover 
for $w\in S_+=\bigcup_{n\ge1} S_n$, $\G_w$ is independent of $n$ in the sense that for any $n$ such that $w\in S_n$, 
$\G_{\iota(w)}=\G_w$ under the standard embedding $\iota:S_n\to S_{n+1}$.
\end{rem}

Fomin and Kirillov give the following monomial expansion of $\G_w$.

\begin{prop} \label{P:fG FK} \cite[Prop. 3.3]{FK} For $w\in S_+$, we have
	\begin{align}\label{E:Groth FK}
		\G_w = \sum_{\substack{ s_{a_1}*s_{a_2}*\dotsm *s_{a_p}=w \\
				1 \le i_1\le i_2\le \dotsm \le i_p \\
				a_k \le a_{k+1} \Rightarrow i_k<i_{k+1} \\
				i_k \le a_k  }} (-1)^{p-\ell(w)} x_{i_1} x_{i_2}\dotsm x_{i_p},
	\end{align}
	where $p$ is arbitrary.
\end{prop}

\begin{prop}\label{P:positive basis}
The set $\{\G_w\mid w\in S_+\}$ is a $\Z$-basis of $\Z[\xp]=\Z[x_i \mid i \in Z_{\ge 0}]$.
\end{prop}
\begin{proof} By Proposition \ref{P:fG FK}, the lowest homogeneous part of $\G_w$ coincides with the Billey-Jockusch-Stanley formula for the Schubert polynomial \cite{BJS}.
The Schubert polynomials form a basis \cite[Thm. 2.7]{LLS:back stable}.
This shows that the $\G_w$ are linearly independent. The Monk Rule of \cite{Len2} gives a finite Grothendieck polynomial expansion of any product $x_k \G_w$ for $k\ge1$ and $w\in S_+$. In particular, iterating the Monk rule one may expand any monomial (times $\G_\id=1$) as a finite linear combination of Grothendieck polynomials.  This proves that $\{\G_w\mid w\in S_+\}$ span $\Z[\xp]$.
\end{proof}

\subsection{Negative Grothendieck polynomials}
Recall the automorphisms $\omega$ of $S_\Z$ and $\omk$ of $R$ from \S \ref{SS:change of variable}.  The automorphism $\omega\in\Aut(S_\Z)$ restricts to an isomorphism $S_-\cong S_+$. For $w\in S_-$, define the \emph{negative Grothendieck polynomial} $\G_w\in R^+$ by 
\begin{align}\label{E:G S_-}
	\G_w := \omk(\G_{\omega(w)})\qquad\text{for $w\in S_-$.}
\end{align}

\begin{ex}
	The polynomials $\G_w$ for $w\in \langle s_{-1},s_{-2} \rangle$ are given below, with $\G_{s_{-1}s_{-2}s_{-1}}$ at the top.
	\begin{equation*}
		\begin{tikzcd}
			& (\ominus x_0)^2(\ominus x_{-1}) \arrow[dl,swap,"\Gdiff_{-1}"] \arrow[dr,"\Gdiff_{-2}"] \\
			(\ominus x_0)(\ominus x_{-1}) \arrow[d,swap,"\Gdiff_{-2}"] && (\ominus x_0)^2 \arrow[d,"\Gdiff_{-1}"] \\[6pt]
			\ominus x_0 \arrow[dr,swap,"\Gdiff_{-1}"] && (\ominus x_0) \oplus (\ominus x_{-1}) \arrow[dl,"\Gdiff_{-2}"] \\
			& 1
		\end{tikzcd}
	\end{equation*}
\end{ex}

\begin{prop} \label{P:negative basis}
The set $\{\G_w\mid w\in S_-\}$ is a $\Z$-basis of $\Z[(1-x_i)^{-1}\mid i\in \Z_{\le0}]$.
\end{prop}
\begin{proof}
	The map $\omk$ restricts to an isomorphism $$\Z[x_i\mid i\in\Z_{>0}]\to \Z[(1-x_i)^{-1}\mid i\in \Z_{\le0}]$$ sending the basis $\{\G_w\mid w\in S_+\}$ to the basis
	$\{\G_w\mid w\in S_- \}$. 
\end{proof}

For $w\in S_{\ne0} = S_- \times S_+$, there is a unique factorization $w=uv$ with $u\in S_-$ and $v\in S_+$. Define
\begin{align} \label{E:G S nonzero}
	\G_w := \G_u \G_v.
\end{align}

\begin{lem}\label{L:fG omk} For $w \in S_\Z$, we have
	\begin{align}\label{E:fG omk}
		\omk(\G_w) = \G_{\omega(w)}\qquad\text{for $w\in S_{\ne0}$.}
	\end{align}
\end{lem}

\begin{prop}\label{P:G poly basis}
The set $\{\G_w\mid w\in S_{\ne0}\}$ is a $\Z$-basis of $R^+=\Z[(1-x_i)\mid i\in\Z_{>0}][(1-x_i)^{-1}\mid i \in \Z_{\le0}]$.
\end{prop}
\begin{proof} Follows by tensoring the bases in Propositions \ref{P:positive basis} and \ref{P:negative basis}.
\end{proof}

\section{Back stable Grothendieck polynomials}

\subsection{Back stable rings}
\label{SS:sym func}
Let 
$$
\xp:=(x_1,x_2,x_3,\dotsc), \qquad \xm=x_{\leq 0} := (x_0,x_{-1},x_{-2},\dotsc).
$$
Let $\La = \La(\xm)$ denote the ring of symmetric functions in the variables $\xm$ with coefficients in $\Z$ and let $\hLa$ denote the graded completion of $\La$.  That is, an element $f \in \hLa$ is a formal linear combination $f = f_1+f_2+\cdots$ where $f_i \in \La$ is homogeneous of degree $i$.  Similarly, we define $\hLa(\xp)$ using positive variables.  When no decorations are present, we assume negative variables are used: $\hLa$ means $\hLa(\xm)$.

Let $p_k, e_k, h_k \in \hLa$ denote the power sum, elementary, and homogeneous symmetric functions.  We have isomorphisms $\hLa \cong \Z[[e_1,e_2,\ldots]]$ and $\hLa \cong \Z[[h_1,h_2,\ldots]$. Despite the fact that $\hLa(\xm) \supsetneq \Z[[p_1,p_2,\ldots]]$, since $\hLa(\xm) \otimes_{\Z} \Q = \Q[[p_1,p_2,\ldots]]$ we define maps from $\hLa$ by describing the images of $p_k$. That is, we define $\Q$-algebra homomorphisms using the topological generators $p_k$ and in all cases they restrict to ring homomorphisms over the integers.  
Define the back stable rings
\begin{align*}
	\bR& := \hLa(\xm) \otimes  R \\
	\bR^+ &:= \hLa(\xm) \otimes R^+.
\end{align*}
The Dynkin shift automorphism $\gamma\in\Aut(A_\Z)$ from \S \ref{SS:Dynkin AZ} induces 
the ring automorphism of $\bR$ defined by $\gamma(x_i)=x_{i+1}$ for $i\in\Z$ and $\gamma(p_k(\xm))=p_k(\xm)+x_1^k$ for all $k\ge1$.

We identify $\hLa(\xp)$ with $\hLa(\xm)$ by setting
\begin{align}\label{E:positive and negative symmetric functions}
p_k(\xp) = - p_k(\xm) \qquad\text{for $k\ge1$.}
\end{align}
This is consistent with the $S_\Z$-action: the element $\sum_{i\in \Z} x_i^k$ is set to zero, and is $S_\Z$-symmetric.  The identification \eqref{E:positive and negative symmetric functions} implies that $e_n(\xp)=(-1)^n h_n(\xm)$ and $h_n(\xp) = (-1)^n e_n(\xm)$, and thus preserves symmetric functions with integral coefficients.

\subsection{Conjugation automorphism}
\label{SS:conj auto bR}
Recall $\omk\in \Aut(R)$ from \S \ref{SS:change of variable}. Since the $p_k$ are algebraically independent topological generators of $\hLa$, we may extend $\omk$ to an automorphism of $\bR$ by
\begin{align}
\label{E:omk power}
	\omk(p_k)  &:= (-1)^{k+1} \sum_{r\ge0} \binom{r+k-1}{k-1} p_{k+r}&&\text{for $k\ge1$.}
\end{align}
Heuristically, we have $\omk(f(x)) = f(\omk(x)) = f(\ominus(x_{1-i}))$ for $f \in \hLa$.  Using \eqref{E:positive and negative symmetric functions}, this induces the map \eqref{E:omk power}.

We call $\omk$ the \emph{($K$-theoretic) conjugation automorphism}.  Note that $\omk$ differs from the automorphism $\omega$ defined in \cite{LLS:back stable}, and restricts to an automorphism of $\hLa$ that differs from the usual conjugation automorphism $\omega$ of symmetric functions.  The automorphisms $\omk$ and $\gamma$ define an action of $\Aut(A_\Z)$ on $\bR$.

\begin{prop}
The maps	 $\gamma$ and $\omk$ define an action of $\Aut(A_\Z)$ on $\bR$. That is,
	\begin{align}
	\omk\circ	\gamma \circ \omk = \gamma^{-1}.
	\end{align}
\end{prop}
\begin{proof} As this relation already holds on $R$ by \eqref{E:omk and gamma}, it is enough to check the identity 
	applied to $p_k$ for $k \geq 1$. We have
\begin{align*}
	 (\gamma \circ \omk)(p_k)  	&=  \gamma((-1)^{k+1} \sum_{r\ge0} \binom{r+k-1}{k-1} p_{k+r})\\
	 &=  (-1)^{k+1} \left( \sum_{r\ge0} \binom{r+k-1}{k-1} p_{k+r} + \sum_{r \geq 0} \binom{r+k-1}{k-1} x_1^{k+r} \right) \\
	 & = \omk(p_k) -\left(\dfrac{-x_1}{1-x_1} \right)^k \\
	 &= \omk(p_k- x_0^k) \\
	&= \omk \circ \gamma^{-1}(p_k). \qedhere
	 \end{align*}
	 \end{proof}
	 The conjugation automorphism $\omk$ restricts to an automorphism of $\bR^+$.

\subsection{Group law negation automorphism}\label{SS:group law negation}
Denote by $\ominus: R \to R$ be the involutive ring automorphism of $R$ defined by $x_i\mapsto \ominus(x_i)$ for all $i$.  We will write $f(\ominus(x))$ for $\ominus(f(x))$ for $f(x)\in R$. The map $\ominus$ restricts to an isomorphism $R^+\cong R^-$.
It satisfies the operator identities
\begin{align}
	\label{E:aut perm}
	\ominus \circ w \circ \ominus &= w \\
	\label{E:aut diff}
	\ominus \circ \Gdiff_i \circ \ominus &= s_i \Gdiff_i s_i =: \Gdiff_i^{\ominus} = A_i(x_i-1).
\end{align}

\begin{lem}\label{L:ddiff and negation}
\begin{align}
	\Gdiff_i^\anti(\G_w(\ominus x)) = \begin{cases} \G_{ws_i}(\ominus x) & \text{if $ws_i<w$} \\
		\G_w(\ominus x) & \text{otherwise.}
	\end{cases}
\end{align}
\end{lem}
We extend $\ominus$ to an automorphism of $\bR$, called the \emph{group law negation automorphism}, by setting 
\begin{align*}
	\ominus(p_k)  &:= -\omk(p_k) = (-1)^{k} \sum_{r\ge0} \binom{r+k-1}{k-1} p_{k+r}&&\text{for $k\ge1$.}
\end{align*}

\subsection{Antipode automorphism}
Define the \emph{antipode} involution $S$ on $\bR$ by $S(p_k)=-p_k$ for all $k\ge1$
and $S(x_i) = x_{1-i}$ for $i\in \Z$. 
With these definitions we have
\begin{align}
	\label{E:S omk commute}
	\ominus = S \circ \omk &= \omk \circ S. 
\end{align}

\begin{prop}
We have $\ominus \circ \gamma = \gamma \circ \ominus$ and $S\circ \gamma^{-1} = \gamma \circ S$.
\end{prop}

\subsection{Back stable Grothendieck polynomials}
\label{SS:bG}
For $w\in S_\Z$ and an interval $[p,q]\subset\Z$ that contains all integers moved by $w$, let $\G^{[p,q]}_w$ be the usual Grothendieck polynomial except computed using variables $x_p,x_{p+1},\dotsc,x_q$ instead of $x_1,x_2,\dotsc$. That is,
\begin{align}\label{E:G interval}
	\G_w^{[p,q]} := \gamma^{p-1} \G_{\gamma^{1-p}(w)}
\end{align}
where $\gamma$ denotes both $\gamma\in \Aut(\bR)$ from \S \ref{SS:sym func} and $\gamma\in\Aut(S_\Z)$ from \S \ref{SS:Dynkin autos}.

The \emph{back stable Grothendieck polynomial} $\bG_w\in \bR$ is defined by
\begin{align}\label{E:G back stable def}
	\bG_w &:= \lim_{\substack{p\to-\infty \\ q \to \infty}} \G_w^{[p,q]}.
\end{align}

It is immediate that we have 
\begin{align}
	\Gdiff_i \G_w = \begin{cases}
		\G_{ws_i} & \text{if $ws_i<w$} \\
		\G_w & \text{otherwise.}
		\end{cases}
\end{align}

Taking the limit of Proposition \ref{P:fG FK} gives the following formula in which the indices $i_k$ are integers.
\begin{prop} \label{P:G back stable} For $w \in S_\Z$, we have
\begin{align}\label{E:G back stable}
	\bG_w = \sum_{\substack{ s_{a_1}*s_{a_2}\dotsm *s_{a_p}=w \\
			i_1\le i_2\le \dotsm \le i_p \\
			a_k \le a_{k+1} \Rightarrow i_k<i_{k+1} \\
			i_k \le a_k  }} (-1)^{p-\ell(w)} x_{i_1} x_{i_2}\dotsm x_{i_p},
\end{align}
where $p$ is arbitrary.
\end{prop}

\begin{ex} Consider $\bG_{s_0}$. For every $p\ge1$, $s_0$ may be obtained as the Hecke product of $p$ copies of $s_0$.
For this Hecke factorization the associated sequence $(i_1,\dotsc,i_p)$ can be any sequence of integers with $i_1 < i_2 <\dotsm< i_p \le 0$. Therefore, we have $\bG_{s_0} = \sum_{p\ge1} (-1)^{p-1} e_p$.
\end{ex}

See Appendix~\ref{S:computations} for more examples of back stable Grothendieck polynomials.

\subsection{$\Aut(A_\Z)$-action on back stable Grothendieck polynomials}
\begin{prop}\label{P:G back stable shift} For $w \in S_\Z$, we have
\begin{align}\label{E:G back stable shift}
	\bG_{\gamma(w)} = \gamma(\bG_w).
\end{align}
\end{prop}
\begin{proof} Follows from Proposition \ref{P:G back stable}.
\end{proof}

\begin{prop}\label{P:omk backstable} For $w \in S_\Z$, we have
\begin{align}
 \bG_{\omega(w)}=\omk(\bG_w). 
\end{align}
\end{prop}
\begin{proof} This holds by Proposition \ref{P:omk double backstable} and Proposition \ref{P:forget Schubert basis}, using the fact that $\For$ (defined in \S \ref{SS:double bG}) and $\omk$ commute.
\end{proof}

\subsection{$K$-Stanley functions}
There is a ring homomorphism $\eta_0:\bR \to \hLa$ given by $x_i\mapsto 0$ for all $i \in\Z$, and $p_k\mapsto p_k$ for all $k\ge1$. This ``sets all $x_i$ to zero except those in $\hLa(\xm)$''.

Define the \emph{$K$-Stanley function} $G_w\in\hLa(\xm)$ for $w\in S_\Z$ by
\begin{align}
	\label{E:K Stanley}
	G_w := \eta_0(\bG_w).
\end{align}

\begin{rem} \label{R:K Stanley name}
We call $G_w$ the $K$-Stanley function because it is the $K$-theoretic analogue of a Stanley function in cohomology \cite{LLS:back stable}.  For $w\in S_+$, the (forward-)stable Grothendieck polynomial of \cite{FK} is the element $G_w(\xp)\in \hLa(\xp)$ defined by
\begin{align*}
	G_w(\xp) &:= \lim_{n\to\infty} \G_{\gamma^n(w)},
\end{align*}
which by \eqref{E:Groth FK} is
\begin{align}\label{E:stable G positive}
	G_w(\xp) = \sum_{\substack{ s_{a_1}*s_{a_2}*\dotsm *s_{a_p}=w \\
			1 \le i_1\le i_2\le \dotsm \le i_p \\
			a_k \le a_{k+1} \Rightarrow i_k<i_{k+1}} }(-1)^{p-\ell(w)} x_{i_1} x_{i_2}\dotsm x_{i_p}.
\end{align}
Note that this gives the same symmetric series as the definition in \eqref{E:K Stanley}, only with variables
$\xp$ rather than $\xm$.

In this work, we use the name ``$K$-Stanley" instead of ``stable Grothendieck" because the latter produces a conflict in the equivariant setting as there are three different versions of limiting double Grothendieck polynomial: (1) a back stable limit (the back stable double Grothendieck polynomial) which is not symmetric in the $x$ variables, (2) a forward stable limit (the super $K$-Stanley function), which is supersymmetric, and (3) an equivariant analogue of \eqref{E:K Stanley} (the double $K$-Stanley functions).
\end{rem}

\begin{prop}\label{P:wrong way shift}  We have
	\begin{align}\label{E:eta0 shift}
		\eta_0 \circ \gamma &= \eta_0 \\
	\label{E:G shift invariant}
	G_{\gamma(w)} &= G_w \qquad\text{for all $w\in S_{\Z}$.}
	\end{align}
\end{prop}
\begin{proof} As both sides are algebra maps one may check \eqref{E:eta0 shift} on the topological algebra generators, which is straightforward. Equation \eqref{E:G shift invariant} follows from equation \eqref{E:eta0  shift}.
\end{proof}

\begin{prop}\label{P:K Stanley} For $w \in S_\Z$, we have
	\begin{align}\label{E:K Stanley FK}
		G_w = \sum_{\substack{ s_{a_1}*s_{a_2}*\dotsm *s_{a_p}=w \\
				i_1\le i_2\le \dotsm \le i_p \le 0 \\
				a_k \le a_{k+1} \Rightarrow i_k<i_{k+1} }} 
		(-1)^{p-\ell(w)} x_{i_1} x_{i_2}\dotsm x_{i_p}.
	\end{align}
\end{prop}
\begin{proof} By \eqref{E:G shift invariant} we may assume that $w\in S_+$ . Since $ws_i>w$ for $i\in\Z_{<0}$,
$\bG_w$ is $S_-$-symmetric. Hence may write $\bG_w = \sum_\alpha f_\alpha \xp^\alpha$ where $f_\alpha\in \hLa(\xm)$ with only finitely many $f_\alpha\ne0$. Applying $\eta_0$ to \eqref{E:G back stable} we have $G_w = f_0$. But $f_0$ equals the right hand side of \eqref{E:K Stanley FK}, the sum of monomials in \eqref{E:G back stable} having only $x_i$ with $i\le0$.
\end{proof}

\begin{prop}\label{P:G inverse omega}  For $w \in S_\Z$, we have
\begin{align}\label{E:G inverse omega}
	G_{\omega(w^{-1})} = G_w. 
\end{align}
\end{prop}
\begin{proof} By \eqref{E:G shift invariant} one may assume $w\in S_+$
so that $w\in S_n$ for some $n$. By \cite[Cor. 5.10]{FK} and the $S_n$-invariance of $G_w(x_1,\dotsc,x_n)$ we have
\begin{align}\label{E:S_n Dynkin reversal}
		G_w(x_1,\dotsc,x_n) &= \sum_{u*v=w} (-1)^{\ell(u)+\ell(v)-\ell(w)} \G_u w_0 (\G_{w_0v^{-1} w_0})\\
\label{E:S_n Dynkin reversal II}
		&= 	\sum_{u*v=w} (-1)^{\ell(u)+\ell(v)-\ell(w)} w_0(\G_u) \G_{w_0v^{-1} w_0}.
\end{align}
Since conjugation by $w_0$ is a length-preserving automorphism of $S_n$ and taking inverses is a length-preserving anti-automorphism of $S_n$, $u*v=w$ if and only if $(w_0 v^{-1} w_0)*(w_0 u^{-1} w_0) = w_0w^{-1}w_0$.
Applying \eqref{E:S_n Dynkin reversal} for $w_0w^{-1}w_0$ and \eqref{E:S_n Dynkin reversal II} for $w$ we have
\begin{align*}
	G_{w_0w^{-1}w_0}(x_1,\dotsc,x_n) &= \sum_{u*v = w} (-1)^{\ell(u)+\ell(v)-\ell(w)} \G_{w_0v^{-1}w_0} w_0(\G_u) \\
	&= G_w(x_1,\dotsc,x_n)
\end{align*}
Letting $n\to\infty$ we deduce that $G_w = G_{w_0w^{-1}w_0}$ in $\hLa$. 

The Dynkin reversal on $S_n$ given by conjugation by $w_0$, sends $s_i$ to $s_{n-i}$,
whereas $\omega$ sends $s_i$ to $s_{-i}$. Both are group homomorphisms, so 
$w_0 w^{-1} w_0= \gamma^n \omega(w)$. The result follows by \eqref{E:G shift invariant}.
\end{proof}
\begin{prop}\label{P:omega G} For all $w\in S_\Z$, we have
\begin{align}
\label{E:omega G perm}
\omk(G_w) &= G_{\omega(w)} \\
\label{E:antipode on G}
S G_w(x) &= G_{w^{-1}}(\ominus x)
\end{align}
\end{prop}
\begin{proof} 
Equation \eqref{E:omega G perm} follows from Proposition \ref{P:omk backstable} and applying the specialization $x_i\mapsto0$, which commutes with $\omk$.
For \eqref{E:antipode on G}, we calculate
\begin{align*}
	S G_{w^{-1}}(\ominus x)
	= (S \circ \ominus) G_{w^{-1}}
	= \omk G_{w^{-1}} 
	= G_{\omega(w^{-1})}
	= G_w
\end{align*}
by \eqref{E:S omk commute}, \eqref{E:omega G perm}, and \eqref{E:G inverse omega}.  Applying $S$ yields \eqref{E:antipode on G}.
\end{proof}

\subsection{Coproduct on symmetric functions}\label{SS:coproduct}
Let $\Delta:\hLa\to \hLa\otimes\hLa$ be the coproduct: $\Delta(p_k) = p_k \otimes 1 + 1 \otimes p_k$.
We identify $\hLa\otimes\hLa$ with symmetric series in two sets of variables, one set for each tensor factor.
If we use $\xm$ for the first factor and $\am$ for the second then $\Delta(p_k)=p_k(\xm)+p_k(\am)$.
\begin{prop}\label{P:G coproduct} For $w \in S_\Z$, we have
	\begin{align}\label{E:G coproduct}
		\Delta(G_w) &= \sum_{u*v =w} (-1)^{\ell(u)+\ell(v)-\ell(w)} G_u \otimes G_v.
	\end{align}
\end{prop}
\begin{proof} Plug two sets of variables into \eqref{E:K Stanley FK}.
\end{proof}

\subsection{Superization}
For $f(x)\in\hLa$ let $f(x/a)$ denote the image of $f$ in $\hLa\otimes\hLa=\hLa(x)\otimes \hLa(a)$ under superization, which is the coproduct $\Delta$ followed by $\id \otimes S$ where the antipode $S$ acts on the ``a" variables. The map $f\mapsto f(x/a)$ is the unique $\Z$-algebra homomorphism sending $p_k(\xm)$ to $p_k(x/a) = p_k(\xm)-p_k(\am)$, which for historical reasons we denote by $p_k(x||a)$. We call $G_w(x/a)$ the \emph{super $K$-Stanley function}.

\begin{prop}\label{P:G super} For $w \in S_\Z$, we have
	\begin{align}\label{E:G super}
		G_w(x/a) &= \sum_{u*v=w} (-1)^{\ell(u)+\ell(v)-\ell(w)} 
		G_{u^{-1}}(\ominus(a)) G_v(x) \\
	\label{E:G super alt}
		&= \sum_{u*v=w} (-1)^{\ell(u)+\ell(v)-\ell(w)} 
		G_{v^{-1}}(\ominus(a)) G_u(x).
	\end{align}
\end{prop}
\begin{proof} This holds by  Propositions 
	\ref{P:omega G} and \ref{P:G coproduct}.
\end{proof}

\subsection{Coproduct formula}
The coproduct $\Delta$ on $\hLa$ can be extended to a map
$\Delta: R \to \bR$ giving $R$ the structure of a $\hLa$-comodule, by
defining $\Delta(x_i) = 1\otimes x_i$ for $i\in\Z$.
It restricts to a map $R^+ \to \bR^+$ that makes $R^+$ into a $\hLa$-comodule.

\begin{thm}\label{T:coproduct} For $w \in S_\Z$, we have
\begin{align}
\label{E:coproduct}
  \Delta(\bG_w) &= \sum_{u*v=w} (-1)^{\ell(u)+\ell(v)-\ell(w)} G_u \otimes \bG_v, \\
\label{E:little coproduct}
  \bG_w &= \sum_{\substack{u*v=w \\ v\in S_{\ne0}}} (-1)^{\ell(u)+\ell(v)-\ell(w)}
  G_u \G_v.
\end{align}
\end{thm}
\begin{proof} We first derive \eqref{E:coproduct} from \eqref{E:little coproduct}. By Proposition \ref{P:G coproduct}, we have
\begin{align*}
	\Delta(\bG_w) &= \Delta\left(\sum_{\substack{u*v=w\\ v\in S_{\ne0}}} (-1)^{\ell(u)+\ell(v)-\ell(w)} G_u \G_v \right) \\
	&= \sum_{\substack{u_1*u_2*v=w \\ v\in S_{\ne0}}} (-1)^{\ell(u_1)+\ell(u_2)+\ell(v)-\ell(w)} G_{u_1} \otimes G_{u_2} \G_v \\
	&= \sum_{u_1*z=w} (-1)^{\ell(u_1)+\ell(z)-\ell(w)} G_{u_1} \otimes \bG_z.
\end{align*}
The equality \eqref{E:little coproduct} can be deduced from \cite[Theorem 4.7]{LLS:coprod} by taking a limit; see \cite[Section 6.2]{LLS:coprod} for an explanation of this limit in the very similar cohomology setting.  Presumably, \eqref{E:little coproduct} could also be deduced from a direct combinatorial argument similar to the proof of the coproduct formula in \cite{LLS:back stable}.
\end{proof}

\begin{ex} We compute $G_{s_1}$ and $G_{s_{-1}}$ using Theorem \ref{T:coproduct}. 
By \eqref{E:G shift invariant} we have $G_{s_1}=G_{s_{-1}}=G_{s_0}$. By \eqref{E:G Grass} we have $G_{s_0} = G_1$, the Grassmannian $K$-Stanley for the partition $(1)$. Therefore
\begin{align*}
	\bG_{s_1} &= G_{s_1} + \G_{s_1} - G_{s_1} \G_{s_1} = G_1 \oplus x_1 \\
	\bG_{s_{-1}} &= G_1 \ominus x_0.
\end{align*}

\end{ex}

\subsection{Grassmannian $K$-Stanley functions}
\label{SS:Grass K Stanley}
Let $S_\Z^0$ be the subset of Grassmannian elements, the set of $w\in S_\Z$ such that $ws_i>w$ for all $i\in \Z\setminus \{0\}$. Let $\P$ be Young's lattice of partitions. There is a bijection $\P\to S_\Z^0$ denoted $\la\mapsto w_\la$ \cite{LLS:back stable}.

\begin{ex} For $\la=(3,2)$ we have $w_\la = s_0s_{-1}s_2s_1s_0$. This can be obtained by a row-reading word from the shape of $\la$ in which 
the reflection $s_{j-i}$ is placed in the box $(i,j)$ in the $i$-th row and $j$-th column in matrix-style indexing.
\begin{align*}
\begin{ytableau}
s_0&s_1&s_2 \\ s_{-1} & s_0
\end{ytableau}
\end{align*}
\end{ex}

Define the \emph{Grassmannian $K$-Stanley function}
\begin{align}\label{E:G Grass}
	G_\la := G_{w_\la} \qquad\text{for $\la\in \P$.}
\end{align}
Some examples of $G_\la$ are given in Appendix \ref{S:computations}.

\begin{lem}\label{L:G stable}
	For $\la \in \P$, we have
	\begin{align}\label{E:G stable Grassmannian}
		G_\la = \bG_{w_\la}.
	\end{align}
\end{lem}
\begin{proof} Since $w_{\la}\in S_{\Z}^0$, it follows that $\bG_{w_\la}$ is $S_{\ne0}$-symmetric, that is, $\bG_{w_\la}\in \hLa$. Thus $G_\la=G_{w_\la}=\eta_0(\bG_{w_\la})=\bG_{w_\la}$, as required.
\end{proof}


\begin{prop} \label{P:omega G partition} For $\la \in \P$, we have
	\begin{align}
		\label{E:omega G partition}
		\omk(G_\la) &= G_{\la'}. 
	\end{align}
\end{prop}
\begin{proof} We have $\omega(w_\la)=w_{\la'}$.  This is a special case of  Proposition \ref{P:omk backstable} for $w_\la$.
\end{proof}

The following is easily deduced from e.g. the Schur expansion of $G_\la$ in \cite[Thm. 2.2]{Len}.

\begin{prop} \label{P:G one row column} For all $r\ge1$, we have
	\begin{align}
	\label{E:row G to s}
		G_r &= \sum_{i\ge0} (-1)^r s_{r,1^i} \\
	\label{E:col G to s}
		G_{1^r} &= \sum_{i\ge0} (-1)^i \binom{i+r-1}{r-1} s_{1^{r+i}}.
	\end{align}
\end{prop}


\section{Back stable double Grothendieck polynomials}
Recall the Laurent polynomial rings $R(a)$ and $R=R(x)$ from \S \ref{SS:change of variable}. Define the $R(a)$-algebras 
\begin{align}
	R(x;a) &= R(x) \otimes R(a) \\
R(x;a)^+ &= R(x)^+ \otimes R(a)
\end{align}
For $w\in S_\Z$, let $w^x$ (resp. $w^a$) denote the action of $w$ on the $x$ (resp. $a$) variables in $R(x;a)$. We use similar superscript notation for other operators.

\subsection{Double Grothendieck polynomials}
For $w\in S_n$, the \emph{double Grothendieck polynomial} $\G_w\in R(x;a)$ is defined by
\begin{align}\label{E:G double}
	\G_{w_0}(x;a) &= \prod_{i+j\le n} (x_i \ominus a_j) \\
	\G_w(x;a) &= \Gdiff_i^x(\G_{ws_i}(x;a)) \qquad\text{if $ws_i>w$.}
\end{align}

\begin{ex} The double Grothendieck polynomials $\G_w(x;a)$ for $w\in S_3$ are given below.
	\begin{equation*}
		\begin{tikzcd}
			& (x_1\ominus a_1)(x_1\ominus a_2)(x_2\ominus a_1) \arrow[dl,swap,"\Gdiff_1"] \arrow[dr,"\Gdiff_2"]   \\
			(x_1\ominus a_1)(x_2\ominus a_1) \arrow[d,swap,"\Gdiff_2"] && (x_1\ominus a_1)(x_1\ominus a_2) \arrow[d,"\Gdiff_1"] \\[5pt]
			x_1\ominus a_1 \arrow[dr,swap,"\Gdiff_1"]&& (x_1\ominus a_1)\oplus (x_2 \ominus a_2)  \arrow[dl,"\Gdiff_2"] \\
			& 1
		\end{tikzcd}
	\end{equation*}
\end{ex}

\begin{rem} \label{R:double Groth convention}
The literature defines $\G_w(x;a)$ using $\oplus$ instead of $\ominus$. This has the convenient feature of avoiding denominators. However our convention is the most natural with respect to localization. 
\end{rem}

One may show that $\G_w(x;a)$ is well-defined for $w\in S_+$.

\begin{lem}\label{L:forget double to single} For $w \in S_+$, we have
	\begin{align}\label{E:forget double to single}
		\G_w(x;0) = \G_w(x).
	\end{align}
\end{lem}

\begin{prop} \label{P:double positive basis}
	The set $\{\G_w(x;a) \mid w\in S_+\}$ is a $R(a)$-basis of $R(a) \otimes \Z[x_i\mid i\in \Z_{>0}]$.
\end{prop}
\begin{proof} Follows by Lemma \ref{L:forget double to single} and Proposition \ref{P:G poly basis}.
\end{proof}

\subsection{Negative double Grothendieck polynomials}\label{SS:negative double G}
The automorphism $\omk$ of $R(x)$ defined by \eqref{E:omk x} can be extended to
a ring automorphism of $R(x;a)$ by letting $\omk(a_i) = \ominus(a_{1-i})$ for all $i\in \Z$. 

Recall that $\omega$ acts on $S_\Z$ by 	\eqref{E:omega s}.  For $w \in S_-$, define the \emph{negative double Grothendieck polynomials} by
\begin{align*}
	\G_w(x;a) := \omk(\G_{\omega(w)}(x;a)). 
\end{align*}

\begin{ex} We have $\G_{s_1}(x;a)=x_1\ominus a_1$.  Thus
	$$\G_{s_{-1}}(x;a) = \omk(\G_{s_1})=\omk(x_1\ominus a_1) = (\ominus x_0)\ominus(\ominus a_0) =
	a_0\ominus x_0.$$
\end{ex}

\begin{prop} \label{P:double negative basis}
	The set $\{\G_w(x;a) \mid w\in S_-\}$ is a $R(a)$-basis of 
	$R(a)[(1-x_i)^{-1}\mid i \in \Z_{\le0}]$.
\end{prop}
\begin{proof} This follows from the fact that $\omk$ restricts to a $R(a)$-algebra isomorphism
	\begin{align*}
		R(a)[x_i\mid i \in \Z_{>0}] &\to R(a)[(1-x_i)^{-1}\mid i \in\Z_{\le0}] \\
		\G_w(x;a)& \mapsto \G_{\omega(w)}(x;a)\qquad\text{for $w\in S_+$.} \qedhere
	\end{align*}	
\end{proof}

For $w\in S_{\ne0}$, write $w=uv$ with $u\in S_+$ and $v\in S_-$. Define
\begin{align*}
	\G_w(x;a) := \G_u(x;a) \G_v(x;a) \qquad\text{for $w\in S_{\ne0}$.}
\end{align*}

\begin{prop} \label{P: double nonzero basis}
The set $\{\G_w(x;a) \mid w\in S_{\ne0} \}$ is a $R(a)$-basis of  $R(x;a)^+$.
\end{prop}
\begin{proof} Follows immediately from Propositions \ref{P:double positive basis} and \ref{P:double negative basis}.
\end{proof}

\begin{prop} \label{P:G double to single} For $w\in S_{\ne0}$, we have
	\begin{align}\label{E:double to single}
		\G_w(x;a) &= \sum_{u*v=w} (-1)^{\ell(u)+\ell(v)-\ell(w)} \G_{u^{-1}}(\ominus a)\G_v(x), \\
		\label{E:double to single inverted}
		\G_w(x) &= \sum_{u*v=w} (-1)^{\ell(u)+\ell(v)-\ell(w)} \G_u(a)\G_v(x;a).
	\end{align}
\end{prop}
\begin{proof} For \eqref{E:double to single}, by factoring and applying $\omk$ one may reduce to the case that $w\in S_+$. Then $w\in S_n$ for some $n$, in which case \eqref{E:double to single} holds by \cite[Prop. 3.2, Lemma 5.5]{FK}.
	For \eqref{E:double to single inverted} apply Proposition \ref{P:left Groth  inversion}.
\end{proof}

\begin{ex} We have $\G_{s_1}(x;a) = x_1 \ominus a_1$.  Using \eqref{E:double to single} with the Hecke factorizations $u*v=s_1$ given by $(u,v)$ equal to $(s_1,\id)$, $(\id,s_1)$, and $(s_1,s_1)$ we have
\begin{align*}
	\G_{s_1}(x;a) &= \G_{s_1}(\ominus a)\G_{\id}(x) + \G_\id(\ominus(a))\G_{s_1}(x)- \G_{s_1}(\ominus a)\G_{s_1}(x) \\
	&= \G_{s_1}(\ominus a) \oplus \G_{s_1}(x) = \ominus(a_1) \oplus x_1 = x_1 \ominus a_1.
\end{align*}
\end{ex}

\begin{prop} \label{P:double G inverse} For $w\in S_{\ne0}$, we have
	\begin{align}\label{E:double G inverse}
		\G_{w^{-1}}(x;a) = \G_w(\ominus a;\ominus  x).
	\end{align}	
\end{prop}
\begin{proof} This follows from \eqref{E:double to single} using that $u*v = w$ if and only if $v^{-1} * u^{-1} = w^{-1}$.
\end{proof}

\begin{prop}\label{P:equivariant ddiff} We have
\begin{align}
	\Gdiff_i^{a,\anti} \G_w(x;a) = \begin{cases}
		\G_{s_iw}(x;a) & \text{if $s_i w<w$} \\
		\G_w(x;a) & \text{otherwise.}
		\end{cases}
\end{align}
\end{prop}
\begin{proof} Follows by Proposition \ref{P:double G inverse} and Lemma \ref{L:ddiff and negation}.
\end{proof}

\subsection{Double symmetric function ring}\label{SS:double symmetric}
Let $\SSS=\Z[[a_i\mid i\in\Z]]$ be the ring of formal power series in variables $a_i$ for $i\in\Z$ with coefficients in $\Z$. Let $\hLa(x||a) := \hLa(x/a) \otimes_\Z \otimes \SSS$ be the $\SSS$-algebra of formal power series in $e_k(x/a)$ with coefficients in $\SSS$.  The ring $\hLa(x||a)$ is a Hopf algebra over $\SSS$ with coproduct such that the elements $p_k(x||a)$ are primitive,
counit $\epsilon_a: \hLa(x||a) \to \SSS$ given by $\epsilon_a(p_k(x||a))=0$ for $k\ge1$, and antipode $S(p_k(x||a)) = - p_k(x||a)$ for all $k\ge1$. 

\subsection{Double back stable rings}
Define the $\SSS$-algebras, called \emph{double back stable rings},
\begin{align}
\label{E:BR}
\bR(x;a) &= 
\hLa(x||a) \otimes_{R(a)} R(x;a)  \\
\label{E:BR+}
\bR(x;a)^+ &= \hLa(x||a) \otimes_{R(a)} R(x;a)^+.
\end{align}
The infinite symmetric group $S_\Z$ has two commuting actions on $\bR(x;a)$, one acting on $x$ variables and the other on $a$ variables,
including ``the variables in $\hLa(x||a)$", where $p_k(x||a)$ is as in \S \ref{SS:coproduct}:
\begin{align*}
	s_i^x (p_k(x||a)) &= p_k(x||a) + \delta_{i0} (x_1^k-x_0^k) \\
	s_i^a (p_k(x||a)) &= p_k(x||a) + \delta_{i0} (a_0^k-a_1^k).
\end{align*}
The algebra $\hLa(x||a)$ is the $\SSS$-subalgebra of  $S^x_{\ne0}$-invariants in $\bR(x;a)$.

\begin{rem} Since formal series in $p_k(x||a)$ are allowed, in order to admit an $S_\Z$-action, series in the $a_i$ must also be allowed. This is why $\SSS$ is used for the coefficient ring rather than $R(a)$.
\end{rem}

\subsection{$\Aut(A_\Z)$-action on double back stable rings}
The group $\Aut(A_\Z)$ of Dynkin automorphisms acts on $\bR(x;a)$ by ring automorphisms.
Define $\gamma:\bR(x;a) \to \bR(x;a)$ by 
$$\gamma(x_i)=x_{i+1}, \qquad \gamma(a_i)=a_{i+1}, \qquad \gamma(p_k(x||a)) = p_k(x||a) + x_1^k-a_1^k.$$ Let $\omk$ be the ring automorphism of $\bR(x;a)$ extending the automorphism $\omk$ of $R(x;a)$ in 
\S \ref{SS:negative double G} by
\begin{align}\label{E:omk power double}
	\omk(p_k(x||a)) = (-1)^{k+1}  \sum_{r\ge0} \binom{r+k-1}{k-1} p_{k+r}(x||a) \qquad\text{for all $k\ge 1$.}
\end{align}
This is consistent with the definition of $\omk(p_k(\xm))$ and the parallel definition of $\omk(p_k(\am))$ in \S \ref{SS:sym func} (using the convention $p_k(\ap)=-p_k(\am)$). 

\begin{prop}\label{P:omk shift}
The maps	 $\gamma$ and $\omk$ define an action of $\Aut(A_\Z)$ on $\bR(x;a)$. That is,
	\begin{align}\label{E:omk shift}
	\omk\circ	\gamma \circ \omk = \gamma^{-1}.
	\end{align}
\end{prop}
\begin{proof} This can be readily verified on the generators.
\end{proof}

\begin{prop} \label{P:omk G super}
	For all $w\in S_\Z$, we have
	\begin{align}
		\omk(G_w(x/a)) &= G_{\omega(w)}(x/a).
	\end{align}
\end{prop}
\begin{proof} 
Follows from Proposition \ref{P:omk double backstable} and Lemma \ref{L:double BG to super K Stanley}.
\end{proof}

\subsection{Back stable double Grothendieck polynomials}
\label{SS:double bG}
Given $w\in S_{\Z}$, let $[p,q]\subset\Z$ be an interval that contains all elements of $\Z$ moved by $w$. Define $\G_w^{[p,q]}(x;a)\in R(x;a)$ by
\begin{align*}
	\G_w^{[p,q]}(x;a) = \gamma^{p-1}(\G_{\gamma^{1-p}(w)}(x;a))
\end{align*}

Define the \emph{back stable double Grothendieck polynomial} $\bG_w(x;a)$ by
\begin{align*}
  \bG_w(x;a) = \lim_{\substack{p\to-\infty \\ q\to\infty}} \G_w^{[p,q]}(x;a).
\end{align*}

It is immediate that we have
\begin{align}
	\Gdiff_i^{x} \G_w(x;a) = \begin{cases}
		\G_{ws_i}(x;a) & \text{if $ws_i<w$} \\
		\G_w(x;a) & \text{otherwise.}
		\end{cases}
\end{align}

\begin{prop}\label{P:equiv ddiff bG} We have
\begin{align}\label{E:equiv ddiff bG}
	\Gdiff_i^{a,\anti}(\bG_w(x;a)) = 
	\begin{cases}
		\bG_{s_iw}(x;a) & \text{if $s_iw<w$} \\
		\bG_w(x;a) &\text{otherwise.}
	\end{cases}
\end{align}
\end{prop}
\begin{proof} Follows from Proposition \ref{P:equivariant ddiff}.
\end{proof}

\begin{prop}\label{P:G back stable double} For $w \in S_\Z$, we have
\begin{align}\label{E:G back stable double}
  \bG_w(x;a) = \sum_{u*v=w} (-1)^{\ell(u)+\ell(v)-\ell(w)} \bG_{u^{-1}}(\ominus(a)) \bG_v(x)
\end{align}
\end{prop}
\begin{proof} This holds by Proposition \ref{P:G double to single} and the definitions.
\end{proof}

\begin{prop}\label{P:G back stable triple}  For $w \in S_\Z$, we have
\begin{align}\label{E:G back stable triple}
	\bG_w(x;a) = \sum_{\substack{u*v*z=w \\ u,z\in S_{\ne0}}} (-1)^{\ell(u)+\ell(v)+\ell(z)-\ell(w)} \G_{u^{-1}}(\ominus a) G_v(x/a) \G_z(x)
\end{align}
\end{prop}
\begin{proof} Let $w=u*v$, $u^{-1} = u_1*v_1$ and $v=u_2*v_2$ with $v_1,v_2\in S_{\ne0}$, so that $w=v_1^{-1}*u_1^{-1}*u_2*v_2$.
By Propositions \ref{P:G back stable double}, \ref{T:coproduct}, and \ref{P:G super} we have
\begin{align*}
  \bG_w(x;a) &= \sum_{u*v=w} (-1)^{\ell(u)+\ell(v)-\ell(w)} 
  \bG_{u^{-1}}(\ominus(a)) \bG_v(x) \\
  &= \sum_{\substack{v_1^{-1}*u_1^{-1}*u_2*v_2=w\\
  		v_1,v_2\in S_{\ne0}}}
  (-1)^{\ell(u_1)+\ell(v_1)+\ell(u_2)+\ell(v_2)-\ell(w)} \\
  &\quad
   \G_{v_1}(\ominus(a)) G_{u_1}(\ominus(a)) G_{u_2}(x) \G_{v_2}(x) \\
  &= \sum_{\substack{v_1^{-1}*y*v_2=w\\
  		v_1,v_2\in S_{\ne0}}} (-1)^{\ell(v_1)+\ell(y)+\ell(v_2)-\ell(w)} \G_{v_1}(\ominus(a)) G_y(x/a) \G_{v_2}(x). \qedhere
\end{align*}
\end{proof}

\begin{ex}\label{ex:doublecoproduct} We have
\begin{align*}
	\bG_{s_0}(x;a) &= G_{s_0}(x/a) = G_1(x/a) \\
	\bG_{s_1}(x;a) &= \G_{s_1}(\ominus a) + G_{s_1}(x/a) + \G_{s_1}(x) \\
	&- \G_{s_1}(\ominus a) G_{s_1}(x/a) - \G_{s_1}(\ominus a) \G_{s_1}(x) - G_{s_1}(x/a) \G_{s_1}(x) \\
	&+ \G_{s_1}(\ominus a) G_{s_1}(x/a) \G_{s_1}(x) \\
	&= \G_{s_1}(\ominus a) \oplus G_1(x/a) \oplus \G_{s_1}(x) \\
	&= \ominus(a_1) \oplus G_1(x/a) \oplus  x_1 = G_1(x/a) \oplus \G_{s_1}(x;a).
\end{align*}
\end{ex}

\begin{prop}For $w\in S_\Z$, we have
	\begin{align}\label{E:double bG inverse}
		\bG_{w^{-1}}(x;a) = \bG_w(\ominus a;\ominus  x).
	\end{align}	
\end{prop}
\begin{proof} This holds by Proposition \ref{P:G back stable double} and the fact that $u*v=w$ if and only if $v^{-1}*u^{-1}=w^{-1}$.
\end{proof}

\begin{lem} \label{L:double BG to super K Stanley}
	For all $w\in S_\Z$, we have
	\begin{align}\label{E:double BG to super K Stanley}
		G_w(x/a) = \bG_w(x;a)|_{x_i\mapsto0, a_i\mapsto0}.
	\end{align}
\end{lem}
\begin{proof} This follows by Proposition \ref{P:G back stable triple} and the fact that $\G_z(x)|_{x_i\mapsto0}=0$ for $z\ne\id$.
\end{proof}

Define the forgetful ring homomorphism $\For:\bR(x;a) \to \bR$ by $\For(p_k(x||a)) = p_k(\xm)$, $\For(x_i)=x_i$, and $\For(a_i)=0$.

\begin{prop} \label{P:forget Schubert basis} For all $w\in S_\Z$, $\For(\bG_w(x;a))=\bG_w(x)$.
\end{prop}
\begin{proof} Follows from Propositions \ref{P:G back stable triple} and \ref{T:coproduct}.
\end{proof}

\subsection{$\Aut(A_\Z)$-action on back stable double Grothendieck polynomials}

$\Aut(A_\Z)$ permutes the back stable double Grothendieck polynomials. By Theorem~\ref{T:infinite basis} and Proposition~\ref{P:Aut GKM} the following hold.


\begin{prop}\label{P:double bG shift} For $w \in S_\Z$, we have
\begin{align}\label{E:double bG shift}
 \bG_{\gamma(w)}(x;a) =	\gamma(\bG_w(x;a)).
\end{align}
\end{prop}

\begin{prop}\label{P:omk double backstable} For $w \in S_\Z$, we have
\begin{align}
	\bG_{\omega(w)}(x;a) = \omk(\bG_w(x;a)). 
\end{align}
\end{prop}

\subsection{Back stable double Grothendieck polynomials are the equivariant Schubert basis}

For $f\in \bR(x;a)$ and $w\in S_\Z$ define $f|_w = f(wa;a) = \epsilon_a(w^x(f(x;a))$. The map $f\mapsto f|w$ is an $\SSS$-algebra homomorphism $\bR(x;a) \to \SSS$.

\begin{prop} \label{P:locs} \
	\begin{enumerate}
		\item 	For $w,y\in S_+$, $\bG_w(x;a)|_y = \G_w(x;a)|_y$.
		\item For $w,y\in S_\Z$ let $M\in\Z_{\ge0}$ be such that $\gamma^M(w),\gamma^M(y)\in S_+$.
		Then 
		\begin{align*}
			\bG_w(x;a)|_y= \gamma^{-M} \G_{\gamma^M(w)}(x;a)|_{\gamma^M(y)}.
		\end{align*}
	\end{enumerate}
\end{prop}
\begin{proof} For (1) consider \eqref{E:G back stable triple}. Since $y\in S_+$, $G_v(x/a)$ is invariant under $y^a$ and vanishes at $x=ya$ unless $v=\id$. Therefore by \eqref{E:double to single} we have
	\begin{align*}
		\bG_w (x;a)|_y &= \sum_{\substack{u*z=y \\ u,z\in S_{\ne0}}} (-1)^{\ell(u)+\ell(z)-\ell(w)}
		\G_{u^{-1}}(\ominus a) \G_z(ya)\\
		&= \G_w(x;a)|_y.
	\end{align*}
	Part (2) follows from part (1) and \eqref{E:double bG shift}.
\end{proof}


\begin{thm} \label{T:infinite basis}
	There is a $\SSS$-algebra and $S_\Z$-equivariant embedding $\res: \bR(x;a) \to \SSS\otimes_{R(a)} \Psit$ defined by $\res(f)(w) = f|_w$ for all $w\in S_\Z$. Moreover,
		$\res(\bG_v)=\psit^v$ for all $v\in S_\Z$.
\end{thm}
\begin{proof} We observe that the values of $f|_w$ are in $\SSS$. Let $f\in \bR(x;a)$. $w^x(f) - (s_{ij}w)^x(f)$ is $s_{ij}^x$-antisymmetric and therefore divisible by $x_i-x_j$. It follows that $f|_w - f|_{s_{ij}w}$ is a multiple of $a_i-a_j$. Thus $\res(f) \in \SSS\otimes_{R(a)} \Psit$ by \eqref{E:transformed GKM}.	

Next we show that $\res(\bG_v(x;a)) = \psit^v$ for all $v\in S_\Z$. By Proposition \ref{P:Schubert loc}, this is equivalent to showing that
\begin{align}\label{E:initial}
\bG_v(a;a)=\delta_{\id,v}\qquad\text{for $v\in S_\Z$} 
\end{align}
 and if $ws_i<w$ then
\begin{align}\label{E:recursive}
	\bG_v(wa;a) = \begin{cases}
		\bG_v(ws_ia;a) & \text{if $vs_i>v$} \\
		\dfrac{1-a_{w(i)}}{1-a_{w(i+1)}} \bG_v(ws_ia;a) + \dfrac{a_{w(i)}-a_{w(i+1)}}{1-a_{w(i+1)}}  \bG_{vs_i}(wa;a)  & \text{if $vs_i<v$.}
	\end{cases}
\end{align}
Equation \eqref{E:initial} follows by Proposition \ref{P:locs} and \eqref{E:cancellation}. 

Suppose $w\ne\id$. Let $i\in I$ be such that $ws_i<w$. Suppose first that $vs_i>v$. Then
$\bG_v(x;a) = \Gdiff_i^x \bG_v(x;a)$. Since the image of $\Gdiff_i^x$ is $s_i^x$-invariant, we deduce that
$\bG_v(wa;a)=\bG_v(ws_ia;a)$ as required. 

Suppose $vs_i<v$. Then $\bG_{vs_i}(x;a) = \Gdiff_i^x \bG_v(x;a)$. This yields
\begin{align*}
	\bG_{vs_i}(wa;a) &= (a_{w(i)}-a_{w(i+1)})^{-1} \left( (1-a_{w(i+1)}) \bG_v(wa;a) - (1-a_{w(i)})\bG_v(ws_ia;a)
	 \right)
\end{align*}
which rearranges to \eqref{E:recursive}, as required.
\end{proof}

\begin{cor} \label{C:Schub basis}
The back stable Grothendieck polynomials $\bG_w$ represent the Schubert classes in the $K$-theory $K(\Flb)$.  The back stable double Grothendieck polynomials $\bG_w(x;a)$ represent the Schubert classes in the torus-equivariant $K$-theory $K_T(\Flb)$.
\end{cor}

\subsection{Double $K$-Stanley functions}\label{sec:doubleK}
Let $\eta_a: \bR(x;a)\to \hLa(x||a)$ be the $\SSS$-algebra homomorphism sending
$x_i\mapsto a_i$ and $p_r(x||a)\mapsto p_r(x||a)$. Define the \emph{double $K$-Stanley function} $G_w(x||a)\in\hLa(x||a)$ by
\begin{align}\label{E:Gw definition}
	G_w(x||a) := \eta_a(\bG_w(x;a)).
\end{align}

\begin{prop}\label{P:G stable double} For $w \in S_\Z$, we have
	\begin{align}\label{E:G stable double}
		G_w(x||a) = \sum_{\substack{u*v*z=w \\ u,z\in S_{\ne0}}} (-1)^{\ell(u)+\ell(v)+\ell(z)-\ell(w)} \G_{u^{-1}}(\ominus(a)) G_v(x/a) \G_z(a)
	\end{align}
\end{prop}
\begin{proof} This follows immediately from Proposition \ref{P:G back stable triple}.
\end{proof}

\begin{ex} Since $G_{s_0s_1}=G_{s_{-1}s_0}=G_{11}$ and $G_{s_0}=G_1$ we have
	\begin{align*}
		G_{s_0s_1}(x||a) &= G_{s_0s_1}(x/a) + G_{s_0}(x/a) \G_{s_1}(a) - G_{s_0s_1}(x/a) \G_{s_1}(a) \\
		&= G_{s_0s_1}(x/a) + G_{s_0}(x/a) a_1 - G_{s_0s_1}(x/a) a_1  \\
	 &= (1-a_1) G_{11}(x/a) + a_1 G_1(x/a).
	\end{align*}
Since $\G_{s_{-1}}(\ominus a) = a_0$ we have
	\begin{align*}
		G_{s_{-1}s_0}(x||a) &= \G_{s_{-1}}(\ominus a) G_{s_0}(x/a) + G_{s_{-1}s_0}(x/a) -\G_{s_{-1}}(\ominus a) G_{s_{-1}s_0}(x/a) \\
		&=(1-a_0) G_{11}(x/a) +  a_0 G_1(x/a).
\end{align*}
\end{ex}

\begin{prop}\label{P:double G negation} For $w \in S_\Z$, we have
\begin{align*}
	\omk(G_w(x||a)) = G_{\omega(w)}(x||a).
\end{align*}
\end{prop}
\begin{proof} 
It is straightforward to check that $\eta_a \circ \omk = \omk$.  Thus \eqref{E:omega Gxa perm} follows from Proposition~\ref{P:double bG shift}.
\end{proof}

\begin{prop}\label{P:omega Gxa} For all $w\in S_\Z$, we have
\begin{equation}
\label{E:omega Gxa perm}
G_{w^{-1}}(x||a) = G_w(x||a)|_{x\mapsto \ominus a, a \mapsto \ominus x}.
\end{equation}
\end{prop}
\begin{proof}
The transformation $x\mapsto \ominus a, a \mapsto \ominus x$ commutes with $\eta_a$.  The result follows from \eqref{E:double bG inverse}.
\end{proof}

\subsection{Coproduct formula}
Recall that the Hopf algebra structure on $\hLa(x||a)$ is defined by letting $p_k(x||a)$ be primitive for all $k\ge 1$. We give $\bR(x;a)$ the structure of a $\hLa(x||a)$-comodule by letting $\Delta$ act on the tensor factor $\hLa(x||a)$ in $\bR(x;a)$.

\begin{thm}\label{T:double coproduct} For $w \in S_\Z$, we have
	\begin{align}	
		\label{E:double G coproduct}
		\Delta(\bG_w(x;a)) &= \sum_{u*v=w} (-1)^{\ell(u)+\ell(v)-\ell(w)} G_u(x||a) \otimes \bG_v(x;a) \\
		\label{E:double G little coproduct}
		\bG_w(x;a) &= \sum_{\substack{u*v=w\\ v\in S_{\ne0}}} (-1)^{\ell(u)+\ell(v)-\ell(w)} G_u(x||a) \G_v(x;a)
	\end{align}
\end{thm}
\begin{proof} Equation \eqref{E:double G little coproduct} follows from Propositions \ref{P:G back stable triple}, \ref{P:G stable double}, and \eqref{E:double to single inverted}.
	
	For \eqref{E:double G coproduct}, using Propositions \ref{P:G back stable triple}, \ref{P:G stable double}, and \ref{P:right Groth inversion} we have
	\begin{align*}
		\sum_{u*v=w} \pm G_u(x||a) \otimes \bG_v(x;a) 
		&= \sum_{\substack{u_1*v_1*z_1*u_2*v_2*z_2=w \\ u_i,z_i\in S_{\ne0}}} \pm \G_{{u_1}^{-1}}(\ominus(a)) G_{v_1}(x/a) G_{z_1}(a) \\
		&\quad \otimes
		\G_{u_2}^{-1}(\ominus(a)) G_{v_2}(x/a) G_{z_2}(x) \\
		&= \sum_{\substack{u_1*v_1*v_2*z_2=w \\ u_1,z_2\in S_{\ne0}}} \pm \G_{{u_1}^{-1}}(\ominus(a)) G_{v_1}(x/a)  \\
		&\quad \otimes
		G_{v_2}(x/a) G_{z_2}(x) \\
		&= \sum_{\substack{u_1*v*z_2=w \\ u_1,z_2\in S_{\ne0}}}
		\pm \G_{{u_1}^{-1}}(\ominus(a)) \Delta(G_v(x/a)) \G_{z_2}(x) \\
		&= \Delta(\bG_{w}(x;a)). \qedhere
	\end{align*}
\end{proof}

\subsection{Grassmannian double $K$-Stanley functions}
Recalling $w_\la$ from \S \ref{SS:Grass K Stanley}, define
\begin{align}\label{E:Grassmannian double K Stanley}
G_\la(x||a) := G_{w_\la}(x;a)\qquad\text{for $\la\in\P$.}
\end{align}

\begin{lem}\label{L:double G stable}
	For $\la \in \P$, we have
	\begin{align}
		G_\la(x||a) = \bG_{w_\la}(x;a).
	\end{align}
\end{lem}
\begin{proof}
For any Hecke factorization $w_\la =u*v*z$ with $u,z \in S_{\neq 0}$, we have $z = \id$.  The result then follows from the definition \eqref{E:Gw definition} by comparing \eqref{E:G back stable triple} and \eqref{E:G stable double}.
\end{proof}

\begin{cor} \label{C:Grass Schubert basis}
The Grassmannian $K$-Stanley functions $\{G_\la\mid \la\in\Y\}$ represent the basis of structure sheaves of opposite Schubert varieties in the $K$-theory $K(\Grb^0)$. The Grassmannian double $K$-Stanley functions $\{G_\la(x||a)\mid \la\in\Y\}$ represent the structure sheaves of opposite Schubert varieties in the torus-equivariant $K$-theory $K_T(\Grb^0)$.
\end{cor}
\begin{proof} Since $K_T(\Grb^0) = (K_T(\Flb))^{S^x_{\ne0}}$, the statement about $G_\la(x||a)$ follows from Corollary \ref{C:Schub basis}, Lemma \ref{L:double G stable}, and the fact that the $G_\la(x||a)$ are precisely the $S^x_{\ne0}$-invariant Schubert basis elements for $K_T(\Flb)$. Applying the forgetful homomorphism 
$K_T(\Grb) \to K(\Grb)$, the equivariant Schubert basis is sent to the Schubert basis. The proof is completed by applying Proposition \ref{P:forget Schubert basis}.
\end{proof}

\begin{prop}\label{P:double G negation Grassmannian} For $\la \in \P$, we have
\begin{align*}
	\omk(G_\la(x||a)) = G_{\la'}(x||a).
\end{align*}
\end{prop}
\begin{proof} Follows from Proposition \ref{P:double G negation} and the fact that
	$\omega(w_\la) = w_{\la'}$.
\end{proof}

A rook strip is a skew shape $\nu/\mu$ which has at most one box in each row and in each column.
Write $\nu/\mu\in\RS$ if $\nu/\mu$ is a rook strip. Let $d(\la)$ be the Durfee square of $\la$, the maximum
$d$ such that the $d\times d$ square is contained in $\la$.

If $\mu\subset \la$, we define 
\begin{align}\label{E:w skew}
	w_{\la/\mu} := w_\la w_\mu^{-1}.
\end{align}

\begin{lem} \label{L:Grassmannian Hecke factorization}
	Let $\mu\subset\la\in \P$ and $u\in S_{\ne0}$. Then $u* w_\mu = w_\la$ if and only if 
	$u=w_{\la/\nu}$ where $\nu\subset\mu$, $d(\nu)=d(\la)$, and $\mu/\nu\in\RS$. 
\end{lem}

The coproduct formula gives the transition matrix between the Grassmannian double $K$-Stanley functions and the super $K$-Stanley symmetric functions.

\begin{prop}\label{P:double to super} For $\lambda \in \P$, we have
\begin{align}
\label{E:double to super}
   G_\la(x||a) &= \sum_{\substack{\nu \subset \mu \subset \la \\ d(\nu)=d(\la)\\ \mu/\nu\in \RS}} (-1)^{|\mu|-|\nu|} \G_{w_{\la/\nu}^{-1}}(\ominus a) G_\mu(x/a) \\
 \label{E:super to double}
   G_\la(x/a) &= \sum_{\substack{\nu \subset \mu \subset \la \\ d(\nu)=d(\la)\\ \mu/\nu\in \RS}} (-1)^{|\mu|-|\nu|} \G_{w_{\la/\nu}}(a) G_\mu(x||a).
\end{align}
\end{prop}
\begin{proof} Equation \eqref{E:super to double} follows from \eqref{E:double to super} using 
Proposition \ref{P:left Groth inversion}. Equation \eqref{E:double to super} is an instance of 
Proposition \ref{P:G back stable triple}. In this application $z=\id$ since $w_\la\in S_\Z^0$. The Proposition follows from Lemma \ref{L:Grassmannian Hecke factorization}.
\end{proof}
In \S\ref{S:detformula}, we give determinantal formulae for $G_\la(x||a)$.

\begin{ex}\label{ex:G double to super}
    By Proposition~\ref{P:double to super} and using that $\G_{s_1}(a) = a_1$, $\G_{s_{-1}}(a) = \ominus(a_0)$, and $\G_{s_1s_{-1}}(a) = a_1 (\ominus(a_0))$, we have
    \begin{align*}
        G_2(x/a) &= (1-\G_{s_1}(a))G_2(x||a) + \G_{s_1}(a) G_1(x||a) = (1 - a_1) G_2(x||a) + a_1 G_1(x||a) \\
        G_{11}(x/a) &= (1-\G_{s_{-1}}(a))G_{11}(x||a) + \G_{s_{-1}}(a) G_1(x||a) = (1-\ominus(a_0)) G_{11}(x||a) + \ominus(a_0)G_1(x||a) \\
        G_{21}(x/a) &= (1-\G_{s_1}(a)-\G_{s_{-1}}(a)+\G_{s_1s_{-1}}(a))G_{21}(x||a) + (\G_{s_{-1}}(a)-\G_{s_1s_{-1}}(a)) G_{2}(x||a) \\ & + (\G_{s_{1}}(a)-\G_{s_1s_{-1}}(a)) G_{11}(x||a)+ \G_{s_1s_{-1}}(a)G_1(x||a) \\
        &= (1-(a_1\ominus a_0)) G_{21}(x||a) +\ominus(a_0)(1-a_1) G_2(x||a) + a_1(1-\ominus(a_0)) G_{11}(x||a) \\
        &+ a_1(\ominus(a_0)) G_1(x||a).
    \end{align*}
\end{ex}

\section{$K$-Bumpless pipedreams}
In \cite{LLS:back stable}, we introduced bumpless pipedreams and showed that back stable (double) Schubert polynomials can be obtained as sums over bumpless pipedreams.  Weigandt \cite{Wei} connected bumpless pipedreams to alternating sign matrices and a formula of Lascoux \cite{Las}, and thereby obtained a bumpless pipedream formula for (double) Grothendieck polynomials.  

\subsection{Back stable double Grothendieck polynomials}
Recall that a bumpless pipedream is a tiling of the plane by the tiles: empty, NW elbow, SE elbow, horizontal line, crossing, and vertical line.
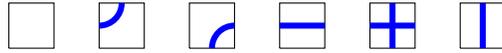
\begin{center}
\begin{tikzpicture}[scale=0.6,line width=0.8mm]
\bbox{-3}{0}
\bbox{-1}{0}
\leftelbow{-1}{0}{blue}
\bbox{1}{0}
\rightelbow{1}{0}{blue}
\bbox{3}{0}
\horline{3}{0}{blue}
\bbox{5}{0}
\cross{5}{0}{blue}{blue}
\bbox{7}{0}
\vertline{7}{0}{blue}
\end{tikzpicture}
 \end{center}
 We shall use {\bf matrix coordinates} for unit squares in the plane.  Thus row coordinates increase from top to bottom, column coordinates increase from left to right, and $(i,j)$ indicates the square in row $i$ and column $j$.  
 A \emph{$K$-bumpless pipedream} is a bumpless pipedream $D$ covering the whole plane, such that for all $N \gg 0$ and all $N \ll 0$, there is a pipe traveling north from $(\infty,N)$ to the square $(N,N)$ where it turns east and travels towards $(N,\infty)$.  The permutation $w(D) \in S_\Z$ of a $K$-bumpless pipedream is obtained as follows.  For each $i \in \Z$, there is a pipe that heads north from $(\infty, i)$.  We follow this pipe until it heads east towards $(j, \infty)$, ignoring all crossings between pairs of pipes that have already crossed (reading pipes from SW to NE).  Then $w(D)$ is determined by $w(j) = i$, as $i \in \Z$ varies.
 
The weight $\wt(D)$ of a $K$-bumpless pipedream $D$ is given by
\begin{equation}\label{eq:Kwt}
\wt(D) := \prod_{{\rm empty} \; {\rm tiles} \;(i,j)} (-(x_i \ominus a_j)) \prod_{{\rm NW-elbows} \;(i,j)} (1- (x_i \ominus a_j)),
\end{equation}
where the first product is over empty tiles $(i,j)$ and the second product is over elbow tiles $(i,j)$ that connect the north and west sides.

\begin{ex} Let $w = s_0s_2$.  In one line notation, $w(-2,-1,0,1,2,3) = (-2,-1,1,0,3,2)$ and the rest are fixed points.  Figure~\ref{fig:bumplessexample} shows a $K$-bumpless pipedream $D$ for $w$, where we have only drawn the region $\{(i,j) \mid i,j \in [-2,3]\}$ (the rest of the pipes head north, turn once and head east).  In the left picture, the empty tiles have been indicated, as have the row and column numbers.  In the right picture, we have indicated the calculation of $w(D)$, labeling each pipe by the column where it enters the picture.  The pipes labeled $0$ and $1$ intersect twice, and the second intersection is ignored when computing $w(D)$.
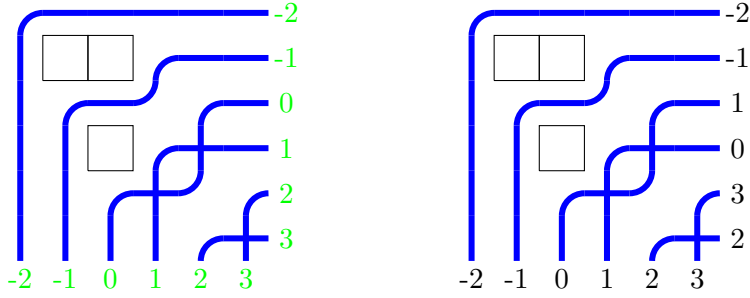
\begin{figure}
\centering
\begin{tikzpicture}[scale=0.6,line width=0.8mm]
\rightelbow{-2}{2}{blue}
\horline{-1}{2}{blue}
\horline{0}{2}{blue}
\horline{1}{2}{blue}
\horline{2}{2}{blue}
\horline{3}{2}{blue}

\vertline{-2}{1}{blue}
\bbox{-1}{1}{blue}
\bbox{0}{1}{blue}
\rightelbow{1}{1}{blue}
\horline{2}{1}{blue}
\horline{3}{1}{blue}

\vertline{-2}{0}{blue}
\rightelbow{-1}{0}{blue}
\horline{0}{0}{blue}
\leftelbow{1}{0}{blue}
\rightelbow{2}{0}{blue}
\horline{3}{0}{blue}

\vertline{-2}{-1}{blue}
\vertline{-1}{-1}{blue}
\bbox{0}{-1}{blue}
\rightelbow{1}{-1}{blue}
\cross{2}{-1}{blue}{blue}
\horline{3}{-1}{blue}

\vertline{-2}{-2}{blue}
\vertline{-1}{-2}{blue}
\rightelbow{0}{-2}{blue}
\cross{1}{-2}{blue}{blue}
\leftelbow{2}{-2}{blue}
\rightelbow{3}{-2}{blue}

\vertline{-2}{-3}{blue}
\vertline{-1}{-3}{blue}
\vertline{0}{-3}{blue}
\vertline{1}{-3}{blue}
\rightelbow{2}{-3}{blue}
\cross{3}{-3}{blue}{blue}

%

\draw[green] (-1.5,-3.4) node {-2};
\draw[green] (-0.5,-3.4) node {-1};
\draw[green] (0.5,-3.4) node {0};
\draw[green] (1.5,-3.4) node  {$1$};
\draw[green] (2.5,-3.4) node  {$2$};
\draw[green] (3.5,-3.4) node  {3};

\draw[green] (4.4,2.5) node {-2};
\draw[green] (4.4,1.5) node {-1};
\draw[green] (4.4,.5) node {0};
\draw[green] (4.4,-.5) node {1};
\draw[green] (4.4,-1.5) node {2};
\draw[green] (4.4,-2.5) node {3};


\begin{scope}[shift={(10,0)}]
\rightelbow{-2}{2}{blue}
\horline{-1}{2}{blue}
\horline{0}{2}{blue}
\horline{1}{2}{blue}
\horline{2}{2}{blue}
\horline{3}{2}{blue}

\vertline{-2}{1}{blue}
\bbox{-1}{1}{blue}
\bbox{0}{1}{blue}
\rightelbow{1}{1}{blue}
\horline{2}{1}{blue}
\horline{3}{1}{blue}

\vertline{-2}{0}{blue}
\rightelbow{-1}{0}{blue}
\horline{0}{0}{blue}
\leftelbow{1}{0}{blue}
\rightelbow{2}{0}{blue}
\horline{3}{0}{blue}

\vertline{-2}{-1}{blue}
\vertline{-1}{-1}{blue}
\bbox{0}{-1}{blue}
\rightelbow{1}{-1}{blue}
\cross{2}{-1}{blue}{blue}
\horline{3}{-1}{blue}

\vertline{-2}{-2}{blue}
\vertline{-1}{-2}{blue}
\rightelbow{0}{-2}{blue}
\cross{1}{-2}{blue}{blue}
\leftelbow{2}{-2}{blue}
\rightelbow{3}{-2}{blue}

\vertline{-2}{-3}{blue}
\vertline{-1}{-3}{blue}
\vertline{0}{-3}{blue}
\vertline{1}{-3}{blue}
\rightelbow{2}{-3}{blue}
\cross{3}{-3}{blue}{blue}

\draw(-1.5,-3.4) node {-2};
\draw (-0.5,-3.4) node {-1};
\draw (0.5,-3.4) node {0};
\draw(1.5,-3.4) node  {$1$};
\draw (2.5,-3.4) node  {$2$};
\draw (3.5,-3.4) node  {3};

\draw (4.4,2.5) node {-2};
\draw (4.4,1.5) node {-1};
\draw (4.4,.5) node {1};
\draw (4.4,-.5) node {0};
\draw (4.4,-1.5) node {3};
\draw (4.4,-2.5) node {2};
\end{scope}
\end{tikzpicture}
\caption{A $K$-bumpless pipedream for $w = s_0 s_2 = \ldots,-2,-1,1,0,3,2,\ldots$ with weight $\wt = -(x_{-1} \ominus a_{-1})(x_{-1} \ominus a_0)(x_1 \ominus a_0)(1- (x_0\ominus a_1) (1-(x_2 \ominus a_2))$.}
\label{fig:bumplessexample}
\end{figure}

\end{ex}

The following result follows from \cite[Theorem 1.1]{Wei}, reproduced as Theorem~\ref{thm:Wei} below.
\begin{thm}\label{thm:Kbumpless} Let $w \in S_\Z$.  Then
$$
\bG_w(x;a) = (-1)^{\ell(w)} \sum_D \wt(D)
$$
where the summation is over all $K$-bumpless pipedreams $D$ with permutation $w(D) = w$.
\end{thm}
\begin{ex} Let $w = s_0$. Then for each $j \geq 0$, there is one $K$-bumpless pipedream $D_j$ for $w$ with one empty tile $(-j,-j)$ and $j$ NW-elbow tiles $(-k,-k)$ for $0\le k < j$ and 
$$\wt(D_j) = -(x_{-j} \ominus a_{-j}) \prod_{ 0 \leq k < j} (1-(x_{-k} \ominus a_{-k})).$$
By Proposition~\ref{P:double to super}, we have $\bG_{s_0}(x;a) = G_1(x/a)$ and by Proposition~\ref{P:G super}, we have $G_1(x/a) = G_1(x) + G_1(\ominus a) - G_1(x) G_1(\ominus a) = 1 - \prod_{j\ge0} (1-x_{-j})/(1-a_{-j})$.
Using this, one checks that indeed $\bG_{s_0}(x;a) = \sum_{j \geq 0} (x_{-j} \ominus a_{-j}) \prod_{ 0 \leq k < j} (1-(x_{-k} \ominus a_{-k}))$.
\end{ex}

\subsection{Pipedreams for Grassmannian double $K$-Stanley functions}
Let $\la \in \Y$.  Recall from \cite{LLS:back stable} that a \emph{$\la$-halfplane pipedream} is a bumpless pipedream in the upper halfplane $ \Z_{\leq 0} \times \Z$ such that the crossing tile is not used, and
\begin{enumerate}
\item there are (unlabeled) pipes entering from the southern boundary in the columns indexed
by $I \subset \Z$;
\item setting $(I_+,I_-) = (I \cap \Z_{>0}, \Z_{\leq 0} \setminus I)$, we have $I_{\pm} = I_{w_\lambda,\pm}$ (see \eqref{E:I+}, \eqref{E:I-});
\item the $i$-th eastmost pipe entering from the south heads off to the east in row $1-i$. (Equivalently, for every row $i \in \Z_{\leq 0}$, there is some pipe heading towards $(i,\infty)$.)
\end{enumerate}
Since crossing tiles are not used, there is no distinction between a halfplane pipedream and a $K$-halfplane pipedream.  

The weight of a halfplane pipedream $D$ is given by \eqref{eq:Kwt} (this is different from the weight used in \cite{LLS:back stable}).

\begin{ex} Let $\la=(5,3,2,2)$.
In Figure \ref{fig:halfplane bumpless pipedream} the Rothe pipedream (see \cite[Section 5.2]{LLS:back stable}) and another $\la$-halfplane bumpless pipedream are depicted.
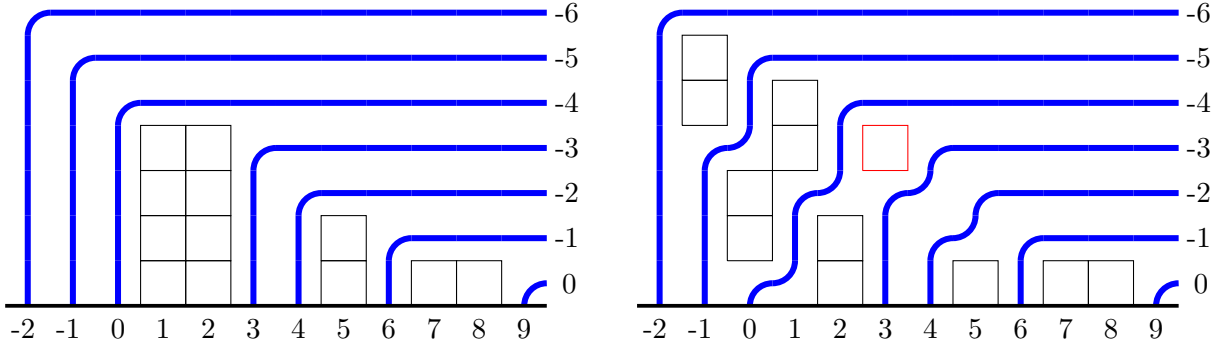
\begin{figure}
\begin{tikzpicture}[scale=0.6,line width=0.8mm]
\rightelbow{-2}{7}{\half}
\horline{-1}{7}{\half}
\horline{0}{7}{\half}
\horline{1}{7}{\half}
\horline{2}{7}{\half}
\horline{3}{7}{\half}
\horline{4}{7}{\half}
\horline{5}{7}{\half}
\horline{6}{7}{\half}
\horline{7}{7}{\half}
\horline{8}{7}{\half}
\horline{9}{7}{\half}

\vertline{-2}{6}{\half}
\rightelbow{-1}{6}{\half}
\horline{0}{6}{\half}
\horline{1}{6}{\half}
\horline{2}{6}{\half}
\horline{3}{6}{\half}
\horline{4}{6}{\half}
\horline{5}{6}{\half}
\horline{6}{6}{\half}
\horline{7}{6}{\half}
\horline{8}{6}{\half}
\horline{9}{6}{\half}

\vertline{-2}{5}{\half}
\vertline{-1}{5}{\half}
\rightelbow{0}{5}{\half}
\horline{1}{5}{\half}
\horline{2}{5}{\half}
\horline{3}{5}{\half}
\horline{4}{5}{\half}
\horline{5}{5}{\half}
\horline{6}{5}{\half}
\horline{7}{5}{\half}
\horline{8}{5}{\half}
\horline{9}{5}{\half}

\vertline{-2}{4}{\half}
\vertline{-1}{4}{\half}
\vertline{0}{4}{\half}
\bbox{1}{4}
\bbox{2}{4}
\rightelbow{3}{4}{\half}
\horline{4}{4}{\half}
\horline{5}{4}{\half}
\horline{6}{4}{\half}
\horline{7}{4}{\half}
\horline{8}{4}{\half}
\horline{9}{4}{\half}

\vertline{-2}{3}{\half}
\vertline{-1}{3}{\half}
\vertline{0}{3}{\half}
\bbox{1}{3}
\bbox{2}{3}
\vertline{3}{3}{\half}
\rightelbow{4}{3}{\half}
\horline{5}{3}{\half}
\horline{6}{3}{\half}
\horline{7}{3}{\half}
\horline{8}{3}{\half}
\horline{9}{3}{\half}

\vertline{-2}{2}{\half}
\vertline{-1}{2}{\half}
\vertline{0}{2}{\half}
\bbox{1}{2}
\bbox{2}{2}
\vertline{3}{2}{\half}
\vertline{4}{2}{\half}
\bbox{5}{2}
\rightelbow{6}{2}{\half}
\horline{7}{2}{\half}
\horline{8}{2}{\half}
\horline{9}{2}{\half}

\vertline{-2}{1}{\half}
\vertline{-1}{1}{\half}
\vertline{0}{1}{\half}
\bbox{1}{1}
\bbox{2}{1}
\vertline{3}{1}{\half}
\vertline{4}{1}{\half}
\bbox{5}{1}
\vertline{6}{1}{\half}
\bbox{7}{1}
\bbox{8}{1}
\rightelbow{9}{1}{\half}

\draw[line width=.5mm] (-2,1)--(10,1);
\draw(-1.6,.5) node {-2};
\draw(-0.6,.5) node {-1};
\draw(.5,.5) node {0};
\draw(1.5,.5) node {1};
\draw(2.5,.5) node {2};
\draw(3.5,.5) node {3};
\draw(4.5,.5) node {4};
\draw(5.5,.5) node {5};
\draw(6.5,.5) node {6};
\draw(7.5,.5) node {7};
\draw(8.5,.5) node {8};
\draw(9.5,.5) node {9};

\draw(10.5,1.5) node {0};
\draw(10.45,2.5) node {-1};
\draw(10.45,3.5) node {-2};
\draw(10.45,4.5) node {-3};
\draw(10.45,5.5) node {-4};
\draw(10.45,6.5) node {-5};
\draw(10.45,7.5) node {-6};

\begin{scope}[shift={(14,0)}]

\rightelbow{-2}{7}{\half}
\horline{-1}{7}{\half}
\horline{0}{7}{\half}
\horline{1}{7}{\half}
\horline{2}{7}{\half}
\horline{3}{7}{\half}
\horline{4}{7}{\half}
\horline{5}{7}{\half}
\horline{6}{7}{\half}
\horline{7}{7}{\half}
\horline{8}{7}{\half}
\horline{9}{7}{\half}

\vertline{-2}{6}{\half}
\bbox{-1}{6}
\rightelbow{0}{6}{\half}
\horline{1}{6}{\half}
\horline{2}{6}{\half}
\horline{3}{6}{\half}
\horline{4}{6}{\half}
\horline{5}{6}{\half}
\horline{6}{6}{\half}
\horline{7}{6}{\half}
\horline{8}{6}{\half}
\horline{9}{6}{\half}

\vertline{-2}{5}{\half}
\bbox{-1}{5}
\vertline{0}{5}{\half}
\bbox{1}{5}
\rightelbow{2}{5}{\half}
\horline{3}{5}{\half}
\horline{4}{5}{\half}
\horline{5}{5}{\half}
\horline{6}{5}{\half}
\horline{7}{5}{\half}
\horline{8}{5}{\half}
\horline{9}{5}{\half}

\vertline{-2}{4}{\half}
\rightelbow{-1}{4}{\half}
\leftelbow{0}{4}{\half}
\bbox{1}{4}
\vertline{2}{4}{\half}
{\color{red}{\bbox{3}{4}}}
\rightelbow{4}{4}{\half}
\horline{5}{4}{\half}
\horline{6}{4}{\half}
\horline{7}{4}{\half}
\horline{8}{4}{\half}
\horline{9}{4}{\half}

\vertline{-2}{3}{\half}
\vertline{-1}{3}{\half}
\bbox{0}{3}
\rightelbow{1}{3}{\half}
\leftelbow{2}{3}{\half}
\rightelbow{3}{3}{\half}
\leftelbow{4}{3}{\half}
\rightelbow{5}{3}{\half}
\horline{6}{3}{\half}
\horline{7}{3}{\half}
\horline{8}{3}{\half}
\horline{9}{3}{\half}

\vertline{-2}{2}{\half}
\vertline{-1}{2}{\half}
\bbox{0}{2}
\vertline{1}{2}{\half}
\bbox{2}{2}
\vertline{3}{2}{\half}
\rightelbow{4}{2}{\half}
\leftelbow{5}{2}{\half}
\rightelbow{6}{2}{\half}
\horline{7}{2}{\half}
\horline{8}{2}{\half}
\horline{9}{2}{\half}

\vertline{-2}{1}{\half}
\vertline{-1}{1}{\half}
\rightelbow{0}{1}{\half}
\leftelbow{1}{1}{\half}
\bbox{2}{1}
\vertline{3}{1}{\half}
\vertline{4}{1}{\half}
\bbox{5}{1}
\vertline{6}{1}{\half}
\bbox{7}{1}
\bbox{8}{1}
\rightelbow{9}{1}{\half}

\draw[line width=.5mm] (-2,1)--(10,1);
\draw(-1.6,.5) node {-2};
\draw(-0.6,.5) node {-1};
\draw(.5,.5) node {0};
\draw(1.5,.5) node {1};
\draw(2.5,.5) node {2};
\draw(3.5,.5) node {3};
\draw(4.5,.5) node {4};
\draw(5.5,.5) node {5};
\draw(6.5,.5) node {6};
\draw(7.5,.5) node {7};
\draw(8.5,.5) node {8};
\draw(9.5,.5) node {9};

\draw(10.5,1.5) node {0};
\draw(10.45,2.5) node {-1};
\draw(10.45,3.5) node {-2};
\draw(10.45,4.5) node {-3};
\draw(10.45,5.5) node {-4};
\draw(10.45,6.5) node {-5};
\draw(10.45,7.5) node {-6};
\end{scope}
\end{tikzpicture}

\caption{The Rothe and another $(5,3,2,2)$-halfplane bumpless pipedream.}
\label{fig:halfplane bumpless pipedream}
\end{figure}
\end{ex}

\begin{thm}\label{thm:Khalfplane}
Let $\la \in \Y$.  Then 
$$G_\la(x||a) = (-1)^{|\la|} \sum_D \wt(D)$$
where the summation is over all $\la$-halfplane pipedreams.
\end{thm}
\begin{proof}[Proof of Theorem \ref{thm:Khalfplane}]
Let $w_\la \in S_\Z$ be a Grassmannian permutation, and set $I = w_\la^{-1}(\Z_{\leq 0}) = I_{w_\la,+} \sqcup (\Z_{\leq 0} \setminus I_{w_\la,-})$ and $I' = \Z \setminus I$.  Let $D$ be a $K$-bumpless pipedream for $w_\la$.  Then the pipes labeled by $I$ head off to the east in the rows labeled by $\Z_{\leq 0}$, while the pipes labeled by $I'$ head off to the east in the rows labeled by $\Z_{>0}$.  Furthermore, pipes of each type do not cross pipes of the same type.  It follows immediately that the part $D_{\leq 0}$ of $D$ that lies in rows indexed by $\Z_{\leq 0}$ contains no crossing tiles and is a halfplane pipedream.

On the other hand, we claim that the bottom half $D_{>0}$ of $D$ that lies in rows indexed by $\Z_{>0}$ depends only on $w_\la$ and furthermore contains no empty tiles, and no NW-elbows.  Indeed, $D_{>0}$ is given as follows: any pipe labeled by $i \in I_+$ travels northward until row $0$, and any pipe labeled by $i \in I_-$ travels northward until row $w^{-1}(i)$, turns and travels eastward.  This description can be proven, for example, by descending induction on the label $i \in \Z$.  (Note that it follows that $D$ is actually a bumpless pipedream -- there are no pipes that cross twice.)  We also deduce that the tophalf $D_{\leq 0}$ is thus a $\la$-halfplane pipedream.

The stated formula for $G_\la(x||a) = \bG_{w_\la}(x||a)$ now follows from Theorem~\ref{thm:Kbumpless}.
\end{proof}

We restate Theorem \ref{thm:Khalfplane} using semistandard tableaux. 
Given $k\in\Z$ and a box $s=(i,j)\in\la$ define 
\begin{align}\label{E:box weight}
\wt(k,s,\la) = (x_k \ominus a_{k+\ell(\la)+j-i})
\end{align}

\begin{cor} \label{C:double G lambda by tableaux}
Let $\la \in \Y$.  Then
\begin{align}\label{E:Grothendieck by tableaux}
G_\la(x||a) = (-1)^{|\la|} \sum_T \prod_{s\in\la} \left(\wt(T(s),s,\la) \prod_{k} (1-\wt(k,s,\la))\right)
\end{align}
where $T$ runs over the semistandard tableaux of shape $\la$ with entries in $\Z_{\le0}$, 
$s$ runs over the boxes in $\la$ and $k\in\Z_{\le0}$ runs over values such that $T(s) < k$ and 
replacing the $s$-th entry of $T$ by $k$ results in a semistandard tableau.
\end{cor}
\begin{proof} Every $\la$-halfplane bumpless pipedream can be obtained from the Rothe bumpless pipedream, the unique one that has no NW elbow tiles. This corresponds to the unique semistandard tableau of shape $\la$ having $\la_i$ copies of the value $1-i$. 
Moreover each droop moves an empty tile one row north and one row west. Since the tiles move diagonally, reading along diagonals from northwest to southeast starting from the southwestmost empty tiles, the row indices of the empty tiles define the values in a corresponding diagonal in the semistandard tableau. This yields a bijection between the $\la$-halfplane bumpless pipedreams and the semistandard tableaux. The NW elbow tiles correspond to entries $k$ that can be increased while preserving semistandardness.
\end{proof}

\begin{ex} The semistandard tableaux for the pipedreams of in Figure \ref{fig:halfplane bumpless pipedream} are given by
\[
\begin{ytableau} -3&-3&-1&0&0 \\ -2 &-2&0\\-1&-1\\0&0 \end{ytableau}\qquad
\begin{ytableau} -5&-4&{\color{red}{-3}}&0&0 \\ -4 &-3&0\\-2&-1\\-1&0 \end{ytableau}
\]
The red empty tile in the second pipedream in Figure \ref{fig:halfplane bumpless pipedream} corresponds to the red $-3$ tableau entry. The fact that this $-3$ can be replaced by the larger elements $-2$ and $-1$ corresponds to the presence of the two NW elbow tiles that are directly to the southeast of the red empty tile. The contribution of the box of the red $-3$ is given by
\begin{align*}
	(x_{-3}\ominus a_3) (1 - (x_{-2}\ominus a_4)) ( 1-  (x_{-1} \ominus a_5)).
\end{align*}
\end{ex}

\begin{rem} \label{R:Grassmannian Grothendieck by pipedreams} 
\begin{enumerate}
\item In the nonequivariant setting the original formula for $G_\la$ uses set-valued tableaux \cite{Bu}. These are in bijection with bumpless pipedreams \cite[\S 7]{Wei}. Our formula ``compresses" the corresponding set-valued tableau formula by grouping a number of set-valued tableaux and bijecting this group with a single semistandard tableau and a corresponding factorized term.
\item We biject from $\la$-halfplane pipedreams to semistandard tableaux with nonpositive entries. Upon replacing entries $i$ by $1-i$ one obtains reverse semistandard tableaux with entries $\{1,2,\dotsc,\}$. In this context our formula is an equivariant upgrade of the 
symmetric Grothendieck special case of \cite[Theorem 1.3]{SY}. This result is also implicit in \cite{BSW}.
\item In the more general case of vexillary permutations, a more complicated version of \eqref{E:Grothendieck by tableaux}
was given in part 3 of the second corollary in \S 1.2 of \cite{KMY:tableau complexes}. The article \cite{SY} makes the observation that using
reverse tableaux often leads to simpler formulas.
\end{enumerate}
\end{rem}

\subsection{Expansion formulae}
Let $w \in S_n$. A \emph{$w$-rectangular $K$-bumpless pipedream}
is a $K$-bumpless pipedream in the $n \times 2n$ rectangular region
$R_n := \{(i, j) \mid i \in [1, n] \text{ and } j \in [1-n, n]\}$.
The pipes are labeled $1-n,2-n,\ldots,0,1,\ldots,n$ entering the south boundary from left to right.
The positively labeled pipes exit the east boundary, and determine $w \in S_n$ using the same prescription as for $K$-bumpless pipedreams.  The nonpositively labeled
pipes exit the north boundary, and these pipes cannot intersect any other pipe.  The weight of a $w$-rectangular $K$-bumpless pipedream is again given by \eqref{eq:Kwt}.  We also associate a partition $\la(D)$ to an $S_n$-rectangular $K$-bumpless pipedream: it is obtained by
reading the north boundary edges from right to left, to then obtain the boundary of a partition
inside a $n\times n$ box, where empty edges correspond to steps to the left, and edges with a pipe exiting
correspond to downward steps.

\begin{thm}\label{thm:Krect}
Let $w \in S_n$.  Then
$$
\bG_w(x;a) = (-1)^{\ell(w)} \sum_D (-1)^{|\la(D)|} \wt(D) G_{\la(D)}(x||a)
$$
where the summation is over $w$-rectangular pipedreams.
\end{thm}
\begin{proof}
Let $D$ be a $K$-bumpless pipedream for $w \in S_n$.  The pipes labeled by $n+1,n+2,\ldots$ travel northward and turn east in the row corresponding to their labels.  The pipes labeled by $1,2,\ldots,n$ travel northward until row $n$, perform one or more turns inside the square $[1,n] \times [1,n]$ and then travel eastward once they exit the square in one of the rows $1,2,\ldots,n$.  The pipes labeled by $0,-1,-2,\ldots,$ cannot cross other pipes.  In particular, the pipe labeled $i \in \Z_{\leq 0}$ travels northward until row $n+i-1$ before it makes its first turn.  Thus the pipes labeled $0,-1,-2,\ldots,1-n$ travel northward until row $n$, possibly perform some turns, and then exit the north boundary of the rectangle $R_n$.  The pipes labeled $-n,-n-1,\ldots$ travel vertically at least until row $0$.

To summarize: (1) the top half $D_{\leq 0}$ of $D$ is a $\la$-halfplane pipedream for some $\la \in \Y$; (2) the interesting part of the bottom half $D_{>0}$ of $D$ is contained in the rectangular region $R_n$, which in particular contains all the empty tiles and NW-elbows of $D_{>0}$.  Since $\wt(D) = \wt(D_{\leq 0})\wt(D_{>0})$, we obtain the stated formula by combining Theorem~\ref{thm:Kbumpless} with Theorem~\ref{thm:Khalfplane}.
\end{proof}

The following result follows immediately from Theorem~\ref{thm:Krect} and the definition of $G_w(x||a)$.
\begin{cor}\label{cor:Gw}
Let $w \in S_n$.  Then
$$
G_w(x||a) = (-1)^{\ell(w)} \sum_D (-1)^{|\la(D)|} \eta_a(\wt(D)) G_{\la(D)}(x||a)
$$
where the summation is over $w$-rectangular pipedreams.
\end{cor}

\begin{ex}\label{ex:rectbumplesseasy}
Let $w = s_1 \in S_2$.  Theorem~\ref{thm:Krect} gives
$$
\bG_{s_1}(x||a) = (x_1\ominus a_1) + (1-(x_1 \ominus a_1))G_1(x||a),
$$
and noting that $G_1(x||a) = G_1(x/a)$ this agrees with Example~\ref{ex:doublecoproduct}.
Corollary~\ref{cor:Gw} gives
$$
G_{s_1}(x||a) = G_1(x||a).
$$
\end{ex}

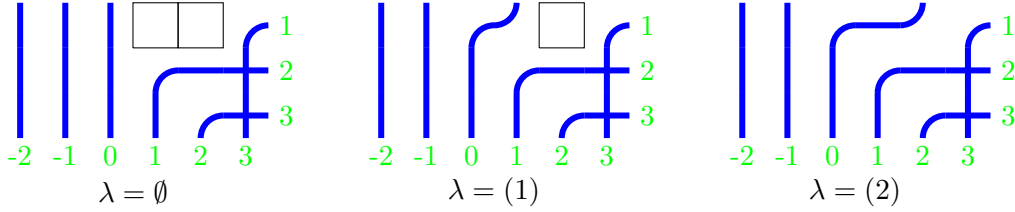
\begin{figure}
\centering
\begin{tikzpicture}[scale=0.6,line width=0.8mm]
\vertline{-2}{-1}{blue}
\vertline{-1}{-1}{blue}
\vertline{0}{-1}{blue}
\bbox{1}{-1}{blue}
\bbox{2}{-1}{blue}
\rightelbow{3}{-1}{blue}

\vertline{-2}{-2}{blue}
\vertline{-1}{-2}{blue}
\vertline{0}{-2}{blue}
\rightelbow{1}{-2}{blue}
\horline{2}{-2}{blue}
\cross{3}{-2}{blue}{blue}

\vertline{-2}{-3}{blue}
\vertline{-1}{-3}{blue}
\vertline{0}{-3}{blue}
\vertline{1}{-3}{blue}
\rightelbow{2}{-3}{blue}
\cross{3}{-3}{blue}{blue}

%

\draw[green] (-1.5,-3.4) node {-2};
\draw[green] (-0.5,-3.4) node {-1};
\draw[green] (0.5,-3.4) node {0};
\draw[green] (1.5,-3.4) node  {$1$};
\draw[green] (2.5,-3.4) node  {$2$};
\draw[green] (3.5,-3.4) node  {3};

\draw[green] (4.4,-.5) node {1};
\draw[green] (4.4,-1.5) node {2};
\draw[green] (4.4,-2.5) node {3};

\draw[black] (1,-4.2) node {$\la=\emptyset$};

\begin{scope}[shift={(8,0)}]
\vertline{-2}{-1}{blue}
\vertline{-1}{-1}{blue}
\rightelbow{0}{-1}{blue}
\leftelbow{1}{-1}{blue}
\bbox{2}{-1}{blue}
\rightelbow{3}{-1}{blue}

\vertline{-2}{-2}{blue}
\vertline{-1}{-2}{blue}
\vertline{0}{-2}{blue}
\rightelbow{1}{-2}{blue}
\horline{2}{-2}{blue}
\cross{3}{-2}{blue}{blue}

\vertline{-2}{-3}{blue}
\vertline{-1}{-3}{blue}
\vertline{0}{-3}{blue}
\vertline{1}{-3}{blue}
\rightelbow{2}{-3}{blue}
\cross{3}{-3}{blue}{blue}

%

\draw[green] (-1.5,-3.4) node {-2};
\draw[green] (-0.5,-3.4) node {-1};
\draw[green] (0.5,-3.4) node {0};
\draw[green] (1.5,-3.4) node  {$1$};
\draw[green] (2.5,-3.4) node  {$2$};
\draw[green] (3.5,-3.4) node  {3};

\draw[green] (4.4,-.5) node {1};
\draw[green] (4.4,-1.5) node {2};
\draw[green] (4.4,-2.5) node {3};

\draw[black] (1,-4.2) node {$\la=(1)$};
\end{scope}

\begin{scope}[shift={(16,0)}]
\vertline{-2}{-1}{blue}
\vertline{-1}{-1}{blue}
\rightelbow{0}{-1}{blue}
\horline{1}{-1}{blue}
\leftelbow{2}{-1}{blue}
\rightelbow{3}{-1}{blue}

\vertline{-2}{-2}{blue}
\vertline{-1}{-2}{blue}
\vertline{0}{-2}{blue}
\rightelbow{1}{-2}{blue}
\horline{2}{-2}{blue}
\cross{3}{-2}{blue}{blue}

\vertline{-2}{-3}{blue}
\vertline{-1}{-3}{blue}
\vertline{0}{-3}{blue}
\vertline{1}{-3}{blue}
\rightelbow{2}{-3}{blue}
\cross{3}{-3}{blue}{blue}

%

\draw[green] (-1.5,-3.4) node {-2};
\draw[green] (-0.5,-3.4) node {-1};
\draw[green] (0.5,-3.4) node {0};
\draw[green] (1.5,-3.4) node  {$1$};
\draw[green] (2.5,-3.4) node  {$2$};
\draw[green] (3.5,-3.4) node  {3};

\draw[green] (4.4,-.5) node {1};
\draw[green] (4.4,-1.5) node {2};
\draw[green] (4.4,-2.5) node {3};
\draw[black] (1,-4.2) node {$\la=(2)$};
\end{scope}
\end{tikzpicture}
\caption{$w$-rectangular $K$-bumpless pipedreams for $w = s_2 s_1$.}
\label{fig:rectbumplessexample}
\end{figure}

\begin{ex}\label{ex:rectbumpless}
Let $w = s_2s_1 \in S_3$.  In one line notation, $w(1,2,3) = (3,1,2)$.  The $w$-rectangular $K$-bumpless pipedreams are shown in Figure~\ref{fig:rectbumplessexample}.  By Theorem~\ref{thm:Krect}, we have
$$
\bG_w(x;a) = (x_1\ominus a_1)(x_1 \ominus a_2) +(x_1\ominus a_2)(1-(x_1\ominus a_1))G_1(x||a) + (1-(x_1\ominus a_2))G_2(x||a),
$$
and by Corollary~\ref{cor:Gw}, we have
$$
G_w(x||a) = (a_1\ominus a_2)G_1(x||a) + (1-(a_1\ominus a_2))G_2(x||a).
$$
\end{ex}

\subsection{Weigandt's formula for double Grothendieck polynomials}
Let $w \in S_n$.  A \emph{$w$-square $K$-bumpless pipedream} is a bumpless pipedream in the $ n \times n$ square region $[n] \times [n]$. The pipes are labeled $1,\ldots,n$ entering the south boundary from left to right, and all pipes exit the east boundary.  The permutation $w \in S_n$ is determined as for $K$-bumpless pipedreams.  The weight $\wt(D)$ of a $w$-square $K$-bumpless pipedream is defined by \eqref{eq:Kwt} as before.  Weigandt's formula \cite[Theorem 1.1]{Wei} for double Grothendieck polynomials is the following.

\begin{thm}\label{thm:Wei}
Let $w \in S_n$.  Then 
$$
\G_w(x;a) = (-1)^{\ell(w)} \sum_D \wt(D)
$$
where the summation is over all $w$-square $K$-bumpless pipedreams $D$.
\end{thm}

\subsection{Hecke bumpless pipedreams}
Let $D$ be a $w$-square $K$-bumpless pipedream.  Following Weigandt \cite{Wei}, we call $D$ a \emph{Hecke bumpless pipedream} if all the empty tiles are in the northeast corner, where they form a partition shape
$\la = \la(D)$, called the shape of $D$.  The following result was proven by Weigandt \cite{Wei}.

\begin{thm}\label{thm:Heckebumpless}
The coefficient $k_w^\la$ of $G_\la$ in $G_w$ (see \eqref{E:stable K stanley coefficients}) is equal to $(-1)^{\ell(w) - |\la|}$ times the number of $w$-Hecke bumpless pipedream with shape $\la$.
\end{thm}
\begin{proof}
Substituting $a = 0$ into Corollary~\ref{cor:Gw}, only rectangular pipedreams with no empty tiles contribute.  Erasing the nonpositively labeled pipes from such a rectangular pipedream gives a $w$-Hecke bumpless pipedream.
\end{proof}

Theorem~\ref{thm:Heckebumpless} is a $K$-theoretic analogue of \cite[Theorem 5.14]{LLS:back stable}, a direct bijective proof of which was given by Fan, Guo and Sun \cite{FGS}.  Weigandt \cite{Wei} gave a more general bijection between Hecke bumpless pipedreams and decreasing tableaux \cite{BKSTY}.

\section{$K$-homology and Hopf structure}\label{sec:homology}

\subsection{Hopf structure on GKM ring}
Let $\Psi_\Gr \subset \Psi$ denote the subspace of functions $\psi$ satisfying $\psi(v) = \psi(w)$ if $v S_{\neq 0} = w S_{\neq 0}$, and similarly define $\Psit_\Gr \subset \Psit$.  Then $\Psi_\Gr$ (resp. $\Psit_\Gr$) has basis (allowing infinite sums) $\{\psi^{w_\la} \mid \la \in \P\}$ (resp. $\{\psit^{w_\la} \mid \la \in \P\}$).  We have $K_T(\Grb^0) \simeq \Psi_\Gr$.

The map $\res: \bR(x;a) \to \SSS \otimes_{R(a)} \Psit$ of Theorem~\ref{T:infinite basis} restricts to a map
$\res: \Gamma(x||a) \to \Psit$ with image given by $\bigoplus_\la R(a) \psit^{w_\la}$.  Following \cite{LLS:back stable}, we now describe a Hopf structure on $\Psi_\Gr$ that is compatible with the bialgebra structure on $\Gamma(x||a)$.  

For any $w\in S_\Z$, let
\begin{align}
	\label{E:I+}
	I_{w,+} &:= \Z_{>0} \cap w(\Z_{\le0})\\
	\label{E:I-}
	I_{w,-} &:= \Z_{\le0} \cap w(\Z_{>0})
\end{align}
The map $w\mapsto (I_{w,+},I_{w,-})$ is a bijection
from $S_\Z^0$  to pairs of finite sets $(I_+,I_-)$ such that $I_+\subset \Z_{>0}$, $I_-\subset \Z_{\le0}$, and $|I_+|=|I_-|$.  There is a partial multiplication map $S_\Z^0 \times S_\Z^0 \to S_\Z^0$.  The product of $x \in S_\Z^0$ and $y\in S_\Z^0$ is equal to $z \in S_\Z^0$ if (1) $I_{x,+} \cap I_{y,+} = \emptyset = I_{x,-} \cap I_{y,-}$ and (2) $I_{x,\pm} \cup I_{y,\pm} = I_{z,\pm}$.  

The following result is proved in the same manner as \cite[Proposition 7.11]{LLS:back stable}.
\begin{prop}\label{prop:locHopf}
	There is a coproduct $\Delta: \Psit_\Gr \to \Psit_\Gr \hat \otimes_{R(a)} \Psit_\Gr$ such that the map $\res: \Gamma(x||a) \to \Psit_\Gr$ is a $R(a)$-bialgebra morphism.  The coproduct satisfies
	\begin{equation}\label{eq:GrHopf}
	\Delta(\psi)|_{x \otimes y} = \psi|_{xy}
	\end{equation}
	whenever $x,y, \in S_\Z^0$ and $xy \in S_\Z^0$ is defined. 	
\end{prop}

\subsection{$K$-homology basis}
\def\tLa{{\tilde \La}}
\def\tdelta{{\tilde \delta}}
\def\tGdiff{{\tilde \pi}}

Let $\La(y)$ denote the $\Z$-algebra of symmetric functions in $y=y_{\le0} = (y_0,y_{-1},y_{-2},\dotsc$ and $\tLa(y||a)=\prod_{\la\in\P} R(a) s_\la(y)$ the completion of $R(a) \otimes \La(y)$ whose elements are formal (possibly infinite) $R(a)$-linear combinations $\sum_{\lambda \in \P} a_\lambda s_\lambda(y)$ of Schur functions, with $a_\lambda \in R(a)$.  The ring $\tLa(y||a)$ is a $R(a)$-Hopf algebra with coproduct
$\Delta(p_k(y)) = 1 \otimes p_k(y) + p_k(y) \otimes 1$.

Define the Cauchy kernel
\begin{align*}
\Omega[(\xm-\am)y] = \prod_{i,j\le0} \dfrac{1-a_iy_j}{1-x_iy_j} = \exp\left(\sum_{k \geq 0} \frac{1}{k}p_k(x||a) p_k(y)\right).
\end{align*}
This induces the structure of dual $R(a)$-Hopf algebras on $\tLa(y||a)$ and (a completion of) $\Gamma(x||a)$.
Write $\pair{\cdot}{\cdot}$ for the corresponding pairing $\Gamma(x||a) \otimes_{R(a)} \tLa(y||a) \to R(a)$. Then by definition
\begin{align}\label{E:equivariant pairing}
	\pair{s_\la(x/a)}{s_\mu(y)} = \delta_{\la,\mu}.
\end{align}
Let $g_\la(y)$ be the \emph{dual stable Grothendieck polynomials} of \cite{LP}.  They are defined by
\begin{align}
	\pair{G_\la(x/a)}{g_\mu(y)} = \delta_{\la,\mu}.
\end{align}
Define the \emph{$K$-Molev functions} $g_\la(y||a)\in\tLa(y||a)$ by duality with the Grassmannian double $K$-Stanley functions $G_\la(x||a)$:
\begin{align}\label{E:dual Schur}
	\pair{G_\la(x||a)}{g_\mu(y||a)} &= \delta_{\la\mu}.
\end{align}
The ring $\tLa(y||a)$ consists of formal $R(a)$-linear combinations of the $g_\la(y||a)$.  At $a = 0$, the polynomials $g_\la(y||a)$ reduce to the dual stable Grothendieck polynomials $g_\la(y)$.  The functions $g_\la(y||a)$ are $K$-theoretic analogues of Molev's dual Schur functions.  

\begin{rem}
Since the Grassmannian double $K$-Stanley functions $G_\la(x||a)$ represent the structure sheaves of Schubert varieties in the equivariant $K$-theory $K_T(\Grb^0)$ of the thick infinite Grassmannian (Corollary \ref{C:Grass Schubert basis}), the functions $g_\la(y||a)$ represent the dual basis of ideal sheaves of boundaries of Schubert varieties in the $K$-group $K^T(\Gr^{(0)})$ of finitely supported equivariant coherent sheaves on the thin infinite Grassmannian. Note that the Schubert varieties of $\Grb^0$ and of $\Gr^{(0)}$ are ``opposite".  See \cite{Kum, BK} for details on the duality between thin and thick $K$-groups, and \cite{LSS} for a discussion of the Hopf structure. 
\end{rem}

Recall the element $w_{\la/\mu} \in S_\Z$ from \eqref{E:w skew}.
Proposition \ref{P:double to super} implies the following.

\begin{prop}\label{prop:dualsuper} For $\mu \in \P$, we have
\begin{align}
   g_\mu(y) &= \sum_{\substack{\nu \subset \mu \subset \la \\ d(\nu)=d(\la)\\ \mu/\nu\in \RS}} (-1)^{|\mu|-|\nu|} \G_{w_{\la/\nu}^{-1}}(\ominus(a)) g_\la(y||a) \\
   g_\mu(y||a) &= \sum_{\substack{\nu \subset \mu \subset \la \\ d(\nu)=d(\la)\\ \mu/\nu\in \RS}} (-1)^{|\mu|-|\nu|} \G_{w_{\la/\nu}}(a) g_\la(y).
\end{align}
\end{prop}

\subsection{$K$-Homology divided difference operators}\label{SS:homology divided differences}
Recall that by convention we have $x:= x_{\leq 0} = (\dotsc,x_{-1},x_0)$ and similarly, $y:= y_{\le 0}$ and $a:= a_{\le 0}$.

For $i\in\Z$, define the operators (see \eqref{E:parh})
\begin{align}
	\ts_i^a &= \Omega[(x-a)y] s_i^a \Omega[(a-x)y] \\
  \delta^{a,\anti}_i &= \Omega[(x-a)y] \Gdiff_i^{a,\anti} \Omega[(a-x)y].
\end{align}
It is clear that these operators, being conjugate to the operators $s_i^a$ and $\Gdiff_i^{a,\anti}$ 
respectively, satisfy the type A braid relations. In operator expressions, a symmetric function or polynomial $f$ denotes left multiplication by $f$. We have
\begin{align}
	\label{E:sdelta}
	\ts_i^a \delta^{a,\anti}_i &= \delta^{a,\anti}_i \\ 
	\label{E:deltaf}
	\delta^{a,\anti}_i f &= \delta^{a,\anti}_i(f) + \ts^a_i(f) \Gdiff_i^{a,\anti}.
\end{align}
Let $\alpha_i = a_i - a_{i+1}$ for $i\in\Z$. 
Since $\Omega[(a-x)y]$ is $s_i^a$ invariant for $i\ne0$ we have
\begin{align}\label{E:nonzero s}
	\ts_i^a &= s_i ^a&\qquad&\text{for $i\ne0$}\\
\label{E:nonzero delta}
	\delta^{a,\anti}_i &= \Gdiff_i^{a,\anti}  &&\text{for $i\ne 0$.} \\
\label{E:zero s}
	\ts_0^a &= \Omega[-\alpha_0 y] s_0^a \\
\label{E:zero delta}
\delta^{a,\anti}_0&= \alpha_0^{-1}(1-\ts_0)(a_0-1).
\end{align}


The diagonal index of a box in row $i$ and column $j$ is by definition $j-i$. For $\la\in \P$ and $d\in \Z$, let $\la+d$ denote the partition obtained by adding a corner to $\la$ in the $d$-th diagonal if such a corner exists, and $\la+d := \la$ if such a partition does not exist. Define $\la-d$ similarly for removal of the corner in diagonal $d$.

By Proposition~\ref{P:equiv ddiff bG}, we have
\begin{align}\label{E:left ddiff double Schur}
	\Gdiff_i^{a,\anti} G_\la(x||a) = G_{\la-i}(x||a) \qquad\text{for all $\la\in\P$ and $i\in\Z$.}
\end{align}
Let $\parh_i^a = \Gdiff_i^a - 1 = (1-a_i)A^a_i$ and $\tau_i = \Omega[(x-a)y] \parh_i^a \Omega[(a-x)y]$.

\begin{lem}\label{L:the ops} For $i \in \Z$, we have
\begin{align}
\label{E:ops}
\parh_i^a &= s_i^a - \Gdiff_i^{a,\anti} \\
\label{E:conj ops}
\tau_i &= \ts_i^a - \delta^{a,\anti}_i.
\end{align}
\end{lem}
\begin{proof} Equation \ref{E:conj ops} follows from \eqref{E:ops} by conjugation by $\Omega[(x-a)y]$. We have
\begin{align*}
	s_i^a-\Gdiff_i^{a,\anti} &= s_i^a - A^a_i(a_i-1) \\
	&=s_i^a- (1 + (a_{i+1}-1)A^a_i) \\
	&= (-\alpha_i -(a_{i+1}-1) A^a_i \\
	&= (1-a_i)A^a_i \\
	\parh_i^a &= \Gdiff_i^a - 1 =- 1+ A^a_i(1-a_{i+1}) \\
	&= (1-a_i) A^a_i .  \qedhere
\end{align*}
\end{proof}

\begin{prop} \label{P:left ddiff homology} For all $\mu\in\P$ and $i\in\Z$, we have
\begin{align} \label{E:left ddiff homology}
	\tau_i(g_\mu(y||a)) &=  g_{\mu+i}(y||a).
\end{align}
\end{prop}
\begin{proof} Using \eqref{E:deltaf}, we have
\begin{align*}
	\Omega[(x-a)y] &= \Omega[(x-a)y] \,\Gdiff_i^{a,\anti}(1)= \delta^{a,\anti}_i(\Omega[(x-a)y]) \\
	&= \delta^{a,\anti}_i \sum_\la  g_\la(y||a) G_\la(x||a) \\
	&= \sum_\la \left (\delta^{a,\anti}_i(g_\la(y||a))G_\la(x||a) + \tilde{s}_i^a(g_\la(y||a)) \Gdiff_i^{a,\anti}(G_\la(x||a))     \right) \\
	&= \sum_\la \left(\delta^{a,\anti}_i(g_\la(y||a))G_\la(x||a) +\tilde{s}_i^a(g_\la(y||a)) G_{\la-i}(x||a)\right)
\end{align*}
Taking the coefficient of $G_\mu(x||a)$ we obtain
\begin{align*}
	g_\mu(y||a) &= \delta^{a,\anti}_i(g_\mu(y||a)) + \ts_i^a(g_{\mu+i}(y||a)).
\end{align*}
Acting by $\ts_i^a$, using \eqref{E:sdelta}, and rearranging we have
\begin{align*}
	g_{\mu+i}(y||a) &= (\ts_i^a-\delta^{a,\anti}_i)(g_\mu(y||a)) = \tau_i(g_\mu(y||a)). \qedhere
\end{align*}
\end{proof}

The $\tau_i$ satisfy the type $A$ braid relations and $\tau_i^2 = -\tau_i$.  Thus $\tau_w=\tau_{i_1}\dotsm\tau_{i_\ell}$ makes sense for any reduced decomposition $w=s_{i_1}\dotsm s_{i_\ell}\in S_\Z$.

\begin{thm} \label{T:create Khom} For any $\la\in\P$, we have $g_{\la}(y||a) = \tau_{w_\la}(1)$.
\end{thm}

\begin{ex}
We have 
\begin{align*}
g_1(y||a) &= \tau_0(1) = \Omega[(x-a)y](1-a_0)A^a_0 \Omega[(a-x)y]  \\
&=(1-a_0)  \alpha_0^{-1} \Omega[-a_0y] (1-s_0^a) \Omega[a_0y] \\
&= (1-a_0) \alpha_0^{-1} (1 - \Omega[(a_1-a_0)y]) \\
&= (1-a_0) \sum_{i,j\ge0} (-a_0)^i a_1^j s_{(j+1,1^i)}(y) .
\end{align*}
In particular, setting all $a_i=0$ we obtain $g_1(y)=s_1(y)$.
Let $Z= \sum_{i,j\ge0} (-a_0)^i a_1^j s_{(j+1,1^i)}(y)$. We have
\begin{align*}
	g_{11}(y||a) &= \tau_{-1} g_1(y||a) = \parh_i^a g_1(y||a) = (1-a_{-1}) A^a_{-1} (1-a_0) Z  \\
	&= (1-a_{-1}) (1+(1-a_{-1})A^a_{-1}) Z\\
	&= (1-a_{-1})(1+(1-a_{-1})A^a_{-1})\left( s_1(y) -a_0 s_{11}(y) + a_1 s_2(y)+\dotsm\right)\\
	&= (1-a_{-1})(s_1(y) - a_0 s_{11}(y) + \dotsm )+ (1-a_{-1})^2(s_{11}(y)\dotsm).
\end{align*}
We obtain $g_{11}(y) = s_{11}+s_1$. These computations of $g_1$ and $g_{11}$ agree with \cite{LP}.
\end{ex}
%

\subsection{Connection to Knutson-Lederer}
Knutson and Lederer \cite{KL} define a deformation of the ring of symmetric functions, denoted $R^{K^S}$.  The ring $R^{K^S}$ has a basis $[X^\la] = [{\mathcal O}_{X^\la}]$ (representing structure sheaves of opposite Schubert varieties in the Grassmannian) as $\la$ ranges over all partitions.  The product structure of $R^{K^S}$ is given by the direct sum operation on Grassmannians.  We refer the reader to \cite{KL} for the details and to \cite[Section 8.3]{LLS:back stable} for a synopsis of the very similar situation in homology.  

Let $h_\la(y||a) := \sum_{\mu \subset \la} g_\la(y||a)$.  Whereas $g_\la(y||a)$ represents the ideal sheaf of a boundary of a Schubert variety in $K^T(\Gr)$, the symmetric function $h_\la(y||a)$ represents the structure sheaf of the same Schubert variety.  Further, let $h_\la(y||\delta)$ be obtained from $h_\la(y||a)$ by the specialization 
$$
a_i \longmapsto \begin{cases} \delta &\mbox{if $i > 0$,}\\
0 & \mbox{if $i \leq 0$.}
\end{cases}
$$
See \cite[Section 8.3]{LLS:back stable} for a more precise description.  The following result is proved in the same manner as \cite[Theorem 8.12]{LLS:back stable}.  

\begin{thm}\label{T:KL}
There is, up to a completion, an isomorphism of $\Z[\delta]$-algebras
$$
R^{K^S} \to \Lambda(y)[\delta] \qquad [X^\la] \mapsto h_\la(y||\delta)
$$
where $\delta$ corresponds to $1-\exp(-t)$ in \cite{KL}.
\end{thm}

\section{The ring of back stable Grothendieck polynomials}

\subsection{The subring of $K$-Stanley functions}
Define 
$$\Gamma=\Gamma(\xm) := \bigoplus_{\la\in\Y} \Z G_\la\subset \hLa$$ 
to be the span of Grassmannian $K$-Stanley functions. The structure of $\Gamma$ was studied by Buch \cite{Bu}.

\begin{prop} \label{P:stable K stanley expansion} \cite{Bu} 
$\Gamma$ is a commutative and cocommutative bialgebra  containing $G_w$ for all $w \in S_\Z$.	
In particular, the expansion
\begin{align}\label{E:stable K stanley coefficients}
	G_w &= \sum_\la k_\la^w G_\la
\end{align}
exists and is finite, and furthermore we have $(-1)^{\ell(w)-|\la|} k_\la^w \in \Z_{\ge0}$.
\end{prop}

For $\la,\mu,\nu \in \P$, define $c_{\mu\nu}^\la\in \Z$ by
\begin{align*}
	G_\mu G_\nu = \sum_{\la} c^\la_{\mu\nu} G_\la.
\end{align*}
By \cite{Bu}, we have the positivity property  $c^\la_{\mu\nu} \in (-1)^{|\mu|+|\nu|-|\la|} \Z_{\ge0}$.

\begin{rem} Explicit combinatorial formulae for the product and coproduct structure constants of $\Gamma$ with respect to the basis $G_\la$ were given by Buch \cite{Bu}.
\end{rem}
\begin{rem}\label{rem:BKSTY}
The following explicit tableau formula is given in \cite[Theorem 1']{BKSTY}: $(-1)^{\ell(w)-|\la|}k_\la^w$ is equal to the number of \emph{decreasing} tableaux (rows strictly decrease from left to right and columns strictly decrease from top to bottom) $T$ of shape $\la$ whose column-reading word is a Hecke word for $w$, that is, if the column-reading word of $T$ is $i_1i_2\dotsm i_\ell$ then $s_{i_1}*s_{i_2}*\dotsm*s_{i_\ell} = w$.
\end{rem}

\subsection{Flagged Grothendieck polynomials and $K$-Stanley functions via symmetrization}
\label{SS:flagged Groth}
A polynomial truncation of the $K$-Stanley function $G_w$ can be obtained from the Grothendieck polynomial $\G_w$ by symmetrization operators. The intermediate polynomials are what we shall call \emph{flagged Grothendieck polynomials}.

Let $\pit_i := A_i \circ x_i (1-x_{i+1})$ for $i\in\Z$. These operators satisfy the braid relations and are idempotent: $\pit_i^2=\pit_i$.  Thus they generate a 0-Hecke algebra.

Say that a sequence $f=(f_1,f_2,\dotsc,f_{n-1},f_n)$ of positive integers is \emph{admissible} if
$1\le f_1 \le f_2 \le\dotsm \le f_{n-1} \le f_n \le n$ and $f_i\ge i$ for all $i$.
Then either
$f=\fmin=(1,2,\dotsc,n-1,n)$ or there is a minimum $i$ such that $f_i>i$.
Let $f^-$ be $f$ with $f_i$ replaced by $f_i-1$. Then define $\sigma_f\in S_n$ by
\begin{align}
	\sigma_f = \begin{cases} \id & \text{if $f=\fmin$} \\
		s_{f_i-1} \sigma_{f^-} & \text{if $i$ is minimum such that $f_i>i$.}
	\end{cases}
\end{align}
Also, define the sequence $f'=(f'_1,f'_2,\ldots,f'_n)$ by $f'_i = \min\{j\mid f_j\ge i \}$.

\begin{ex} Let $n=7$ and $f=(3,4,4,5,7,7,7)$. We have $f'=(1,1,1,2,4,5,5)$ and
$$\sigma_f  =(s_2s_1)(s_3s_2) (s_3) (s_4)(s_6s_5) (s_6).$$
We illustrate the construction.  In the following diagram the $j$-th column has size $f_j$ and a box in row $i+1$ below the diagonal has corresponding simple reflection $s_i$.  We have $f'_1=1$ and for $i\ge2$, in the diagram $f'_i$ is the leftmost column in the $i$-th row containing a green or gray square. 
	\begin{align*}
		\begin{ytableau}
			*(gray)&*(gray) &*(gray) &*(gray) &*(gray) &*(gray) &*(gray)  \\
			*(green) 1& *(gray) & *(gray) & *(gray) & *(gray) & *(gray) & *(gray) \\
			*(green) 2&*(green) 2&*(gray) &*(gray) &*(gray) &*(gray) &*(gray) \\
			&*(green) 3&*(green) 3&*(gray) &*(gray) &*(gray) &*(gray) \\
			&&&*(green) 4&*(gray) &*(gray) &*(gray) \\
			&&&&*(green) 5&*(gray) &*(gray) \\
			&&&&*(green) 6&*(green) 6&*(gray) \\
		\end{ytableau}
	\end{align*}
\end{ex}

Let $f=(f_1,f_2,\dotsc,f_{n-1},f_n)$ be admissible.
Define the \emph{flagged Grothendieck polynomial} $\G_{w,f}$ by
\begin{align*}
	\G_{w,f} := \pit_{\sigma_f}(\G_w).
\end{align*}

\begin{prop}\label{P:flagged G} For $w \in S_+$, we have
	\begin{align}\label{E:flagged G}
		\G_{w,f} &= \sum_{\substack{s_{a_1}*s_{a_2}*\dotsm*s_{a_p}=w\\1\le  i_1\le i_2\le \dotsm \le i_p \\ a_k\le a_{k+1} \Rightarrow i_k<i_{k+1} \\ i_k\le f_{a_k}}} (-1)^{p-\ell(w)} x_{i_1}x_{i_2}\dotsm x_{i_p}, 
	\end{align}
	where $p$ is arbitrary. 
\end{prop}

Note that compared to \eqref{E:G back stable}, only the bound $i_k \le a_k$ has been changed to 
$i_k \le f_{a_k}$. 

\begin{cor} \label{C:symmetrize G} For $w\in S_n$, we have
	\begin{align}\label{E:symmetrize G}
		\pit_{w_0}(\G_w) = G_w(x_1,\dotsc,x_n).
	\end{align}
\end{cor}
\begin{proof} Apply Proposition \ref{P:flagged G} with $f=(n,n,\dotsc,n)$, which satisfies $\sigma_f=w_0$, and compare with Proposition \ref{P:K Stanley}.
\end{proof}

To prove Proposition \ref{P:flagged G} we use 0-Hecke algebra generating functions as in \cite{FK}.  Consider the algebra $\bA$ over $P=\Z[x_1,\dotsc,x_n]$ generated by elements $u_i$ for $i$ in the Dynkin node set $ I = \{1,2,\dotsc,n-1\}$ 
of type $A_{n-1}$, which satisfy the type $A_{n-1}$ braid relations and $u_i^2= - u_i$ for all $i\in I$. For $w \in S_n$, the element $u_w = u_{i_1} u_{i_2}\dotsm u_{i_\ell}$ is well-defined, where $w=s_{i_1}s_{i_2}\dotsm s_{i_\ell}$ is a reduced expression.
Then 
\begin{align}
	\bA = \bigoplus_{w\in S_n} P  u_w.
\end{align}
We have $(1+au_i)(1+bu_i)=1+(a\oplus b)u_i$. In particular, $1+x_iu_j$ is invertible in $\bA_F:= F \otimes_P\bA$, where $F = \Frac(P)$ denotes the fraction field of $P$.  By \cite{FK}, we have the following identity:
\begin{align}\label{E:FK}
	\prod_{j=1}^{n-1} \prod_{i=n-1}^j (1 + x_j u_i) = \sum_{w\in S_n} \G_w u_w
\end{align}
where in the inner product the index $i$ goes from $n-1$ down to $j$ going from left to right. Taking the coefficient of $u_w$ one obtains the monomial expansion \eqref{E:Groth FK}.

\begin{prop} \label{P:flagged G GF}
	For admissible $f$, we have
	\begin{align*}
		\G_f(x,u)&:= \prod_{i=1}^n \prod_{j=n-1}^{f'_i} (1+x_i u_j)
		= \sum_{w\in S_n} \G_{w,f} u_w.
	\end{align*}
\end{prop}
\begin{proof} Follows by Lemmas \ref{L:diff flagged G} and \ref{L:pit GF}.
\end{proof}

\begin{proof}[Proof of Proposition~\ref{P:flagged G}] This follows from  by taking the coefficient of $u_w$ in $\G_f(x,u)$, noting that $j \ge f'_i$ if and only if $f_j \ge i$.
\end{proof}

\begin{lem}\label{L:diff flagged G}
	Suppose the value $0<k<n$ occurs exactly once in $f$. Then
	\begin{align}
		\pit_k(\G_{w,f}) = \G_{w,g},
	\end{align}
	where $g$ is obtained from $f $ by replacing $k$ by $k+1$ in $f$.
\end{lem}
\begin{proof} This is easily proved by induction.
\end{proof}

\begin{lem} \label{L:pit GF} With the same assumptions as in Lemma \ref{L:diff flagged G}, we have
	\begin{align*}
		\pit_k(\G_f(u,x)) = \G_{g}(u,x).
	\end{align*}
	Here, the operator $\pit_k$ is acting on the coefficients of elements in $\bA$.
\end{lem}
\begin{proof} Let $j$ be such that $f_j=k$. We are assuming that $f_{j+1}>k$.
	Thus $f'_k=j$ and $f'_{k+1}=j+1$.
	The operator $\pit_k$ commutes with all operators $1 + x_i u_j$ except when
	$i\in \{k,k+1\}$. It therefore suffices to show that
	\begin{align*}
		&\pit_k \left((1+x_k u_{n-1})  \dotsm (1+x_k u_j)
		(1+x_{k+1}u_{n-1})\dotsm (1+x_{k+1}u_{j+1}) \right) \\
		&= (1+x_k u_{n-1})  \dotsm (1+x_k u_j)
		(1+x_{k+1}u_{n-1})\dotsm (1+x_{k+1}u_j).
	\end{align*}
	Without loss of generality we may assume $k=1$ and $j=1$.
	Let $$h = \prod_{i=3}^{n-1} \prod_{j=n-1}^i (1+x_i u_j).$$ This element is invertible in $\bA_F$ and commutes both with
	$\pit_1$ and $1+x_2u_1$. It therefore suffices to show 
	\begin{align*}
		\pit_1 \left(\prod_{i=1}^{n-1} \prod_{j=n-1}^i (1+x_iu_j)\right) = \left(\prod_{i=1}^{n-1} \prod_{j=n-1}^i (1+x_iu_j) \right)(1+x_2u_1)
	\end{align*}
	as this is the required identity after multiplication by $h$. We have 
	$\pit_1 = A_1 x_1 (1-x_2) = (1-x_2) + x_2 A_1 (1-x_2) = (1-x_2)+x_2 \pih_1$. Using \eqref{E:FK}, we compute
	\begin{align*}
		\pit_1 \sum_{w\in S_n} \G_w u_w 
		&= (1-x_2)\sum_{w\in S_n} \G_w u_w + x_2 \sum_{\substack{w\in S_n \\ ws_1<w}} \G_{ws_1} u_w 
		+ x_2 \sum_{\substack{w\in S_n \\ ws_1>w}} \G_w u_w \\
		&= \sum_{w\in S_n} \G_w u_w + \sum_{\substack{w\in S_n \\ ws_1<w}} (-x_2 \G_w+ x_2 \G_{ws_1}) u_w,
	\end{align*}
	and
	\begin{align*}
		\left(\sum_{w\in S_n} \G_w u_w\right) (1+x_2u_1) &= \sum_{w\in S_n} \G_wu_w + x_2 \sum_{w\in S_n} \G_w u_w u_1 \\
		&= \sum_{w\in S_n} \G_wu_w +x_2 \sum_{\substack{w\in S_n\\ ws_1>w}} \G_w u_{ws_1} -x_2 \sum_{\substack{w\in S_n\\ ws_1<w}} \G_w u_w  \\
		&= \sum_{w\in S_n} \G_wu_w +x_2 \sum_{\substack{w\in S_n\\ ws_1<w}} \G_{ws_1}  u_w - x_2 \sum_{\substack{w\in S_n\\ ws_1<w}} \G_w u_w, 
	\end{align*}
	as required.
\end{proof}

\subsection{The subring of back stable Grothendieck polynomials}
The back stable Grothendieck polynomials $\bG_w$ are defined as elements of the ring $\bR^+$ of \S \ref{SS:sym func}. They are linearly independent because their lowest degree components are back stable Schubert polynomials, which are linearly independent. However, they do not span $\bR^+$, or even the ``finite" subring $\Gamma \otimes R^+ \subset \bR^+$.  The following example shows that $x_1$ and $s_1$ are not finite linear combinations of $\bG_w$.

\begin{ex} By \eqref{E:col G to s}, we have 
	\begin{align*}
		s_1 &= \sum_{p\ge1} G_{1^p} = \sum_{p\ge 1} \bG_{s_{1-p}\dotsm s_{-1}s_0}.
	\end{align*}
	Applying $\gamma$, we obtain
	\begin{align*}
		s_1 + x_1 &= \sum_{p\ge1} \bG_{s_{2-p}\dotsm s_{-1}s_0s_1}.
	\end{align*}
	Subtracting we obtain $x_1$ as an infinite linear combination of back stable Grothendieck polynomials.
\end{ex}

Let 	
\begin{align}
		\bGring := \bigoplus_{w\in S_\Z} \Z \bG_w
	\end{align} 
denote the subspace of $\bR^+$ spanned by back stable Grothendieck polynomials.
For $u,v,w\in S_+$ define $c_{uv}^w\in \Z$ by
\begin{align*}
	\G_u \G_v = \sum_{w\in S_+} c^w_{uv} \G_w.
\end{align*}

It is shown in \cite{Br} that 
\begin{align}\label{E:positivity structure constants}
(-1)^{\ell(u)+\ell(v)-\ell(w)}  c^w_{uv}\in  \Z_{\ge0}.
\end{align}

\begin{thm} \label{T:bG basis} 
	For all $u,v\in S_\Z$, there are constants $\bc^w_{uv}\in \Z$ such that 
	\begin{align}\label{E:structure constants}
		\bG_u \bG_v = \sum_{w\in S_\Z} \bc^w_{uv} \bG_w
	\end{align}
	with only finitely many $\bc^w_{uv}$ nonzero. That is, $\bGring$
	is a $\Z$-subalgebra of $\bR$.
	Moreover,
	\begin{align} \label{E:bcpos}
		(-1)^{\ell(u)+\ell(v)-\ell(w)} \bc^w_{uv} \in\Z_{\ge0}.
	\end{align}
\end{thm}
\begin{proof}
For sufficiently large $q$, we have $\bc^w_{uv} = c^{\gamma^q(w)}_{\gamma^q(u) \gamma^q(v)}$; see the more general Proposition~\ref{P:cbc} below.  Thus \eqref{E:bcpos} follows from \eqref{E:positivity structure constants}.  Applying $\eta_0$ to \eqref{E:structure constants}, we obtain
$$
G_u G_v = \sum_w \bc^w_{uv} G_w.$$
By Proposition~\ref{P:stable K stanley expansion}, we have a finite expansion $G_w= \sum_\la k^w_\la G_\la$.
Thus
\begin{align}
	G_u G_v &= \sum_{\mu,\nu,\la} k^u_\mu k^v_\nu c^\la_{\mu\nu} G_\la \\
	\sum_w \bc^w_{uv} G_w &= \sum_{w,\la} \bc^w_{uv} k^w_\la G_\la \\
\label{E:coefs1}
	\sum_{\mu,\nu} k^u_\mu k^v_\nu c^\la_{\mu\nu} &= \sum_w \bc^w_{uv} k^w_\la\qquad\text{for all $\la$.}
\end{align}
 By Proposition~\ref{P:stable K stanley expansion} the LHS of \eqref{E:coefs1} is finite and equals $0$ for all but finitely many $\la$.  
Since $(-1)^{|\mu|+|\nu|-|\la|}  c^\la_{\mu\nu} \in \Z_{\ge0}$ and $(-1)^{\ell(w)-|\la|} k_\la^w \in \Z_{\ge0}$ and we have \eqref{E:bcpos}, all terms on both sides of \eqref{E:coefs1} have the same sign $(-1)^{\ell(u) + \ell(v) - |\la|}$.  If $\bc^w_{uv}\ne0$ for infinitely many $w$, then either the RHS of \eqref{E:coefs1} is nonzero for infinitely many $\la$, or the RHS of \eqref{E:coefs1} is infinite for some $\la$, either of which is a contradiction.
\end{proof}

The following result follows from Theorem~\ref{T:coproduct}.
\begin{prop} \label{P:backstable coproduct}
We have $\Delta(\bGring) \subset \Gamma \otimes \bGring$, giving $\bGring$ the structure of a $\Gamma$-comodule.
\end{prop}

\subsection{Adjoining $\Omega$}
The ring $\bGring$ has basis $\{\bG_w \mid w \in S_\Z\}$.  The ring $\Gamma \otimes R^+$ has basis $\{G_\la \otimes \G_v \mid \la \in \P \text{ and } v \in S_{\neq 0}\}$.  By Theorem~\ref{T:coproduct}, we have a strict containment $\bGring \subsetneq \Gamma \otimes R$.  In this section, we show that the containment becomes an equality by adjoining the element $\Omega:=\Omega[\xm]$.  

\begin{ex}\label{ex:Omega}
By \eqref{E:little coproduct} we have
\begin{align*}
	\bG_{s_1} &= G_{s_1} + \G_{s_1} - G_{s_1} \G_{s_1} = G_1 + (1-G_1) x_1 \\
	x_1 &= \dfrac{1}{1-G_1} (\bG_{s_1}-\bG_{s_0}) 
	= \Omega  (\bG_{s_1} -\bG_{s_0}).
\end{align*}
where we have used \eqref{E:col G to s} for $r=1$ to obtain $1-G_1 = \Omega[-\xm]$ so that $(1-G_1)^{-1} = \Omega$.  The infinite expansion $\Omega = 1 + G_1 + G_1^2 + \cdots$ shows that $\Omega \notin \Gamma \otimes R$.
\end{ex}

The computation of Example~\ref{ex:Omega} suggests the following result.
\begin{thm} \label{T:Omega span}
	Every element of $\Gamma[\Omega] \otimes R^+$ 
	is a finite linear combination of  $\Omega^k \bG_w$ for $(k,w)\in \Z_{\ge0} \times S_\Z$.  Thus $\bGring[\Omega] \cong \Gamma[\Omega] \otimes R^+$.
\end{thm}

Define 
\begin{align}\label{E:parh}
\parh_i = \Gdiff_i - 1 = (1-x_i)A_i.
\end{align}
The operators $\parh_i$ satisfy the type $A_\Z$-braid relations. The following identity is standard.

\begin{lem} For all $w\in S_n$, we have
\begin{align}
	\Gdiff_w = \sum_{v\le w} \parh_v.
\end{align}
\end{lem}

Let $\rho^{(n)} = (n-1,n-2,\dotsc,1,0)\in\Z^n$. We have $\G_{w_0^{(n)}}=x^{\rho^{(n)}}$, where 
$w_0^{(n)}\in S_n$ is the longest element.

\begin{lem}\label{L:alt G} We have the identity
\begin{align}\label{E:alt G}
	\sum_{w\in S_n} (-1)^{\ell(w)} \G_w = \prod_{i=1}^{n-1} (1-x_i)^{n-i}.
\end{align}	
\end{lem}
\begin{proof} We have
\begin{align}
	\sum_{w\in S_n} (-1)^{\ell(w)} \G_w &= 	\sum_{w\in S_n} (-1)^{\ell(w)} \Gdiff_{w^{-1}w_0} \G_{w_0} \\
	&= (-1)^{\ell(w_0)} \left(\sum_{w\in S_n}  (-1)^{\ell(w)} \Gdiff_w \right)  x^{\rho^{(n)}}.
\end{align}
Now,
\begin{align*}
	(-1)^{\ell(w_0)} \sum_{w\in S_n} (-1)^{\ell(w)} \Gdiff_w &= (-1)^{\ell(w_0)} \sum_{w\in S_n} (-1)^{\ell(w)} 
	\sum_{v\le w} \parh_v \\
	&= (-1)^{\ell(w_0)} \sum_{v\in S_n} (-1)^{\ell(v)}\parh_v \sum_{w\ge v} (-1)^{\ell(w)-\ell(v)} \\
	&= (-1)^{\ell(w_0)} \sum_{v\in S_n} (-1)^{\ell(v)} \parh_v \delta_{v,w_0} \\
	&= \parh_{w_0}.
\end{align*} 
Using $\parh_i(x_i f) = (1-x_i)f$ for $s_i$-invariant $f$ we have
\begin{align*}
	\parh_{n-1}\dotsm \parh_2\parh_1(x^{\rho^{(n)}}) = (1-x_1)(1-x_2)\dotsm (1-x_{n-1}) x^{\rho^{(n-1)}}.
\end{align*}
Since $\prod_{i=1}^{n-1}(1-x_i)$ is $S_{n-1}$-invariant, induction completes the proof.
\end{proof}

\begin{lem}\label{L:signed GSF} We have the identity
	\begin{align*}
		\sum_{w\in S_n} (-1)^{\ell(w)} G_w = \Omega[-\xm] ^{n-1}.
	\end{align*}
\end{lem}
\begin{proof} For $N\gg n$, applying $\pit_{w_0^{(N)}}$ to \eqref{E:alt G}, by Corollary \ref{C:symmetrize G} we obtain
	\begin{align*}
		\sum_{w\in S_n} (-1)^{\ell(w)} G_w(x_1,\dotsc,x_N) = ((1-x_1)(1-x_2)\dotsm(1-x_N))^{n-1}.
	\end{align*}
	Letting $N\to\infty$ we have
	\begin{align*}
		\sum_{w\in S_n} (-1)^{\ell(w)} G_w(\xp) &= \prod_{i>0} (1-x_i)^{n-1} = \Omega[-\xp]^{n-1}.
	\end{align*} 
Now replace $\xp$ by $\xm$.
\end{proof}

We will need the left weak order $\leq_L$.  For $v,w\in S_\Z$, we have $v\le_L w$ if $\ell(v)+\ell(wv^{-1})=\ell(w)$.

\begin{lem} \label{L:first factor}
	Let $w\in S_+$ and define $J:=\{i\in \Z\mid s_i w < w\}$. Then
	\begin{align*}
		\bG_w &\in \Omega[-\xm]^{|J|} \G_w + \sum_{v<_L w} \Gamma \G_v.
	\end{align*}
\end{lem}
\begin{proof} It is not hard to check that 
	\begin{align*}
		\{u\in S_\Z\mid u*w=w\} = S_J = \langle s_i \mid i\in J\rangle.
	\end{align*}
	By \eqref{E:little coproduct}, we have
	\begin{align*}
		\bG_w\in  \left(\sum_{u\in S_J} (-1)^{\ell(u)} G_u \right) \G_w + \sum_{v<_L w} \Gamma \G_v.
	\end{align*}
	But $S_J$ is the direct product of symmetric groups, so the result follows by applying Lemma \ref{L:signed GSF} for each factor of the product.
\end{proof}

\begin{proof}[Proof of Theorem \ref{T:Omega span}]
Let $\mathcal{S}$ be the span of $\Omega^k \bG_w$ for $(k,w)\in \Z_{\ge0} \times S_\Z$.
We show that $\G_w\in \mathcal{S}$ for all $w\in S_+$ by induction on $\ell(w)$.
Clearly $\G_\id = 1=\bG_{\id}\in \mathcal{S}$. Let $\id\ne w\in S_+$.
By Lemma \ref{L:first factor} and its notation, we have
\begin{align*}
	\G_w \in \Omega^{|J|} \bG_w + \sum_{v<_Lw} \Gamma \G_v.
\end{align*}
By induction, we conclude that $\G_v\in\mathcal{S}$.
But elements of $\Gamma$ are finite linear combinations of $G_\la = \bG_{w_\la}$. By Theorem \ref{T:bG basis}, we deduce that $\G_w\in \mathcal{S}$. 

Applying $\omk$, we deduce that $\G_w\in\mathcal{S}$ for $w\in S_{\ne0}$. 
\end{proof}

\subsection{The algebra of double $K$-Stanley functions}
Let $\Gamma(x||a):= \bigoplus_{\lambda} R(a) G_\la(x||a)$ be the $R(a)$-subspace of $\hLa(x||a)$ spanned by the double $K$-Stanleys $G_\la(x||a)$.

\begin{prop}\label{P:doublewla}
\begin{enumerate}
    \item For $\la \in \P$, we have $G_\la(x/a) \in \Gamma(x||a)$.
    \item For $w \in S_\Z$, we have $G_w(x||a) \in \Gamma(x||a)$.
\end{enumerate}

\end{prop}
\begin{proof}
The first statement follows from \eqref{E:super to double}. 
For the second statement, Proposition~\ref{P:G stable double} expresses $G_w(x||a)$ as a finite $R(a)$-linear combination of the $G_v(x/a)$, and Proposition~\ref{P:stable K stanley expansion} implies that $G_v(x/a)$ is a finite $\Z$-linear combination of the $G_\la(x/a)$.  
\end{proof}

For $w \in S_\Z$ and $\la\in\Y$ let $k^w_\la(a)\in R(a)$ be defined by
\begin{align} \label{E:doublewla}
  G_w(x||a) &= \sum_\la k^w_\la(a) G_\la(x||a).
\end{align}

\begin{ex}\label{ex:kpos}
    Let $w = s_0s_2$.  We compute the expansion of $G_w(x||a)$ in terms of the Grassmannian double $K$-Stanley functions $G_\la(x||a)$ using Proposition~\ref{P:G stable double}.  First, by Remark~\ref{rem:BKSTY}, we have $G_{s_0s_2} = G_2 + G_{11} - G_{21}$.  By Proposition~\ref{P:G stable double}, and noting that $\G_{s_2}(a) \oplus \G_{s_2}(\ominus a) = 0$, we have
    \begin{align*}
    G_{s_0s_2}(x||a) &= G_{s_0}(x/a)(\G_{s_2}(a) \oplus \G_{s_2}(\ominus a)) + G_{s_0s_2}(x/a)(1 - \G_{s_2}(a) \oplus \G_{s_2}(\ominus a)) \\
    &= G_2(x/a) + G_{11}(x/a) - G_{21}(x/a).
    \end{align*}
Using the computations in Example \ref{ex:G double to super} we have
\begin{align*}
G_{s_0s_2}(x||a)&= - (1-(a_1\ominus a_0)) G_{21}(x||a)  + (1-(a_1\ominus a_0)) G_2(x||a)  \\
&\quad\,+ (1-(a_1\ominus a_0)) G_{11}(x||a)-(-(a_1\ominus a_0)) G_1(x||a) .
\end{align*}
\end{ex}

\begin{ex}\label{ex:kpos2} 
Let $w=s_2s_1$. The coefficient of $G_2(x/a)$ in $G_{s_2s_1}(x||a)$ is $(1-\G_{s_2}(\ominus(a)))(1-\G_{s_1}(a))$ and the coefficient of $G_2(x||a)$ in $G_2(x/a)$ is $(1-\G_{s_1}(a))$. Since these are the highest degree terms, we have 
\begin{align*}
k_2^{s_2s_1}= (1-\G_{s_2}(\ominus(a)))(1-\G_{s_1}(a))^2 = \dfrac{1-a_1}{1-a_2} = 1 - (a_1\ominus a_2).
\end{align*}
This agrees with Example~\ref{ex:rectbumpless}.
\end{ex}

\begin{thm}\label{T:Gammaxa}
The $R(a)$-submodule $\Gamma(x||a)\subset \hLa(x||a)$ is a $R(a)$-bialgebra.
\end{thm}

\begin{proof}
By Proposition~\ref{P:doublewla} and \eqref{E:double G coproduct}, the coproduct $\Delta(G_\la(x||a))$ is finite and thus belongs to $\Gamma(x||a) \otimes_{R(a)} \Gamma(x||a)$.  We now consider the product.  By \eqref{E:double to super}, $G_\la(x||a)$ is a finite $R(a)$-linear combination of $G_\mu(x/a)$.  The structure constants for the family $\{G_\mu(x/a)\}$ are the same as for $\{G_\mu\}\subset \Gamma$, and thus finite by Proposition~\ref{P:G stable double}.  It follows that $G_\la(x||a) G_\nu(x||a) \in \Gamma(x||a)$ for any $\la,\nu\in\P$.  
\end{proof}

Define the total ordering $\prec$ on $\Z$ by 
$$
1 \prec 2 \prec 3 \prec \cdots \prec -2 \prec -1 \prec 0.
$$

\begin{conjecture}\label{conj:kpos}For $w \in S_\Z$ and $\la \in \P$, we have
\begin{equation}\label{eq:kpos}
(-1)^{\ell(w) - |\la|} k_\la^w(a) \in \Z_{\geq0}[-(a_i \ominus a_j) \mid i \prec j].
\end{equation}

Furthermore, if $d_\nu^{\la \mu}(a)$ denotes the coproduct structure constants of the basis $G_\la(x||a)$ of $\Gamma(x||a)$, then 
\begin{equation}\label{eq:dpos}
(-1)^{|\la|+|\mu|-|\nu|} d_\nu^{\la \mu}(a) \in \Z_{\geq0}[-(a_i \ominus a_j) \mid i \prec j]
\end{equation}
\end{conjecture}
For product structure constants, see \eqref{eq:bcpos}.
For double Stanley symmetric functions an analogous positivity is proven in \cite[Theorem 4.22]{LLS:back stable}.  Similarly to \cite[(5.4)]{LSS}, we have 
\begin{equation}\label{eq:LSS}
d_\nu^{\la \mu}(a) = \sum_{\substack{w \in S_\Z \\w * w_\mu = w_\nu}} (-1)^{|\nu|-|\mu|-\ell(w)} k_\la^w(a).
\end{equation}
In particular, \eqref{eq:dpos} follows from \eqref{eq:kpos}.
\begin{ex}
Example~\ref{ex:kpos} exhibits the conjectured positivity of $k_\la^w(a)$ for $w = s_0 s_2$.  Now let $\la = (1)$, $\mu = (2)$, and $\nu = (21)$.  By \eqref{eq:LSS}, we can calculate the coproduct structure constant in two ways:
\begin{align*}
d^{1,2}_{21} &= k^{s_{-1}}_1 - k^{s_{-1}s_1}_1 = 1 - (a_1 \ominus a_0) \\
d^{2,1}_{21} &= k^{s_{-1}s_1}_{2} - k^{s_{-1}s_1s_0}_{2} = (1 - (a_1 \ominus a_0)) - 0.
\end{align*}
We have used Example~\ref{ex:kpos}, and the equality $G_{s_{-1}s_1}(x||a) = G_{s_0s_2}(x||a)$ that can be verified using Proposition~\ref{P:G stable double}.
\end{ex}

\subsection{The subring of back stable double Grothendieck polynomials}
Let $$
\bGring(x;a) = \bigoplus_{w\in S_\Z} R(a) \bG_w(x;a)
$$
be the $R(a)$-submodule spanned by the back stable double Grothendieck polynomials $\bG_w(x;a)$.
For $u,v,w\in S_\Z$ define $\bc^w_{uv}(a)\in R(a)$ by the formal expansion
\begin{align}\label{E:equivariant structure constants}
	\bG_u(x;a) \bG_v(x;a) = \sum_{w\in S_\Z}\bc^w_{uv}(a) \bG_w(x;a).
\end{align}
The existence of the expansion \eqref{E:equivariant structure constants} follows from Theorem~\ref{T:infinite basis} and Proposition~\ref{P:GKM basis}.  By Proposition \ref{P:double bG shift}, we have
\begin{align}\label{E:shift LR}
	\bc^{\gamma(w)}_{\gamma(u)\gamma(v)}(a) = \gamma  \bc^w_{uv}(a).
\end{align}

For $u,v,w\in S_+$, define $c^w_{uv}(a)\in R(a)$ by the expansion
\begin{align}\label{E:equivariant finite structure constants}
	\G_u(x;a) \G_v(x;a) = \sum_{w\in S_+}c^w_{uv}(a) \G_w(x;a).
\end{align}

\begin{prop}\label{P:cbc}
Let $u,v,w \in S_\Z$.  For $q$ such that $\gamma^q(u), \gamma^q(v), \gamma^q(w) \in S_+$, we have
$$
\bc^w_{uv}(a) = \gamma^{-q} c_{\gamma^q(u)\gamma^q(v)}^{\gamma^q(w)}(a).
$$
\end{prop}
\begin{proof}
Consider the GKM ring $\Psi_+$ consisting of the set of $\psi\in \Fun(S_+,R(T))$ such that
\begin{align}\label{E:GKM condition positive}
	1 - e^\alpha \mid \psi(s_\alpha w) - \psi(w)\qquad\text{for all $\alpha\in \Phi$ and $w\in S_+$.}
\end{align}
The ring $\Psi_+$ is an $R(T)$-subalgebra of $\Fun(S_+,R(T))$ such that $\Psi_+ = \prod_{z\in S_+} R(T) \psi^z_+$, where $\psi^z_+$ is the restriction of $\psi^z \in \Psi$ to $S_+ \subset S_\Z$.  Using Proposition~\ref{P:locs}, we see that an analogue of Theorem~\ref{T:infinite basis} holds $\Psi_+$ replacing $\Psi$ and $\{\G_z \mid z \in S_+\}$ replacing $\{\bG_z \mid z \in S_\Z\}$.  By triangularity \eqref{E:loc support}, the coefficient of $\psi^w$ in the product $\psi^u \psi^v$ can be obtained from a finite computation that only involves the values of the values of various $\psi^z$-s on the lower order ideal in $S_\Z$ generated by $u,v,w$.  When $\gamma^q(u), \gamma^q(v), \gamma^q(w) \in S_+$, this lower order ideal is contained in $S_+$, so $c_{\gamma^q(u)\gamma^q(v)}^{\gamma^q(w)}(a) = \bc_{\gamma^q(u)\gamma^q(v)}^{\gamma^q(w)}(a)$.  The result then follows from \eqref{E:shift LR}.
\end{proof}

It is shown in \cite{AGM} that 
$
(-1)^{\ell(u)+\ell(v) - \ell(w)} c^w_{uv}(a) \in \Z_{\geq0}[-(a_j \ominus a_i) \mid i<j].
$
\begin{ex}
We have $\G_{s_1}(x;a) = x_1 \ominus a_1$ and $\G_{s_2s_1}(x;a) = (x_1 \ominus a_1)(x_1 \ominus a_2)$.  Thus
\begin{align*}
\G_{s_1}(x;a)^2 &= (x_1 \ominus a_1)^2 \\
&= (x_1 \ominus a_1)( (x_1 \ominus a_2) \oplus (a_2 \ominus a_1)) \\
&= (1-(a_2 \ominus a_1))\G_{s_2s_1}(x;a) + (a_2 \ominus a_1)\G_{s_1}(x;a).
\end{align*}
We have $(1-(a_2 \ominus a_1)) \in  \Z_{\geq0}[-(a_j \ominus a_i) \mid i<j]$ and $(-1)^1(a_2 \ominus a_1) \in  \Z_{\geq0}[-(a_j \ominus a_i) \mid i<j]$.
\end{ex}

It follows from Proposition~\ref{P:cbc} that we have
\begin{equation}\label{eq:bcpos}
(-1)^{\ell(u)+\ell(v) - \ell(w)} \bc^w_{uv}(a) \in \Z_{\geq0}[-(a_j \ominus a_i) \mid i < j].
\end{equation}

\begin{conjecture}\label{conj:bG double basis}
For fixed $u,v \in S_\Z$, only finitely many $\bc^w_{uv}$ nonzero. That is, 
\begin{equation}\label{eq:bGringxa}
\mbox{$\bGring(x;a)$ is a $R(a)$-subalgebra of $\bR(x;a)$.}
\end{equation}
\end{conjecture}

\begin{prop} \label{P:bG double basis} 
Suppose that the positivity \eqref{eq:kpos} holds.  Then Conjecture~\ref{conj:bG double basis} holds.
\end{prop}
\begin{proof}
The argument is the same as the proof of Theorem~\ref{T:bG basis}.  Applying $\eta_a$ to \eqref{E:equivariant structure constants}, we obtain
$$
G_u(x||a) G_v(x||a) = \sum_w \bc^w_{uv}(a) G_w(x||a).$$

By equation \eqref{E:doublewla} we have
\begin{align}
	G_u(x||a) G_v(x||a) &= \sum_{\mu,\nu,\la} k^u_\mu(a) k^v_\nu(a) c^\la_{\mu\nu}(a) G_\la(x||a) \\
	\sum_w \bc^w_{uv}(a) G_w(x||a) &= \sum_{w,\la} \bc^w_{uv} k^w_\la(a) G_\la(x||a) \\
\label{E:coefs}
	\sum_{\mu,\nu} k^u_\mu(a) k^v_\nu(a) c^\la_{\mu\nu}(a) &= \sum_w \bc^w_{uv}(a) k^w_\la (a)\qquad\text{for all $\la$.}
\end{align}
By Proposition~\ref{P:doublewla} and Theorem~\ref{T:Gammaxa}, the left side of \eqref{E:coefs} is finite and equals $0$ for all but finitely many $\la$. Since $c^\la_{\mu\nu}(a) \in (-1)^{|\mu|+|\nu|-|\la|} \Z_{\geq0}[a_i \ominus a_j \mid i < j]$, and by assumption $ k_\la^w(a) \in (-1)^{\ell(w)-|\la|} \Z_{\geq0}[a_i \ominus a_j \mid i < j]$, and we have \eqref{eq:bcpos}, all terms on both sides of \eqref{E:coefs} belong to $(-1)^{\ell(u) + \ell(v) - |\la|} \Z_{\geq0}[a_i \ominus a_j \mid i < j]$.  If $\bc^w_{uv}(a)\ne0$ for infinitely many $w$, then either the RHS of \eqref{E:coefs} is nonzero for infinitely many $\la$, or the RHS of \eqref{E:coefs} is infinite for some $\la$, either of which is a contradiction.
%
\end{proof}

The following follows from Theorem~\ref{T:double coproduct}.
\begin{prop} The algebra $\Gamma(x||a)[\Omega(x/a)]$ is a Hopf algebra over $R(a)$.  The $R(a)$-algebra $\bGring(x;a)$ is a $\Gamma(x||a)[\Omega(x/a)]$ Hopf-comodule. 
\end{prop}

\subsection{Adjoining $\Omega(x/a)$}
In this subsection, we assume that Conjecture~\ref{conj:bG double basis} holds, that is, $\bGring(x;a)$ is a $R(a)$-subalgebra of $\bR(x;a)$.  We then compare the two subalgebras $\Gamma(x||a)\otimes_{R(a)} R(x;a)^+$ and $\bGring(x;a)$ by adjoing the element $\Omega(x/a)$.

\begin{rem}
The element $\Omega(x/a)$ has the following geometric interpretation: it is the class $[\mathcal{L}_{\Lambda_0}]$ (in an appropriate equivariant $K$-group of $\Flb$) of the line bundle with fundamental weight $\Lambda_0$.  Indeed, in general the class $[\mathcal{O}_{X_{s_i}}]$ is equal to $1 - [\mathcal{L}_{-\Lambda_i}]$.  Thus $[\mathcal{L}_{\Lambda_0}] = (1-\bG_{s_0}(x;a))^{-1}=\Omega(x/a)$.
\end{rem}

\begin{prop} \label{P:invert equivariant omega}
Assume \eqref{eq:bGringxa}.  Then
every element of $\Gamma(x||a)[\Omega(x/a)] \otimes_{R(a)} R(x;a)^+$ is a finite $R(a)$-linear combination of $\Omega(x/a)^k \bG_w(x;a)$ for $(k,w)\in \Z_{\ge0}\times S_\Z$. Thus $\bGring(x;a)[\Omega] \cong \Gamma(x||a)[\Omega(x/a)] \otimes_{R(a)} R(x;a)^+$.
\end{prop}
\begin{proof} 
Let $R'$ be the $R(a)$-span of $\Omega(x/a)^k \bG_w(x;a)$ for $(k,w)\in \Z_{\ge0}\times S_\Z$. Since $\Gamma(x||a)$ is spanned by $G_\la(x||a)=\bG_{w_\la}(x;a)$ and assuming \eqref{eq:bGringxa} the proof reduces to showing that $\G_w(x;a)\in R'$ for $w\in S_{\ne0}$. Using $\omk$ one may further reduce to $w\in S_+$ which we now assume.

Arguing as in Lemma \ref{L:first factor} and using its notation, let $J=\{i\mid s_i w < w\}$ and let $S_J = \langle s_j\mid j\in J\rangle$. By \eqref{E:double G little coproduct} we have
\begin{align}
	\bG_w(x;a) \in\G_w(x;a) \sum_{u\in S_J} (-1)^{\ell(u)} G_u(x||a) + \sum_{v <_L w} \Gamma(x||a) \G_v.
\end{align}
By Proposition \ref{P:G stable double} we have
\begin{align*}
	\sum_{u\in S_J} (-1)^{\ell(u)} G_u(x||a) &= 
	\sum_{u\in S_J} (-1)^{\ell(u)} \sum_{u_1*u_2*u_3=u} (-1)^{\ell(u_1)+\ell(u_2)+\ell(u_3)-\ell(u)} \times  \\
	& \qquad G_{u_1^{-1}}(\ominus(a)) G_{u_2}(x/a) G_{u_3}(a) \\
	&= \left(\sum_{u_1\in S_J} (-1)^{\ell(u_1)} G_{u_1^{-1}}(\ominus(a))\right) \left( \sum_{u_2\in S_J} (-1)^{\ell(u_2)} G_{u_2}(x/a)\right) \\
	&\quad\times \left( \sum_{u_3\in S_J} (-1)^{\ell(u_3)} \G_{u_3}(a) \right) \\
	&= \sum_{u_2\in S_J} (-1)^{\ell(u_2)} G_{u_2}(x/a) \\
	&= \Omega(a/x)^{|J|}
\end{align*}
where in the last step we have applied the superization of Lemma \ref{L:signed GSF}. Therefore $\bG_w(x;a) \in \Omega(a/x)^{|J|} \G_w(x;a) + \sum_{v<_L w} \Gamma(x||a) \G_v(x;a)$. The Proposition holds by induction.
\end{proof}

\section{Determinantal formulae for $G_\la(x||a)$} \label{S:detformula}
In this section, we recover the equality of the Grassmannian double $K$-Stanley functions as Grassmannian back stable Grothendieck polynomials, with the determinantal formulae in the literature.

\begin{prop}\label{P:G double row col} For $r \geq 1$, we have
	\begin{align*}
		G_r(x||a) &= \gamma_a^{r-1} G_r(x/a) \\
		G_{1^r}(x||a) &= \gamma_a^{1-r} G_{1^r}(x/a).
	\end{align*}
\end{prop}
\begin{proof} We prove the second formula as the first follows from it by 
	Propositions \ref{P:double G negation Grassmannian}, \ref{P:omk G super},
	and \ref{P:omk shift}.
		
	Let $c_r := w_{(1^r)} = s_{1-r} \dotsm s_{-1} s_0$ 
	and $c_{r/p}=w_{(1^r)/(1^p)}=s_{1-r}\dotsm s_{-p}$ for $
	0\le p\le r$. The $0$-Hecke factorizations of $c_r$ with left factor in $S_{\ne0}$ are given by $c_{r/p}*c_p$  for $1\le p\le r$ and	
	$c_{r/(p-1)}*c_p$ for $2\le p\le r$. Letting
	$A=(a_0,a_{-1},\dotsc,a_{2-r})$ we have
	\begin{align*}
		\G_{c_{r/p}^{-1}}(\ominus a) &= \G_{s_p s_{p+1}\dotsm s_{r-1}}[A] = G_{r-p}[A]	
	\end{align*}
By \eqref{E:G back stable triple} we have
	\begin{align*}
		G_{(1^r)}(x||a) 
		&= 
		\sum_{p=1}^r \G_{c_{r/p}^{-1}}(\ominus a) G_{(1^p)}(x/a) - \sum_{p=2}^r \G_{c_{r/(p-1)}^{-1}}(\ominus a) G_{(1^p)}(x/a) \\
		&= \sum_{p=1}^r (G_{r-p}[A]-G_{r-p+1}[A])  G_{(1^p)}(x/a)
	\end{align*}
	using the fact that $G_{1^r}[A]=0$ because $A$ consists of $r-1$ variables. We have $\gamma_a^{1-r}(p_1(x||a)) = p_1(x||a)  + A$ and
	\begin{align*}
		\gamma_a^{1-r} G_{1^r}(x/a)  &= G_{1^r}[x_{\le0}-a_{\le0} + A].
	\end{align*}
	By \cite{Bu} we have
	\begin{align*}
		\Delta(G_{1^r}) &= \sum_{p=0}^r G_{1^p}  \otimes G_{1^{r-p}} - \sum_{p=1}^{r} G_{1^p} \otimes G_{1^{r+1-p}}.
	\end{align*}
	Since we are using the coproduct such that $p_r(x/a)$ are primitive we may superize this formula,  replacing $K$-Stanleys by superized $K$-Stanleys. We obtain
	\begin{align*}
		\gamma_a^{1-r} G_{1^r}(x/a) &= \sum_{p=0}^r G_{1^p}(x/a)  G_{1^{r-p}}[A] - \sum_{p=1}^{r} G_{1^p}(x/a)  G_{1^{r+1-p}}[A] \\
		&= \sum_{p=1}^r G_{1^p}(x/a)   G_{1^{r-p}}[A] - \sum_{p=1}^{r} G_{1^p}(x/a) G_{1^{r+1-p}}[A]
	\end{align*}
as required.
\end{proof}

\begin{lem} If $s_i^a g = g$ then $\Gdiff_i^{a,\anti} (fg) = \Gdiff_i^{a,\anti}(f) g$.
\end{lem}

\begin{lem} \label{L:Gdiff on shifted row G}
	For $p,q\in \Z$ and $r\in\Z_{>0}$ we have
	\begin{align}
		\Gdiff_p^{a,\anti} \gamma_a^q G_r(x/a) = 
		\begin{cases}
			\gamma_a^q G_r(x/a) & \text{if $p\ne q$} \\
			\gamma_a^{q-1} G_{r-1}(x/a)&\text{if $p=q$.}
			\end{cases}
	\end{align}
\end{lem}
\begin{proof}
We have
\begin{align*}
	\Gdiff_p^{a,\anti} \gamma_a^q G_r(x/a) &=
	\gamma_a^q \Gdiff_{p-q}^{a,\anti} G_r(x/a).
\end{align*}
If $p\ne q$ then $s_{p-q}^a G_r(x/a)=G_r(x/a)$ and 
$\gamma_a^q \Gdiff_{p-q}^{a,\anti} G_r(x/a)=\gamma_a^q G_r(x/a)$. If $p=q$ then
\begin{align*}
	\gamma_a^q \Gdiff_0^{a,\anti} G_r(x/a) &= 
		\gamma_a^q \Gdiff_0^{a,\anti} \gamma_a^{1-r} \gamma_a^{r-1}G_r(x/a) \\
		&= \gamma_a^{q+1-r} \Gdiff_{r-1}^{a,\anti} G_r(x||a) \\
		&= 	\gamma_a^{q+1-r} G_{r-1}(x||a) \\
		&= \gamma_a^{q-1} G_{r-1}(x/a). \qedhere
\end{align*}
\end{proof}

In the following, we denote by $\binom{n}{k}$ the signed binomial coefficient, the coefficient of $x^k$ in the power series $(1+x)^n$ for any integer $n$.
In particular for $n\ge0$,  $\binom{-n}{k} =(-1)^k \binom{k+n-1}{n-1}$.

Formulas \eqref{E:G det} and \eqref{E:e det} are due to \cite{HIMN} and \cite{And} respectively.

\begin{prop}\label{P:G det}
	With $\ell=\ell(\la)$,
\begin{align}\label{E:G det}
G_\la(x||a) &= \det \left(\gamma_a^{\la_i-i} \sum_{k\ge0} (-1)^k \binom{i-j}{k}  G_{\la_i-i+j+k}(x/a)\right)_{1\le i,j\le \ell} \\
\label{E:e det}
	G_\la(x||a) &= \det \left(\gamma_a^{i-\la'_i} \sum_{p\ge0} (-1)^p \binom{p+\la'_i-1}{\la'_i-1} e_{\la'_i-i+j+p}(x/a) \right)_{1\le i,j\le \la_1}
\end{align}
\end{prop}
\begin{proof} 
Assuming \eqref{E:G det} we derive \eqref{E:e det} as follows. Applying $\omk$ to \eqref{E:G det} for $G_{\la'}(x||a)$, by Propositions \ref{P:double G negation Grassmannian}, \ref{P:omk shift}, \ref{P:omega G partition}, and \ref{P:G one row column} we have 
\begin{align*}
	G_\la(x||a) &= \omk \det \left(\gamma_a^{\la'_i-i} \sum_{k\ge0} (-1)^k \binom{i-j}{k}  G_{\la'_i-i+j+k}(x/a)\right)\\
	&= \det \left(\gamma_a^{i-\la'_i} \sum_{k\ge0} (-1)^k \binom{i-j}{k}  G_{1^{\la'_i-i+j+k}}(x/a)\right)\\
	&= \det \left(\gamma_a^{i-\la'_i} \sum_{k\ge0} (-1)^k \binom{i-j}{k} \sum_{m \ge0} (-1)^m \binom{m+\la'_i-i+j+k-1}{m}  e_{m+\la'_i-i+j+k}(x/a)\right)\\
	&= \det \left(\gamma_a^{i-\la'_i} \sum_{p\ge0} (-1)^p  e_{\la'_i-i+j+p}(x/a) \sum_{k+m=p} \binom{i-j}{k}  \binom{p+\la'_i-i+j-1}{m}  \right)\\
&=	\det \left(\gamma_a^{i-\la'_i} \sum_{p\ge0} (-1)^p  e_{\la'_i-i+j+p}(x/a) \binom{p+\la'_i-1}{p} \right).
\end{align*}

To prove \eqref{E:G det} let
\begin{align*}
	L^{(m)}_{ij} &= \gamma_a^{m-i} \sum_{k\ge0}(-1)^k \binom{i-j}{k} G_{m-i+j+k}(x/a) 
\end{align*}
so that $H_\la := \det L^{(\la_i)}_{ij}$ is the right hand side of \eqref{E:G det}.
By induction it suffices to show that 
\begin{align}\label{E:det suffices}
\Gdiff_p^{a,\anti} H_\la = \begin{cases} H_\mu & \mbox{if $p = \la_i - i$ and $\la/\mu$ is a box in the $i$-th row} \\
H_\la &  \mbox{otherwise.} \\
\end{cases}
\end{align}
By Lemma \ref{L:Gdiff on shifted row G}, every entry $L_{ij}^{(\la_i)}$ in the $i$-th row is $s_p^a$-invariant for all $p$ except $p=\la_i-i$.

Suppose $p\ne \la_i-i$ for all $i$. We have $s_p^a(H_\la)=H_\la$ and $\Gdiff_p^{a,\anti} H_\la=H_\la$ which agrees with \eqref{E:det suffices}.

Otherwise let $p = \la_i-i$ for some $i$. Since this $p$ is necessarily unique it follows that $\Gdiff_p^{a,\anti}$ fixes all of the entries in the rows of the determinant other than the $i$-th. By Lemma \ref{L:Gdiff on shifted row G} we have
\begin{align}\label{E:special box}
\Gdiff_{\la_i-i}^{a,\anti} L_{ij}^{(\la_i)} = L_{ij}^{(\la_i - 1)}.
\end{align}	
If $\la$ has a removable box in the $i$-th row and $\mu$ is the partition obtained by removing a box in row $i$ then
we have $L_{qj}^{(\mu_q)} = L_{qj}^{(\la_q)}$ for $q\ne i$. By multilinearity of the determinant, \eqref{E:det suffices} follows from \eqref{E:special box}.

Otherwise, let $\ell := \la_{i+1} = \la_i$.  Define $F_m:= \gamma_a^{-1} G_m(x/a)$.  Then by \eqref{E:gammad} with $d = 1$, and using $G_r[a_0] = a_0^r$, we have
$$
G_m(x/a) = \frac{1}{1-a_0}(F_m - a_0 F_{m-1}).
$$
Thus
\begin{align*}
\gamma_a^{i-\ell} L_{ij}^{(\ell)} &= \sum_{k\ge0}(-1)^k \binom{i-j}{k} G_{\ell-i+j+k}(x/a)= \frac{1}{1-a_0} \sum_{k\ge0}(-1)^k \binom{i-j}{k}(F_{\ell-i+j+k}-a_0 F_{\ell-i+j+k-1})\\
\gamma_a^{i-\ell} L_{i+1,j}^{(\ell)} &= \gamma_a^{-1} \sum_{k\ge0}(-1)^k \binom{i+1-j}{k} G_{\ell-i+j+k-1}(x/a) =  \sum_{k\ge0}(-1)^k \binom{i+1-j}{k} F_{\ell-i+j+k-1}
\end{align*}
Then $\gamma_a^{i-\ell} L_{ij}^{(\ell)}  + \frac{1}{1-a_0}\gamma_a^{i-\ell} L_{i+1,j}^{(\ell)}$ equals
\begin{align*}
& \frac{1}{1-a_0} \sum_{k\ge0}(-1)^k \left( \binom{i-j}{k}(F_{\ell-i+j+k}-a_0 F_{\ell-i+j+k-1})+\binom{i+1-j}{k} F_{\ell-i+j+k-1}\right)
\\
&=\frac{1}{1-a_0}  \sum_{k\ge0}(-1)^k \left( \binom{i-j}{k} +a_0 \binom{i-j}{k+1} -\binom{i+1-j}{k+1} \right) F_{\ell-i+j+k} + \frac{ (-a_0 F_{\ell-i+j-1}+ F_{\ell-i+j-1}) }{1-a_0}\\
&= F_{\ell-i+j-1} + \sum_{k\ge0}(-1)^{k+1} \binom{i-j}{k+1}  F_{\ell-i+j+k}\\
&= \gamma_a^{i-\ell} L_{ij}^{(\ell-1)}.
\end{align*}
Applying $\gamma_a^{\ell-i}$, we obtain $L_{ij}^{(\ell-1)} = L_{ij}^{(\ell)}  + (\gamma_a^{\ell-i} \frac{1}{1-a_0}\gamma_a^{i-\ell}) L_{i+1,j}^{(\ell)}$, showing that $\Gdiff_p^{a,\anti} H_\la= H_\la$ in this case.
\end{proof}

\begin{ex} 
Let $\la=(1,1)$. We have
\begin{align*}
	G_{11}(x||a) &= \begin{vmatrix}
		G_1 & G_2+G_3+G_4+\dotsm \\
		\gamma_a^{-1}(1-G_1) & \gamma_a^{-1}(G_1) .
	\end{vmatrix}
\end{align*}
Note that $0$ is not the residue of a corner box. We have
\begin{align*}
	\Gdiff_0^{a,\anti} G_{11}(x||a) &= 
	\begin{vmatrix}
		1  & \gamma_a^{-1}(G_1+G_2+G_3+\dotsm) \\
		\gamma_a^{-1}(1-G_1) & \gamma_a^{-1}(G_1) .
	\end{vmatrix}
\end{align*}
The proof of Proposition~\ref{P:G det} shows that 
\begin{align*}
1 &= G_1 + \frac{1}{1-a_0}\gamma_a^{-1}(1-G_1)  \\
\gamma_a^{-1}(G_1+G_2+G_3+\dotsm)  &= (G_2+G_3+G_4+\dotsm)+ \frac{1}{1-a_0}\gamma_a^{-1}(G_1) 
\end{align*}
so $\Gdiff_0^{a,\anti} G_{11}(x||a)  = G_{11}(x||a)$.
\end{ex}

\begin{ex} Let $\la=(1,1)$. We check \eqref{E:G det}.  First, suppose that $a=0$. Using e.g. \cite[Thm. 5.4]{Bu} we have $G_r G_1 = G_{r+1} + G_{r,1} - G_{r+1,1}$ for all $r\ge 1$.
	\begin{align*}
			\det \begin{pmatrix} G_1 & G_2+G_3+\dotsm\\
	1 - G_1 & G_1
	\end{pmatrix} &= G_1^2 - (1-G_1)(G_2+G_3+\dotsm) \\
	&= (G_2 + G_{11} - G_{21}) - (G_2+G_3+\dotsm)\\
	&+ ((G_3+G_{21}-G_{31})+(G_4+G_{31}-G_{41}) +\dotsm) \\
	&= G_{11}.
	\end{align*}
Now consider the general double case. By e.g. \cite[\S 6]{Bu} and superizing we have
\begin{align*}
\Delta(G_r(x/a)) = \sum_{p=0}^r G_p(x/a) \otimes G_{r-p}(x/a) - \sum_{p=1}^r G_p(x/a) \otimes G_{r+1-p}(x/a).
\end{align*}
For $d\ge 0$ and letting $A=a_0+a_{-1}+\dotsm+a_{1-d}$ we have
\begin{align}\label{E:gammad}
	\gamma_a^{-d}(G_r(x/a)) &= \sum_{p=0}^r G_p(x/a) G_{r-p}[A] - \sum_{p=1}^r G_p(x/a) G_{r+1-p}[A].
\end{align}
For $d=1$ we have
\begin{align*}
	\gamma_a^{-1}(G_1(x/a)) &= G_1[a_0] + G_1(x/a) - G_1(x/a) G_1[a_0] = a_0 + (1-a_0)G_1(x/a).
\end{align*}
\begin{align*}
&\det \begin{pmatrix}
G_1(x/a) & (G_2(x/a)+G_3(x/a)+\dotsm) \\
\gamma_a^{-1}(1-G_1(x/a))& \gamma_a^{-1}(G_1(x/a))
\end{pmatrix} \\
&= G_1(x/a) ( a_0 + (1-a_0)G_1(x/a)) \\
&- (1-a_0-(1-a_0)G_1(x/a))) (G_2(x/a)+G_3(x/a)+\dotsm) \\
&= a_0G_1(x/a) + (1-a_0) G_1(x/a)^2\\& -(1-a_0)(1-G_1(x/a))(G_2(x/a)+G_3(x/a)+\dotsm) \\
&= a_0 G_1(x/a) + (1-a_0) G_{11}(x/a).
\end{align*}
Let us compare with Proposition \ref{P:G double row col}.
\begin{align*}
	\Delta(G_{1^r}(x/a)) &= \sum_{p=0}^r G_{1^p}(x/a) \otimes G_{1^{r-p}}(x/a) - \sum_{p=1}^r G_{1^p}(x/a) \otimes G_{1^{r+1-p}}(x/a).
\end{align*}
\begin{align*}
G_{11}(x||a) &= \gamma_a^{-1} G_{11}(x/a) \\
&= G_1(x/a) G_1[a_0] + G_{11}(x/a) - G_1(x/a) G_{11}[a_0] - G_{11}(x/a) G_1[a_0] \\
&= a_0 G_1(x/a) + (1-a_0) G_{11}(x/a)
\end{align*}
which agrees with the determinant.
\end{ex}

\section{Back stable double Grothendieck polynomials via degeneracy loci}
\label{S:degeneracy loci}
By \cite[Theorem 2.1]{Bu:K quiver} it is known that double Grothendieck polynomials are the universal formulas for $K$-classes of 
certain degeneracy loci based on quivers. Since back stable double Grothendiecks are certain limits of double Grothendieck polynomials, they can also be computed by universal quiver locus formulas. Following the suggestion of Buch \cite{Buch:email} we apply such a quiver formula and recover one of our formulas for back stable double Grothendieck polynomials.

Without loss of generality we take $w\in S_+$ and then $w\in S_n$.
Let $m\ge0$ be a nonnegative integer and let $x_-^{(m)}=(x_{1-m},\dotsc,x_{-1},x_0)$ and $a_-^{(m)}=(a_{1-m},\dotsc,a_{-1},a_0)$ be sets of $m$ variables with nonpositive indices.
Let $\gamma^m(w)$ be the $m$-fold forward shift of $w$. By definition
\begin{align}
  \bG_w = \lim_{m\to\infty} \G_{\gamma^m(w)}(x_-^{(m)},x;a_-^{(m)},a).
\end{align}
We now apply a formula for double Grothendieck polynomials which is a variant of 
\cite[Theorem 4]{BKTY} but whose overall form more closely follows \cite[Theorem 4]{BKTY:Schub}.
We use a set of $m$ variables followed by $n-1$ sets of one variable each.
For the $x$ variables we use $x_-^{(m)}$ and then $x_1$, $x_2$, $\dotsc$, $x_{n-1}$.
The equivariant $a$ variables are similarly grouped.
This grouping of the variables is \emph{compatible} with $\gamma^m(w)$ and $\gamma^m(w^{-1})=(\gamma^m(w))^{-1}$ in the 
language of \cite{BKTY:Schub}. Our variant of \cite[Theorem 4]{BKTY} states that 
\begin{align}\label{E:G from quiver}
\G_{\gamma^m(w)}(x;a) &= \sum_{\lad} (-1)^{|\lad|-\ell(w)} c^w_\lad G_{\la^{(1-n)}}[-a_{n-1}] \dotsm G_{\la^{(-1)}}[-a_1] 
G_{\la^{(0)}}[x_-^{(m)}/a_-^{(m)}] \times \\ \notag
&\qquad G_{\la^{(1)}}[x_1]\dotsm G_{\la^{(n-1)}}[x_{n-1}]
\end{align}
where $\lad=(\la^{(j)}\mid 1-n\le j\le n-1)$ runs over tuples of partitions and $c^w_\lad$ is the number of tableau tuples 
$T^\bullet = (T^{(1-n)},\dotsc,T^{(n-1)})$ such that $T^{(j)}$ is a \emph{decreasing} tableau (one whose rows strictly decrease from left to right and whose columns strictly decrease from top to bottom) of shape $\la^{(j)}$ for $1-n\le j\le n-1$ such that the juxtaposition of the column-reading words $T^{(1-n)}\dotsm T^{(0)}\dotsm T^{(n-1)}$ is 0-Hecke equivalent to $w$, and for $1\le i\le n-1$, the entries of $T^{(\pm i)}$ are at least $i$.

\begin{rem} To compare the formulas it is better to look at \cite[Theorem 4]{BKTY:Schub}.
Here we use decreasing tableaux instead of increasing, which allows us to avoid transposing the shapes of 
the tableaux. We are using a different form of double Grothendieck polynomial than \cite{BKTY}; see Remark \ref{R:double Groth convention}.
\end{rem}

We observe that for a single variable $z$, $G_\la[z]=0$ unless $\la$ is a single row, say, $(r)$, in which case
$G_r[z] = z^r$. Similarly, $G_\la[-z] = 0$ unless $\la$ is a single column, say $(1^r)$, in which case one may show that
$G_{1^r}[-z] = (\ominus z)^r$. Therefore we may assume that $\la^{(j)}$ is a single row for $1\le j\le n-1$ and a single 
column for $1-n \le j\le -1$. Next we observe that the Fomin-Kirillov formula for Grothendieck polynomials can be rewritten as follows. For $z\in S_n$ we have
\begin{align*}
  \G_z(x_1,\dotsc,x_{n-1}) = \sum_{T^\bullet} (-1)^{\ell(T^\bullet)-\ell(z)} \prod_{j=1}^{n-1} x_j^{\ell(T^{(j)})}
\end{align*}
where the sum runs over $T^\bullet = (T^{(1)},\dotsc,T^{(n-1)})$ where $T^{(j)}$ is a single row decreasing tableau whose entries are at least $j$, and such that $T^{(1)}\dotsm T^{(n-1)}$ is Hecke equivalent to $z$.

Applying inverses and evaluating at $(\ominus(a_1),\dotsm,\ominus(a_{n-1}))$, for $u\in S_n$ we see that
\begin{align*}
\G_{u^{-1}}(\ominus(a_1),\dotsm,\ominus(a_{n-1})) = \sum_{U^\bullet} (-1)^{\ell(U^\bullet)-\ell(u)} \prod_{j=1}^{n-1} (\ominus a_j)^{\ell(U^{(j)})} \end{align*}
where $U^\bullet=(U^{(1)},\dotsc,U^{(n-1)})$ runs over tuples with $U^{(j)}$ a decreasing tableau of single column shape with entries at least $j$, with $U^{(n-1)}\dotsm U^{(1)}$ Hecke equivalent to $u$.

Finally we use Remark \ref{rem:BKSTY} for the expansion of a super $K$-Stanley function $G_v(x_-^{(m)}/a_-^{(m)})$ into 
Grassmannian super $K$-Stanleys $G_\la(x_-^{(m)}/a_-^{(m)})$.

Combining the above, we see that \eqref{E:G from quiver} becomes 
\begin{align*}
  \G_{\gamma^m(w)}(x_-^{(m)},x; a_-^{(m)},a) = 
  \sum_{u*v*z=w} (-1)^{\ell(u)+\ell(v)+\ell(z)-\ell(w)} \G_{u^{-1}}(\ominus(a)) G_v(x_-^{(m)}/a_-^{(m)}) \G_z(x).
\end{align*}
Sending $m$ to infinity we recover the triple coproduct formula Proposition \ref{P:G back stable triple}.

\section{Further directions}\label{S:further}
\subsection{Ideal sheaf basis}
Define 
\begin{align*}
	\G^\partial_w(x;a) := \sum_{v\le w} (-1)^{\ell(v)} \G_v(x;a).
\end{align*}
These polynomials represent the ideal sheaf basis in the equivariant cohomology of finite flag varieties.  They can be generated by the operators $\Gdiff_i^x - 1$, with $\G^\partial_{w_0}(x;a) :=  \prod_{i+j\le n} (x_i \ominus a_j -1)$.  The back stable double ideal sheaf basis is given by
\begin{align*}
	\bG^\partial_w(x;a) := \sum_{v\le w} (-1)^{\ell(v)} \bG_v(x;a).
\end{align*}
Clearly, $\{\bG^\partial_w(x;a) \mid w \in S_\Z\}$ and $\{\bG_w(x;a) \mid w \in S_\Z\}$ span the same $R(a)$-module.  We expect most of the results of this work to have analogues for the polynomials $\{\bG^\partial_w(x;a) \mid w \in S_\Z\}$.


\subsection{$K$-Peterson subalgebra}
In \cite{LLS:back stable}, we defined a commutative subalgebra of the infinite nilHecke algebra and showed that it provided a model for the equivariant homology of the infinite Grassmannian.  This commutative subalgebra is an analogue of the \emph{Peterson subalgebra} of an affine nilHecke algebra.  We expect that this construction can be extended to $K$-homology of the infinite Grassmannian.  The analogous \emph{$K$-Peterson subalgebra} for modeling the equivariant $K$-homology of the affine Grassmannian is constructed in \cite{LSS}.  We remark that recent work of Kato \cite{Kt} (see also \cite{LLMS}) relate the $K$-homology of the affine Grassmannian with quantum $K$-theory of flag varieties.

\subsection{Relation to back stable Schubert polynomials}
Lenart \cite{Len3} showed that Grothendieck polynomials expand into Schubert polynomials with alternating coefficients, and gave a combinatorial interpretation of the coefficients.  It would be interesting to study the expansion of back stable (double) Grothendieck polynomials $\bG_w$ into back stable (double) Schubert polynomials $\bS_w$, and the relation to our other expansion formulae such as the coproduct formula (Theorem~\ref{T:coproduct}).

\subsection{Relation to $K$-theory affine Schubert calculus}
Affine Schubert calculus is one of our main motivations to study back stable Schubert calculus.  We expect many interesting relations between these subjects.  In particular, we expect that a wealth of combinatorics can be found in the expansion coefficients of Schubert classes of the infinite flag variety (or infinite Grassmannian) in terms of Schubert classes of the affine flag variety (or affine Grassmannian) \cite{LLS:coprod, LSS,Mor}.  In cohomology, these expansion coefficients are known as \emph{$k$-branching coefficients}.

\appendix
\section{Grothendieck Inversion}

Let $J\subset \Z\setminus\{0\}$ and let $S_J$ be the subgroup of $S_\Z$
generated by $s_j$ for $j\in J$.
Let $v,w\in S_{\Z}$ be such that $ S_J v= S_Jw$.
Say $v \leleft{J} w$ if there is a $u\in S_J$ such that $u*v=w$.

\begin{prop} \label{P:left Groth inversion}
Let $W'$ be a $\leleft{J}$ interval in $S_{\Z}$ and let
$\{f_w\mid w\in W'\}$ and $\{g_w\mid w\in W'\}$ be families of elements.
Then
\begin{align}\label{E:inverse change}
	g_w &= \sum_{u*v=w} (-1)^{\ell(u)+\ell(v)-\ell(w)} \G_u(a) f_v
\end{align}
holds if and only if
\begin{align}\label{E:change}
	f_w = \sum_{u*v=w} (-1)^{\ell(u)+\ell(v)-\ell(w)} \G_{u^{-1}}(\ominus a) g_v
\end{align}
does.
\end{prop}

This is equivalent to the following. Let $W'$ be as above.
Define the $W' \times W'$ matrices
\begin{align}
	A_{vw} &= \sum_{\substack{u\in S_J\\ u*v=w}} (-1)^{\ell(u)+\ell(v)-\ell(w)} \G_u(a) \\
	B_{vw} &= \sum_{\substack{u\in S_J\\ u*v=w}} (-1)^{\ell(u)+\ell(v)-\ell(w)} \G_{u^{-1}}(\ominus(a)).
\end{align}
Then $A$ and $B$ are inverse.

Let 
\begin{align}
	G &= \sum_{u\in S_J} (-1)^{\ell(u)} \G_u(a) \Gdiff_u \\
	H &= \sum_{u\in S_J} (-1)^{\ell(u)} \G_{u^{-1}}(\ominus(a)) \Gdiff_u.
\end{align}
Here we are working in the $0$-Hecke algebra with coefficients in the field $\Q(a_i\mid i\in \Z)$ where the field elements commute with the operators $\Gdiff_i$. The statement that $A$ and $B$ are inverse,
is equivalent to
\begin{align}
   \tripair{\pi_w}{GH}{\Gdiff_v} &= \delta_{v,w} \\
   \tripair{\pi_w}{HG}{\Gdiff_v} &= \delta_{v,w}
\end{align}
for all $v \leleft{J} w$. The notation $\tripair{f}{g}{h}$ means apply $g$ to $h$ and take the coefficient of $f$.

Working on all cosets at once,
this is further equivalent to 
\begin{align}
\label{E:pairGH}
	\tripair{\Gdiff_w}{GH}{\Gdiff_\id} &= \delta_{\id,w} \\
\label{E:pairHG}
	\tripair{\Gdiff_w}{HG}{\Gdiff_\id} &= \delta_{\id,w}
\end{align}
for all $w\in S_J$ since
\begin{align*}
\tripair{\Gdiff_w}{GH}{\Gdiff_v}= \sum_{\substack{u \\ u*v=w}} 
\tripair{\Gdiff_u}{GH}{\Gdiff_\id}
\end{align*}
and similarly for $HG$ instead of $GH$.

Thus it suffices to prove \eqref{E:pairHG}, as \eqref{E:pairGH} holds by formal properties of inverses.
But by Proposition \ref{P:G double to single} equation \eqref{E:pairHG} is equivalent to
\begin{align}
\label{E:cancellation}
	\delta_{\id,w} &= \G_w(a;a).
\end{align}
Equation \eqref{E:cancellation}  follows from the fact that $\Psi^{-1}(\G_w(va;a)) = [{\mathcal O}^w]|_v$ is the localization at $v$, of the equivariant class of the structure sheaf $\mathcal{O}^w$ of the opposite Schubert variety $\overline{B_-wB/B}$ in $K_T^*(\Fl_n)$.

We state a ``right-handed" version of Proposition~\ref{P:left Groth inversion}.  With $J$ as above let $v,w\in S_{\Z}$ be such that $vS_J=w S_J$.
Say $v \leright{J} w$ if there is a $u\in S_J$ such that $v*u=w$.

\begin{prop} \label{P:right Groth inversion}
	Let $W'$ be a $\leright{J}$ interval in $S_{\Z}$ and let
	$\{f_w\mid w\in W'\}$ and $\{g_w\mid w\in W'\}$ be families of elements. Then
	\begin{align}\label{E:right change}
		g_w &= \sum_{v*u=w} (-1)^{\ell(u)+\ell(v)-\ell(w)} \G_u(a) f_v \intertext{holds if and only if}  \label{E:right inverse change}
		f_w &= \sum_{v*u=w} (-1)^{\ell(u)+\ell(v)-\ell(w)} \G_{u^{-1}}(\ominus a) g_v
	\end{align}
	does.
\end{prop}

\section{Some computations}
\label{S:computations}

Recall that $\xm = (x_0,x_{-1},x_{-2},\dotsc)$.  We write $X_i$ for $(x_1,x_2,\dotsc,x_i)$.

\subsection{Back stable Grothendieck polynomials}\
We compute $\bG_{s_0}$ from the definition.
\begin{align*}
\begin{split}
	\G_{s_1} &= e_1[X_1] \\
	\G_{s_2} &= e_1[X_2]-e_2[X_2] \\ 
	\G_{s_3} &= e_1[X_3]-e_2[X_3]+e_3[X_3] \\
	\G_{s_n} 
	&= e_1[X_n] - e_2[X_n]+e_3[X_n]-\dotsm+(-1)^{n-1}e_n[X_n]
\end{split}
\end{align*}
Therefore,
$$
	\bG_{s_0} = G_{s_0} 
	= \lim_{n\to-\infty} \G_{s_n}(x_{1-n},\dotsc,x_{-1},x_0) 
	= e_1 - e_2 + e_3 - \dotsm
	= 1 - \Omega[-x].
$$

We compare Proposition~\ref{P:G back stable shift} and Theorem~\ref{T:coproduct} for $\bG_{s_1}, \bG_{s_2}, \bG_{s_{-1}}$.
\begin{align*}
	\bG_{s_1} &= \gamma(\bG_{s_0}) \\
	&= e_1[\xm + x_1] - e_2[\xm+x_1] + e_3[\xm+x_1] - \dotsm \\
	&= (x_1 + e_1) - (x_1e_1+e_2) + (x_1e_2+e_3)-\dotsm \\
	&= x_1 - x_1(e_1-e_2+e_3-\dotsm) +(e_1-e_2+e_3-\dotsm) \\
	&= \G_{s_1} - G_{s_1} \G_{s_1} + G_{s_1} \\
	\bG_{s_2} &= \gamma^2(\bG_{s_0}) \\
	&= e_1[\xm + X_2] - e_2[\xm+X_2] + e_3[\xm+X_2] - \dotsm \\
	&= (e_1[X_2]+e_1) - (e_2[X_2]+e_1[X_2]e_1+e_2) + (e_2[X_2]e_1+e_1[X_2]e_2+e_3)-\dotsm \\
	&= (e_1[X_2]-e_2[X_2]) + (1-e_1[X_2]+e_2[X_2]) e_1 + (-1+e_1[X_2]-e_2[X_2]) e_2 \\
	&+(1-e_1[X_2]+e_2[X_2])e_3+\dotsm \\
	&= \G_{s_2} + G_{s_2} - G_{s_2} \G_{s_2} \\
	\bG_{s_{-1}} &= \gamma^{-1}(\bG_{s_0}) \\
	&= e_1[\xm-x_0] - e_2[\xm-x_0] + e_3[\xm-x_0] -\dotsm \\
	&= (e_1-x_0) - (e_2-x_0 e_1+x_0^2) + (e_3-x_0e_2+x_0^2e_1-x_0^3)-\dotsm\\
	&= G_{s_{-1}} + \dfrac{-x_0}{1-x_0} - \dfrac{-x_0}{1-x_0} G_{s_{-1}} \\
	&= G_{s_{-1}} + \G_{s_{-1}} - G_{s_{-1}} \G_{s_{-1}}.
\end{align*}
Here,
\begin{align*}
	\G_{s_{-1}} = \omega(G_{s_1}) = \omega(x_1) = \dfrac{-x_0}{1-x_0}.
\end{align*}
In particular for $k > 0$ we have
\begin{align*}
  \bG_{s_k} &= G_1 \oplus x_1 \oplus x_2 \oplus \dotsm \oplus x_k \\
  \bG_{s_{-k}} &= G_1 \ominus x_0 \ominus x_1 \ominus \dotsm \ominus x_{1-k}.
\end{align*}

We compute $\bG_{s_{-1}s_0}$ from the limit definition:
\begin{align*}
	\G_{s_1s_2} &= e_2[X_2] \\
	\G_{s_2s_3} &= e_2[X_3] - 2 e_3[X_3] \\
	\G_{s_3s_4} &= e_2[X_4] - 2 e_3[X_4] + 3 e_4[X_4]
\end{align*}
Thus,
\begin{align*}
	\bG_{s_{-1}s_0} &= G_{s_{-1}s_0} = G_{11} = e_2 - 2 e_3 + 3 e_4 -\dotsm.
\end{align*}
Using Proposition~\ref{P:G back stable shift}, we have
\begin{align*}
	\bG_{s_0s_1} &= e_2[\xm+x_1] - 2 e_3[\xm+x_1] + 3e_4[\xm+x_1] \\
	&= (e_2+x_1 e_1) - 2(e_3+x_1 e_2) + 3(e_4+x_1e_3) - \dotsm \\
	&= x_1(e_1-2e_2+3e_3-\dotsm)+ (e_2-2e_3+3e_4-\dotsm) \\
	&= x_1 G_{s_0} - x_1 (e_2-2e_3+\dotsm) + G_{s_0s_1} \\
	&= G_{s_0}\G_{s_1} - G_{s_0s_1} \G_{s_1} + G_{s_0s_1},
\end{align*}
agreeing with Theorem~\ref{T:coproduct}.

With $s_\la$ denoting a Schur polynomial, we compute $\bG_{s_1s_0}$ from the limit definition:
\begin{align*}
	\G_{s_2s_1} &= s_2[X_1] \\
	\G_{s_3s_2} &= s_2[X_2]-s_{21}[X_2] \\
	\G_{s_4s_3} &= s_3[X_3]-s_{21}[X_3]+s_{211}[X_3] \\
	\bG_{s_1s_0} &= G_{s_1s_0} = s_2 - s_{21} + s_{211} - s_{2111} + \dotsm
\end{align*}
Therefore, using Proposition~\ref{P:G back stable shift}, we have
\begin{align*}
	\bG_{s_0s_{-1}} &= \gamma^{-1}(\bG_{s_1s_0}) \\
	&= s_2[\xm-x_0] - s_{21}[\xm-x_0] + s_{211}[\xm-x_0]-s_{2111}[\xm-x_0] +\dotsm \\
	&= (s_2-x_0 s_1) - (s_{21}-x_0s_2-x_0s_{11}+x_0^2 s_1) \\
	&+ 
	(s_{211}-x_0 s_{21}-x_0 s_{111} +x_0^2 s_2 +x_0^2 s_{11} -x_0^3 s_1 ) \\
	&- (s_{2111}-x_0 s_{211}-x_0s_{1111}+x_0^2 s_{21}+x_0^2 s_{111}-x_0^3s_2 -x_0^3 s_{11}+x_0^4 s_1) + \dotsm \\
	&= G_{s_0s_{-1}} + G_{s_0} \G_{s_{-1}} - G_{s_0s_{-1}} \G_{s_{-1}}
\end{align*}

\subsection{$K$-Stanley polynomials}\


\noindent
We give some formulae for various $G_\la$.  See also Proposition~\ref{P:G det} and Proposition~\ref{P:G one row column}.
\begin{align*}
	G_r &= s_r - s_{r,1} + s_{r,1,1} - s_{r,1,1,1} + \dotsm \\
	G_{1^r} &= e_r - \binom{r}{r-1} e_{r+1} + \binom{r+1}{r-1} e_{r+2} - \dotsm \\
	G_{21} &= s_{21} - s_{22} -2 s_{211} + 2 s_{221} - s_{222} + 3 s_{2111} - 3s_{2211} +2 s_{2221} - s_{2222} + \dotsm \\
	G_{31} &= s_{31} - s_{32} - 2s_{311} + 2 s_{321} - s_{322}+\dotsm \\
	G_{211} &= s_{211} - s_{221} + s_{222} - 3 s_{2111} + 3 s_{2211} - 3s_{2221}+2s_{2222} + \dotsm 
\end{align*}
We give some products of $K$-Stanley functions; see Proposition~\ref{P:stable K stanley expansion}.
\begin{align*}
	G_1 G_1 &= G_2 + G_{11} - G_{21} \\
	G_2 G_1 &= G_3 + G_{21} - G_{31} \\
	G_{11} G_1 &= G_{21}+G_{111} - G_{211}.
\end{align*}

 
 
 

\end{document}